\documentclass[11pt, british]{article}

\usepackage{declare-style-article}
\usepackage{declare-maths-article}
\usepackage{declare-text-article}
\usepackage{declare-theorems-article}
\hypersetup{
	pdfauthor 
		= {Nikita Nikolaev},
	pdftitle 
		= {Geometry and Resurgence of WKB Solutions of Schrödinger Equations},
	pdfsubject
		= {},
	pdfcreator
		= {},
	pdfproducer
		= {},
	pdfkeywords
		= {},
} 

\DeclareSymbolFont{bbold}{U}{bbold}{m}{n}
\DeclareSymbolFontAlphabet{\mathbbold}{bbold}

\usepackage{appendix}
\usepackage[font={footnotesize}]{caption}
\newcommand{\pto}[1]{{\scriptscriptstyle(#1)}}

\newcommand{\inv}{{\scriptscriptstyle{\text{--}1}}}

\newcommand{\sqroot}{{\scriptscriptstyle{1\mskip-2mu{/}\mskip-1mu{2}}}}
\newcommand{\invsqroot}{{\scriptscriptstyle{\text{--}1\mskip-2mu{/}\mskip-1mu{2}}}}
\newcommand{\cubesqroot}{{\scriptscriptstyle{3\mskip-2mu{/}\mskip-1mu{2}}}}

\fancypagestyle{frontpage}{

\cfoot{}
\lfoot{}
}

\usepackage{titletoc}
\usepackage{subcaption}
\usepackage{changepage}

\newcommand{\MSCSubjectCode}[1]{\href{https://zbmath.org/classification/?q=cc\%3A#1}{#1}}

\newcommand{\sfSigma}{\sf{\Sigma}}
\newcommand{\sfDelta}{\sf{\Delta}}

\usetikzlibrary{decorations.markings}
\usetikzlibrary{calc}
\usetikzlibrary{arrows.meta}
\tikzset{-latex-/.style={decoration={
  markings,
  mark=at position #1 with {\arrow{latex}}},postaction={decorate}}}

\makeatletter
\renewcommand\@dotsep{200}
\makeatother

\makeatletter
\renewcommand{\dddot}[1]{%
  {\mathop{\kern\z@#1}\limits^{\vbox to-1.4\ex@{\kern-\tw@\ex@
   \hbox{\footnotesize ...}\vss}}}}
\makeatother

\begin{document}


\newgeometry{top=1cm, bottom=0.5cm, left=3.5cm, right=3.5cm}

\title{Geometry and Resurgence\\of WKB Solutions of Schrödinger Equations}

\author{Nikita Nikolaev}

\affil{\small School of Mathematics, University of Birmingham, United Kingdom\vspace{-10pt}}

\date{22 October 2024}

\vspace{-0.5cm}

\maketitle
\thispagestyle{frontpage}

\begin{abstract}
\noindent
We prove that formal WKB solutions of Schrödinger equations on Riemann surfaces are resurgent.
Specifically, they are Borel summable in almost all directions and their Borel transforms admit endless analytic continuation away from a discrete subset of singularities.
Our approach is purely geometric, relying on understanding the global geometry of complex flows of meromorphic vector fields using techniques from holomorphic Lie groupoids and the geometry of spectral curves.
This framework provides a fully geometric description of the Borel plane, Borel singularities, and the Stokes rays.
In doing so, we introduce a geometric perspective on resurgence theory.
\end{abstract}

{\footnotesize

\vspace{-0.25cm}

\textbf{2020 MSC:}
\MSCSubjectCode{34M60} (primary),	
\MSCSubjectCode{34M40}, 
\MSCSubjectCode{34M30},	
\MSCSubjectCode{34M35},	
\MSCSubjectCode{34E20},	
\MSCSubjectCode{30F20},	
\MSCSubjectCode{32M25},	
\MSCSubjectCode{32S65}	
}

\vspace{-0.25cm}

{\begin{spacing}{0.9}
\small
\setcounter{tocdepth}{2}
\tableofcontents
\end{spacing}

\restoregeometry
\setcounter{section}{-1}
\section{Introduction}
\counterwithout{equation}{section}
\enlargethispage{20pt}

Consider the Schrödinger equation
\vspace{-5pt}
\begin{eqntag}
\label{240711180437}
	\hbar^2 \del_x^2 \Psi (x, \hbar) = \Q (x, \hbar) \Psi (x, \hbar)
\vspace{-5pt}
\end{eqntag}
where $x$ is a local complex variable on a compact Riemann surface $X$,~ $\hbar \in \CC$ is a small complex perturbation parameter, and $\Q$ is a polynomial in $\hbar$ with meromorphic coefficients.
It is well-known that away from turning points, which are the zeros and simple poles of the quadratic differential $\phi_0 \coleq \Q_0 (x) \d{x}^2$, this equation has a basis of \textit{formal WKB solutions} of the form
\vspace{-5pt}
\begin{eqntag}\label{230516183430}
	\hat{\Psi} (x, \hbar)
		= e^{- \S (x) / \hbar} \hat{\A} (x, \hbar)
\qqtext{where}
	\hat{\A} (x, \hbar) \coleq \sum_{k=0}^\infty \A_k (x) \hbar^k
\fullstop{,}
\vspace{-5pt}
\end{eqntag}
enumerated by the two sheets of the \textit{spectral curve} $\sfSigma$ which is the Riemann surface of the square root of $\phi_0$.
The function $\S$ is the integral of the Liouville one-form $\lambda$ along local paths on one of two sheets of $\sfSigma$.

These formal solutions are extremely convenient not least because they are explicitly computable using little more than simple algebraic manipulations that involve no differential equations at all.
However, they are not true analytic solutions because the power series $\hat{\A}$ is almost always divergent.
A formal WKB solution $\hat{\Psi}$ can always be lifted to a true analytic solution $\Psi$ which is asymptotic to $\hat{\Psi}$ as $\hbar \to 0$, but such lifts are highly non-unique and typically non-constructive.
The aim of \textit{exact perturbation theory} is to construct such asymptotic lifts in a canonical and explicit way using the \textit{Borel resummation} method.
This lifts the divergent series $\hat{\A}$ to a holomorphic function $\A_{\alpha}$, defined in a halfplane sector of the $\hbar$-plane bisected by a ray $\alpha$, which is given by a Laplace integral, yielding so-called \textit{exact WKB solutions}:
\vspace{-5pt}
\begin{eqntag}
\label{240801130628}
	\Psi_{\alpha} (x, \hbar) 
	= e^{-\S (x) / \hbar} \A_{\alpha} (x, \hbar)
	= e^{-\S (x) / \hbar} \int_0^{\infty e^{i \alpha}} \hat{\Phi} (x, t) e^{-t/\hbar} \d{t}
\fullstop
\vspace{-5pt}
\end{eqntag}
The \textit{Borel transform} $\hat{\Phi}$ is an explicit power series in $t$ easily obtained from $\hat{\A}$.
Whether or not this Laplace integral is well-defined is really a question about the analytic continuation of $\hat{\Phi}$ in the variable $t$.
Directions $\alpha$ such that the integration contour encounters a singularity of $\hat{\Phi}$ are called \textit{Stokes rays}.
As $\alpha$ varies across a Stokes ray, exact WKB solutions exhibit a beautiful Stokes phenomenon which is often called the \textit{WKB connection formula} in the physics literature.
The significance of this Stokes phenomenon in mathematics and physics can hardly be overstated.

\textit{Resurgence} is the study of the global structure of singularities of the Borel transform, and can be broadly summed up as global complex analysis for divergent series.
In essence, \textit{resurgent series} are a class of divergent series which are \textit{Borel summable} in almost all directions (i.e., the Laplace integral in \eqref{240801130628} is well-defined for almost all $\alpha$) and whose Stokes phenomenon across Stokes rays can be described in terms of the Borel resummation of new divergent power series extracted from the singularities of the Borel transform.
In particular, the Borel transform must admit an \textit{endless analytic continuation}; i.e., it can be analytically continued along any path that avoids a discrete set of \textit{Borel singularities}.
Surprisingly, this class of divergent series is both vast --- capturing many phenomena of contemporary interest --- and strikingly special --- exhibiting exceptional properties not typical of general analytic functions.

\paragraph{Main Results.}
It has been a major open problem since at least the early 1980s to establish if formal WKB solutions are resurgent in $\hbar$ and, if so, to describe the location and nature of their Borel singularities.
This property is sometimes called \textit{quantum resurgence} of the Schrödinger equation.
In this article, we give a complete geometric solution to this problem which can be roughly formulated as follows.

\vspace{-2.5pt}

Consider a Schrödinger equation on an arbitrary compact Riemann surface $X$ with potential $\Q$ that has poles along a (possibly empty) effective divisor $D$, and assume that all zeros of the associated quadratic differential $\phi_0$ are simple.
Let $\pi : \sfSigma \to X$ be the spectral double cover with canonical involution $\iota : \sfSigma \to \sfSigma$, and let $\sfDelta \subset \sfSigma$ be the preimage of $D$ (excluding the simple poles of $\phi_0$).
For any point $x \in X$ that is neither a pole nor a turning point, pick one of its preimages $p \in \sfSigma$, and let $\hat{\Psi} (x, \hbar)$ be a formal WKB solution corresponding to the sheet of $\sfSigma$ containing $p$.

\vspace{-2.5pt}

Consider the Riemann surface $\tilde{\sfSigma}_{\sfDelta, p}$, called the \textit{Borel surface}, defined as the universal cover of the punctured spectral curve $\sfSigma_\sfDelta \coleq \sfSigma \smallsetminus \sfDelta$ based at $p$.
We view $\tilde{\sfSigma}_{\sfDelta, p}$ as the space of fixed-endpoint homotopy classes of paths on $\sfSigma_\sfDelta$ starting at $p$.
As such, $\tilde{\sfSigma}_{\sfDelta, p}$ has a canonical basepoint given by the constant path $1_p$, and there is a canonical map to the $t$-plane, called the \textit{central charge}:
\vspace{-7.5pt}
\begin{eqn}
	\Z_p : \tilde{\sfSigma}_{\sfDelta, p} \too \CC
\qqtext{sending}
	\gamma \mapstoo t = \int_\gamma (\lambda - \iota^\ast \lambda)
\fullstop
\vspace{-5pt}
\end{eqn}
This is an infinitely-sheeted branched covering map with only algebraic branch points.
Its ramification locus $\Gamma_p \subset \tilde{\sfSigma}_{\sfDelta, p}$ is a countable discrete subset consisting of fixed-endpoint homotopy classes of \textit{critical paths} which are paths that start at $p$ and terminate at ramification points $R_\sfDelta \subset \sfSigma_\sfDelta$; see \autoref{241021102954}.
As a set, $\Gamma_p$ is isomorphic to $n$ copies of $\pi_1 (\sfSigma_\sfDelta, p)$ where $n$ is the total number of turning points.
There is a distinguished subset $\Gamma_p^0 \subset \Gamma_p$ of elements representable by \textit{critical trajectories}; i.e., paths along which the central charge $\Z_p$ has a constant phase.
\textit{Saddle trajectories} on $\sfSigma_\sfDelta$ are trajectories which begin and end at ramification points, whilst \textit{divergent trajectories} are those whose closure has nonempty interior.

\vspace{-2.5pt}

\enlargethispage{10pt}

{\bfseries \autoref{240801161941}} (summary).
{\itshape
The formal WKB solution $\hat{\Psi}$ is resurgent, uniformly in $x$ in an appropriate sense.
More specifically:

\vspace{-2pt}

\textup{\bfseries(1)}
The Borel transform $\hat{\Phi} (x,t)$ of $\hat{\A} (x,\hbar)$ is a convergent power series in $t$.
It has a canonical lift to a holomorphic germ at the basepoint $1_p$ of the Borel surface $\tilde{\sfSigma}_{\sfDelta, p}$ which has endless analytic continuation away from the discrete subset $\Gamma_p \subset \tilde{\sfSigma}_{\sfDelta, p}$.

\vspace{-2pt}

\textup{\bfseries(2)}
The image of $\Gamma_p^0 \subset \Gamma_p$ under the central charge $\Z_p$ is a discrete subset $\Xi_p^0 \subset \CC$ whose points determine a countable set of directions called \textup{Stokes rays}.
The complement of all radial cuts emerging from $\Xi_p^0$ is the largest open star-shaped domain $\EE_p \subset \CC$ around the origin which lifts biholomorphically to a neighbourhood of $1_p$ in $\tilde{\sfSigma}_{\sfDelta, p}$.

\vspace{-2pt}

\textup{\bfseries(3)}
The Laplace transform of $\hat{\Phi} (x,t)$ is well-defined along every infinite ray $e^{i \alpha} \RR_+$ contained in $\EE_p$, provided that $\alpha$ is not the phase of a divergent trajectory passing through $p$.
Thus, $\hat{\Psi} (x,\hbar)$ is Borel summable in every such direction $\alpha$.

\vspace{-2pt}

\textup{\bfseries(4)}
The discontinuous jump in the Borel resummation of $\hat{\Psi} (x,\hbar)$ across a Stokes ray $\alpha$ is the countable sum of exponentially small, Borel resummed contributions coming from the critical trajectory $\gamma \in \Gamma^0_p$ with phase $\alpha$ and all finite chains of saddle trajectories with phase $\alpha$ composable with $\gamma$, provided that $\gamma$ is not in the closure of a divergent trajectory with phase $\alpha$.
}

\pagebreak

\begin{figure}[p]
\centering
\begin{adjustwidth}{-1cm}{-1cm}
\begin{subfigure}{0.65\textwidth}
\includegraphics[width=\textwidth]{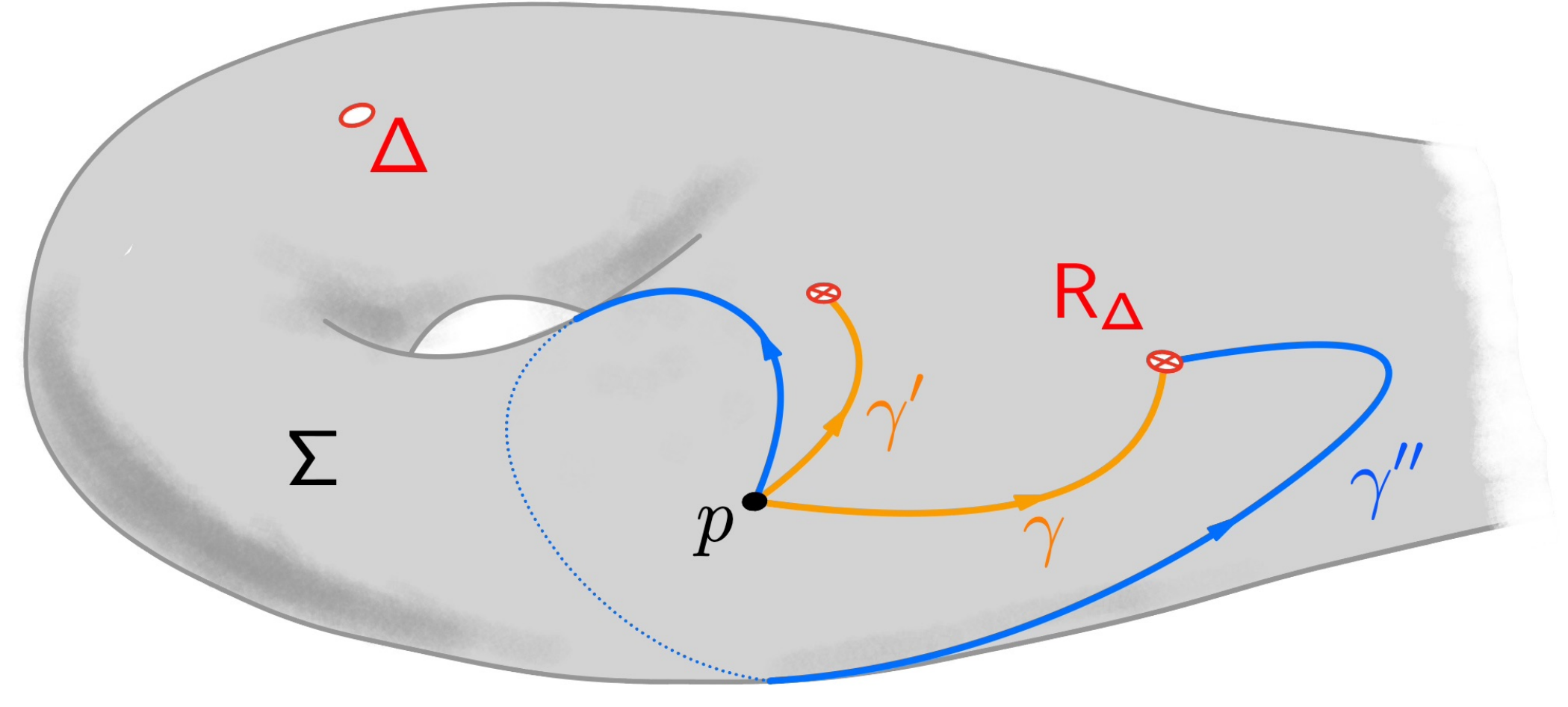}
\caption{}
\label{241022113833}
\end{subfigure}
\quad
\begin{subfigure}{0.4\textwidth}
\includegraphics[width=\textwidth,trim={0 2.5cm 0 0},clip]{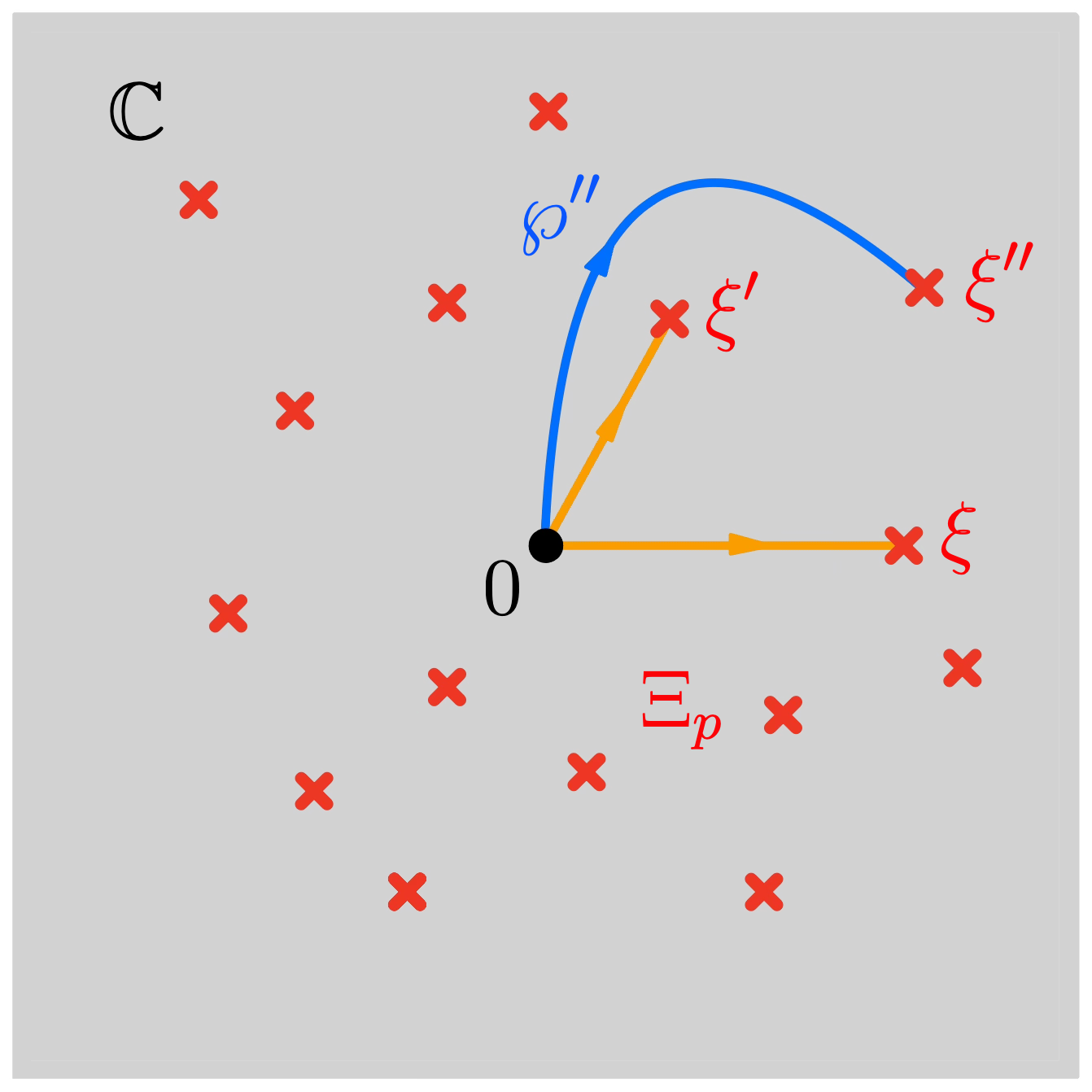}
\caption{}
\label{241022113847}
\end{subfigure}
\smallskip\\
\begin{subfigure}{0.65\textwidth}
\includegraphics[width=\textwidth]{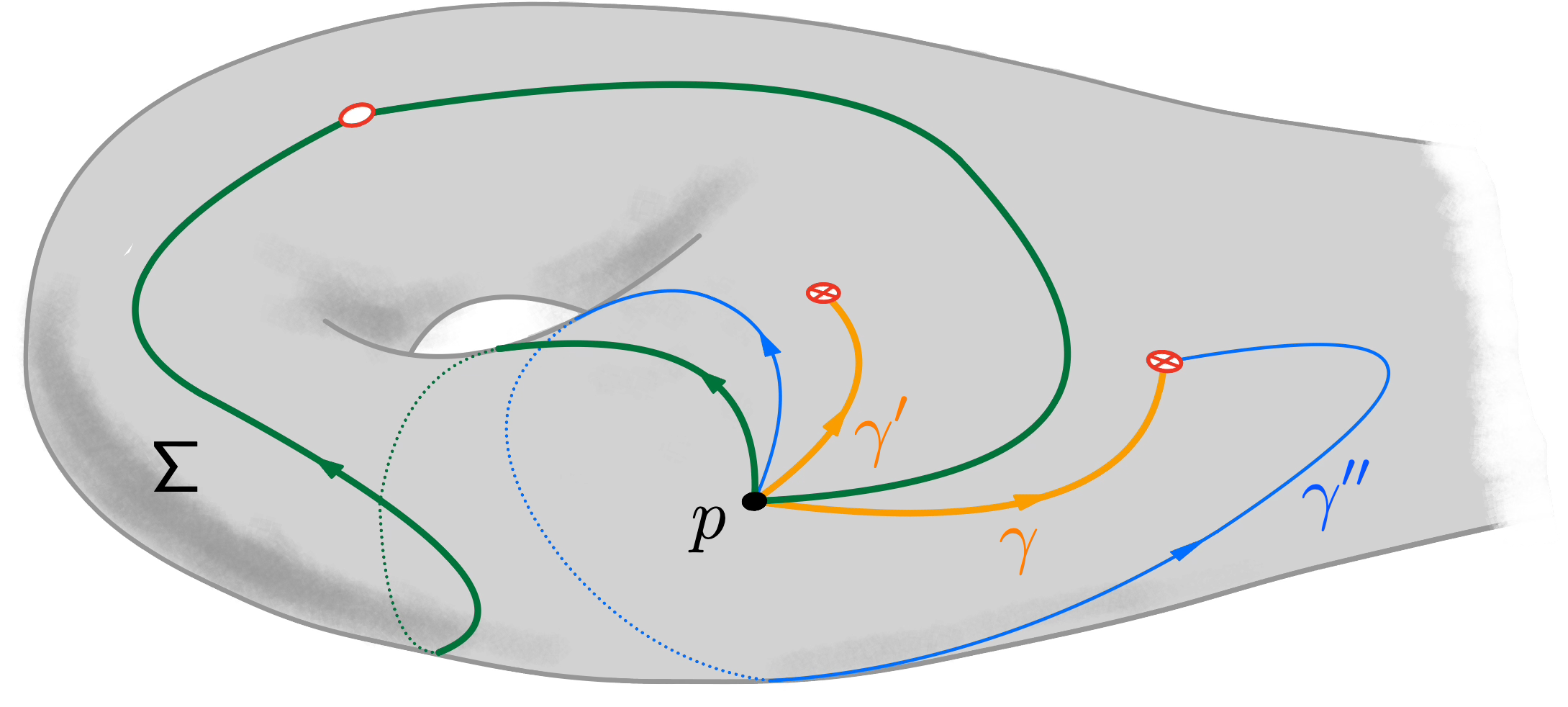}
\vspace{5pt}
\caption{}
\label{241022114907}
\end{subfigure}
\quad
\begin{subfigure}{0.4\textwidth}
\includegraphics[width=\textwidth]{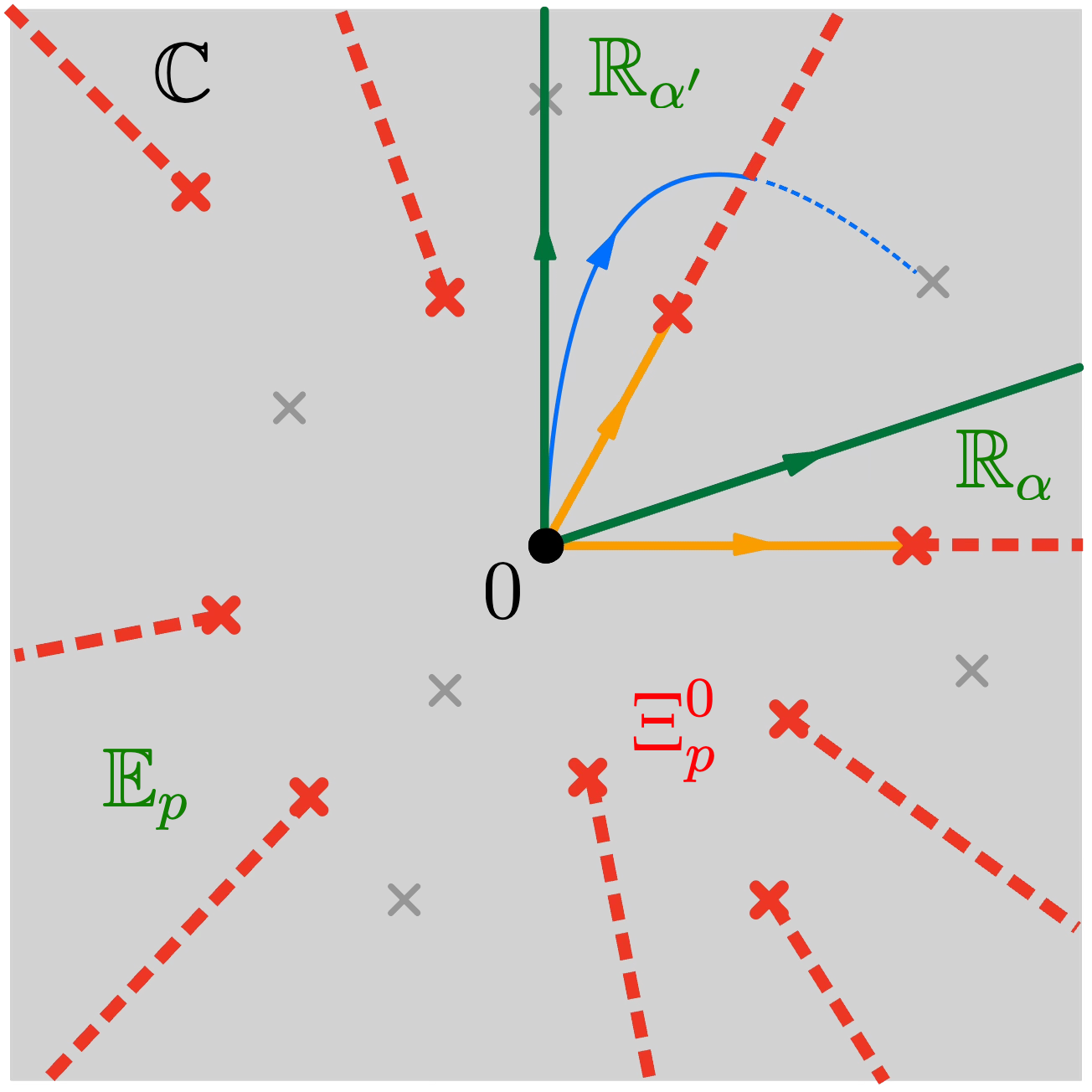}
\caption{}
\label{241022115526}
\end{subfigure}
\smallskip\\
\begin{subfigure}{0.65\textwidth}
\includegraphics[width=\textwidth]{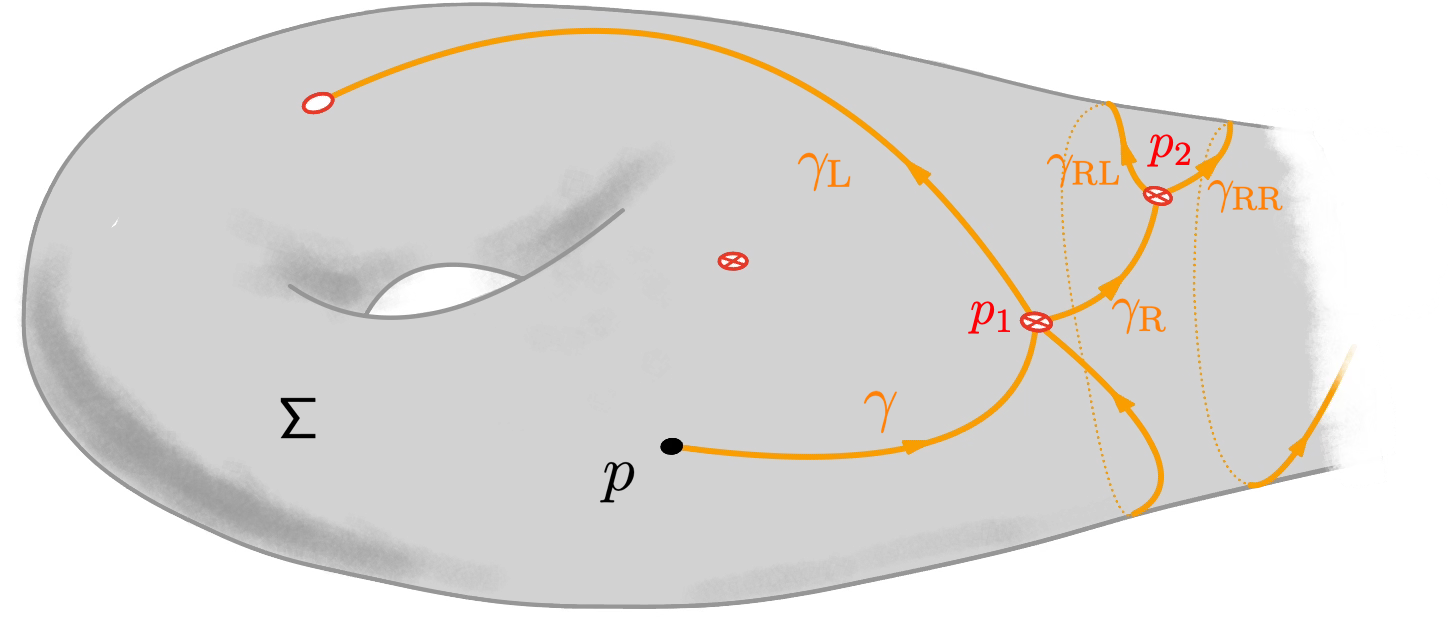}
\caption{}
\label{241022125828}
\end{subfigure}
\quad
\begin{subfigure}{0.4\textwidth}
\includegraphics[width=\textwidth]{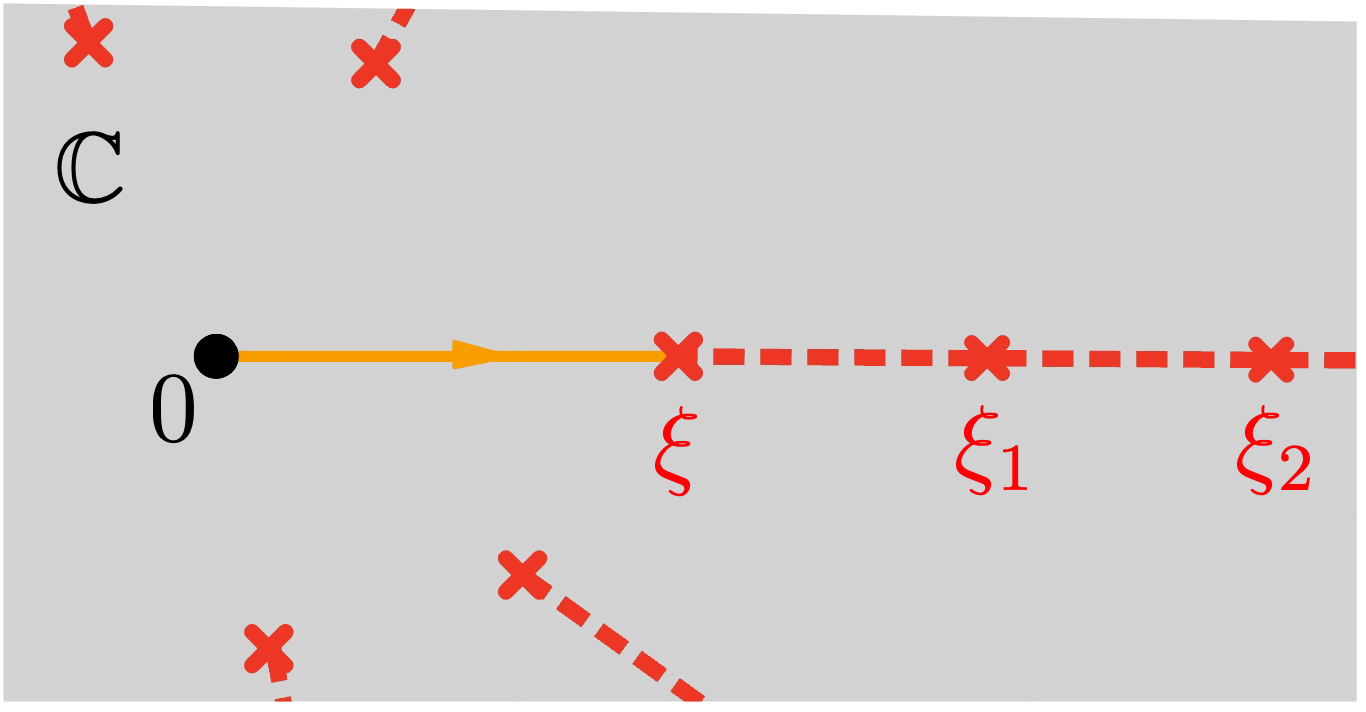}
\vspace{5pt}
\caption{}
\label{241022125902}
\end{subfigure}
\end{adjustwidth}
\caption{All three paths $\gamma,\gamma',\gamma''$ pictured in (a) are critical paths on the punctured spectral curve $\sfSigma_\sfDelta$: they start at $p$ and terminate at ramification points depicted by circled crosses, so they all determine distinct elements of $\Gamma_p$.
The concatenated path $\gamma^\inv \circ \gamma''$ is a loop in $\sfSigma_\sfDelta$ based at $p$, hence an element of the fundamental group $\pi_1 (\sfSigma_\sfDelta, p)$, yet $\gamma$ and $\gamma'$ are not related by any loop on $\sfSigma_\sfDelta$.
This exhibits the fact that $\Gamma_p$ is a disjoint union of several $\pi_1 (\sfSigma_\sfDelta, p)$-torsors.
\smallskip\\
The paths $\gamma$ and $\gamma'$ are critical trajectories because their projections to $\CC$ via the central charge $\Z_p$ pictured in (b) are straight line segments $[0,\xi]$ and $[0,\xi']$, so $\gamma, \gamma' \in \Gamma_p^0$.
On the other hand, the path $\gamma''$ is not a critical trajectory because its projection $\wp''$ is not straight.
Moreover, $\gamma' \notin \Gamma_p^0$ because $\wp''$ is not homotopic to the straight line segment $[0,\xi'']$ due to the obstructing singular point $\xi'$.
\smallskip\\
The discrete subset $\Xi_p^0 = \Z_p (\Gamma_p^0)$ is depicted in (d) by red crosses, and there are radial cuts emerging from each point.
Their complement $\EE_p$ contains infinite rays, such as $\RR_\alpha \coleq e^{i \alpha} \RR_+$ and $\RR_{\alpha'}$, which have a (locally) biholomorphic lift to infinitely long trajectories on the spectral curve $\sfSigma_\sfDelta$ pictured in (c).
In the depicted situation, these trajectories both fall into a pole in $\sfDelta$.
\smallskip\\
Notice that the path $\wp''$ from $0$ to $\xi''$ crosses the radial cut emerging from $\xi'$.
Thus, the singularity represented by the path $\gamma''$ is not visible from the point $p$ in the sense that it cannot be reached along a trajectory, which is why the singular value at $\xi''$ (as well as many other singular values pictured in (b)) has been greyed out and there is no radial cut emerging from it.
\smallskip\\
(Caption continues on the next page.) 
}
\label{241021102954}
\end{figure}

\addtocounter{figure}{-1}
\begin{figure}[t]
\centering
\caption{(Caption continued from previous page.)
The infinite rays that bypass the singularity $\xi$ slightly on the left and slightly on the right yield a pair of Borel resummed WKB solutions.
To describe the discontinuous jump across the Stokes ray in the direction of $\xi$, we must investigate the presence of relevant saddle trajectories as depicted in (e).
The corresponding critical trajectory $\gamma$ hits the ramification point $p_1$ which can be circumvented on the left or on the right.
Continuing on the left, we follow the trajectory $\gamma_\rm{L}$ which is regular because it falls into a pole in $\sfDelta$.
\smallskip\\
On the other hand, if we continue on the right, we follow the trajectory $\gamma_\rm{R}$ which is critical because it hits another ramification point $p_2$.
Thus, $\gamma_\rm{R}$ is a saddle trajectory, and the central charge of the concatenated path $\gamma_\rm{R} \circ \gamma$ is the singular value $\xi_1$ which lies on the radial cut emerging from $\xi$, as depicted in (f).
Notice that this singular value was not visible when we chose to circumvent $\xi$ on the left, but it is visible if we do so on the right.
In the depicted situation, these trajectories both fall into a pole in $\sfDelta$.
\smallskip\\
Upon reaching $p_2$, we must again make a choice of on what side to circumvent the singularity.
If we choose the right side, we follow the trajectory $\gamma_{\rm{RR}}$ which we imagine is regular because it continues  uninterrupted and, say, falls into another pole in $\sfDelta$.
Consequently, the trajectory $\gamma_\rm{RR}$ projects via $\Z_p$ to an straight halfline line in $\CC$ emanating from $\xi_1$.
On the other hand, if choose the left side, we follow the trajectory $\gamma_{\rm{RL}}$ which is another saddle trajectory.
The central charge of the concatenated path $\gamma_{\rm{RL}} \circ \gamma_{\rm{R}} \circ \gamma$ is the new singular value $\xi_2$ which also lies on the cut emerging from $\xi_1$.
Again, this singularity is only visible if we make the correct sequence of turns at each singular encountered singularity (in this case, ``first right, then left''), and it corresponds to a specific finite chain of composable saddle trajectories (in this case, $\gamma_{\rm{RL}} \circ \gamma_{\rm{R}}$).
\smallskip\\
The saddle trajectory $\gamma_{\rm{RL}}$ hits a ramification point, and we can repeat the analysis.
In this case, this ramification is actually again the point ramification point $p_1$; i.e., the path $\gamma_{\rm{RL}} \circ \gamma_{\rm{R}}$ is in fact a homology cycle of $\sfSigma_\sfSigma$.
So, in a situation like this, we can add an arbitrary multiple of this homology cycle and therefore discover countably-infinitely many new singular values located on the cut emerging from $\xi$, each a constant integer multiple translate of another.
}
\label{241022132559}
\end{figure}

Furthermore, if there are no divergent trajectories at all, then an even stronger form of resurgence holds; see \autoref{241007171114}.
In order to prove \autoref{240801161941}, we develop a coordinate-invariant formalism for Schrödinger equations and the WKB method.
In particular, we show in \autoref{240916114621} that, given a Schrödinger equation, there is a canonical formal $\hbar$-series-valued differential form $\hat{\Lambda} (\hbar)$ (which we call the \textit{formal WKB differential}) on the spectral curve $\sfSigma$ that extends the canonical one-form $\lambda$.
The formal WKB solutions are then canonically expressible, in a coordinate-independent manner, in terms of the formal WKB differential, and their resurgence properties follow from that of $\hat{\Lambda} (\hbar)$.
So, in \autoref{240918164353}, we formulate an analogous result asserting that $\hat{\Lambda} (\hbar)$ is a resurgent series of differential forms on $\sfSigma$.

\paragraph{Context.}
The resurgence property of WKB solutions for Schrödinger equations with polynomial potential (i.e., $\Q (x, \hbar) = \Q_0(x)$ is a polynomial in $x$) was first conjectured by Voros in \cite{MR729194,MR728983} and Écalle in \cite[p.40]{EcalleCinqApplications}; see also~\cite[Comment on p.32]{MR1704654}.
In fact, questions of an essentially equivalent nature go back even further to Leray \cite{MR103328} and Hamada \cite{MR0249815}.
Numerous attempts at this problem exist in the literature, notably \cite{MR1162634,MR1232828,MR1209700,MR1270144,MR1278057,MR1397029,MR1654701,MR2118034,MR3050812,MR3185464,MR3266130,MR4088739,MR4226390,MR4441156,MY210623112236}.
However, these approaches are either mathematically incomplete, offer only a limited understanding of the resurgent structure, or study extremely special examples where analytic WKB solutions can be expressed by a more or less explicit formula that can be used to analyse the singularities of their Borel transform.
Thus, to the best of our knowledge, a proof of this conjecture has not yet appeared.

Our main result not only resolves this conjecture but also offers a positive answer in a much broader class of problems.
These include Schrödinger equations with rational potential that may depend on $\hbar$ in a polynomial way, as well as more general Schrödinger equations with or without poles defined over arbitrary compact Riemann surfaces $X$.
We are also able to treat the traditionally difficult, non-generic situations involving simple poles of the quadratic differential or the presence of divergent trajectories.
Our solution gives a complete geometric description of the resurgent structure in full generality.

We have also discovered that the nature of resurgence in the WKB analysis of Schrödinger equations is inherently better behaved than in other problems.
Namely, although the Borel singular values in the $t$-plane may be dense, they actually arise as the projection of a discrete set along a branched covering map which --- although infinitely-sheeted --- has only algebraic singularities.
This type of near-algebraic resurgence (which we call \textit{log-algebraic}) was expected to hold for polynomial potentials, but the fact that general Schrödinger equations share this remarkable property has come as a surprise to us.

This paper demonstrates that the question of resurgence of divergent series is really a question of global and hence geometric nature.
However, most attempts to such questions in the literature rely on predominantly local methods (with notable exceptions such as Kontsevich-Soibelman \cite{MR4403710} or Aniceto-Crew \cite{aniceto2024algebraicstudyparametricstokes}).
We believe this may be the reason that resurgence of WKB solutions had remained unsolved for so long.
In contrast, our approach is rooted in the geometry of the spectral cover and uses techniques from holomorphic Lie groupoids.

In recent years, resurgence has attracted a lot of attention in both mathematics and physics and the literature is rapidly becoming vast, with a few recent highlights including \cite{MR4287192,1811.05376,MR4403710,kontsevich2024holomorphic,Alim2024,kuwagaki2020sheaf,bridgel2024resurgence}.
Our interest in this subject stems from the study of moduli spaces of meromorphic connections, drawing inspiration from concepts such as spectral networks developed by Gaiotto-Moore-Neitzke \cite{MR3115984,MR3003931}; in particular, our goal is to unify and generalise the exact WKB method and the abelianisation of meromorphic connections~\cite{1902.03384}.

\paragraph{Generalised global flows.}
\label{241022140529}
A key new tool that we introduce in \autoref{240725160524} is \textit{generalised global flows} of meromorphic vector fields on a Riemann surface $M$.
Such vector fields, of course, do not generate a global flow $M \times \CC \to M$.
However, we have noticed that their local flow $\Omega \subset M \times \CC \to M$ actually admit a canonical holomorphic extension to the fundamental groupoid $\Pi_1 (M)$ (or more generally a Stokes groupoid in the sense of Gualtieri-Li-Pym \cite{MR3808258}; see \autoref{240725205848}) where it coincides with the groupoid target map $\Pi_1 (M) \to M$, which is a holomorphic surjective submersion.
Flowouts along such vector fields can then be understood as propagation along the source fibres of the groupoid followed by the target projection.
Consequently, any meromorphic vector field generates a global flow in this generalised sense.
A concrete application of this idea --- which is a crucial step in the proof of our main result --- appears in \autoref{240708173915}.

\paragraph{Narrative outline of the proof.}
Our proof of the Main Result proceeds as follows.
Starting with the WKB ansatz, we focus on the corresponding Riccati equation:
\begin{eqntag}
\label{240801205513}
	\Psi (x, \hbar) = \exp \left( - \frac{1}{\hbar} \int_{x_0}^x \Y (x', \hbar) \d{x'} \right)
\quad\implies\quad
	\hbar \del_x \Y = \Y^2 - \Q
\fullstop
\end{eqntag}
This equation depends on the local coordinate $x$ and the integrand $\Y \d{x}$ does \textit{not} determine a global differential form.
The first major observation in this paper is that it is possible to transform this local Riccati equation into another Riccati equation which makes global invariant sense on the spectral double cover $\pi : \sfSigma \to X$.
It is obtained by doing a change of the unknown variable that removes the leading- and sub-leading-order parts of $\Y$ in $\hbar$, resulting in an equation of the form
\begin{eqntag}
\label{240801211624}
	\hbar \V (f) - f = \hbar \big( f^2 + w f + \W \big)
\fullstop{,}
\end{eqntag}
where $f$ is the new globally well-defined unknown function on $\sfSigma$, the coefficients $w, \W$ are global meromorphic functions on $\sfSigma$, and $\V$ is the unique meromorphic vector field on $\sfSigma$ dual to the meromorphic one-form $\sigma = \lambda - \iota^\ast \lambda$ in the sense that $\sigma (\V) = 1$.
The poles and zeros of $\V$ are respectively the zeros and poles of $\sigma$.

Our second observation is that the meromorphic vector field $\V$ generates a global complex flow but only in the generalised sense, as described in \autoref{241022140529}, using the fundamental groupoid $\Pi_1 (\sfSigma_\sfDelta)$ of the punctured spectral curve $\sfSigma_\sfDelta \coleq \sfSigma \smallsetminus \sfDelta$ where $\sfDelta$ is the polar locus of $\sigma$.
From this point of view, the differential form $\sigma$ measures this vector field's complex flow time by integrating to a groupoid $1$-cocycle
\begin{eqntag}
\label{240801220010}
	\Z : \Pi_1 (\sfSigma_\sfDelta) \to \CC
\qqtext{sending}
	\gamma \mapsto t = \int_\gamma \sigma
\fullstop
\end{eqntag}
If we fix the starting point of integration $p \in \sfSigma$ and consider the groupoid source fibre $\tilde{\sfSigma}_{\sfDelta,p} = \rm{s}^\inv (p)$, then $\Z$ restricts to a covering map $\Z_p : \tilde{\sfSigma}_{\sfDelta,p} \to \CC$ which is ramified at every\footnote{For simplicity of exposition, let us assume in this subsection that $\phi_0$ has no simple poles.
In the general case, $\Z_p$ is actually unramified at those $\gamma$ corresponding to preimages of simple poles.} point represented by a path $\gamma$ that starts at $p$ and ends at a ramification point in $R_\sfDelta \subset \sfSigma_\sfDelta$.
Any path $\gamma \in \Pi_1 (\sfSigma_\sfDelta)$ which ends at a ramification point can be composed with its reversed reflection $\check{\gamma} = \iota(\gamma)^\inv$ on the other sheet of $\sfSigma$, where $\iota : \sfSigma \to \sfSigma$ is the involution.
The result is a path $\delta \coleq \iota(\gamma)^\inv \circ \gamma$ which exchanges sheets of the spectral cover and hence projects down to a loop on $X$ that goes around a turning point.
Since $\sigma = \lambda - \iota^\ast \lambda$, we find that $\Z (\gamma) = \int_\delta \lambda = \S (\delta)$; i.e., the branch locus of $\Z_p$ in $\CC$ consists of periods with respect to $\lambda$ of all the paths on $\sfSigma_\sfDelta$ that project down to loops that go around turning points in $X$.

The Riccati equation \eqref{240801211624} has a unique formal $\hbar$-power series solution $\hat{f}$ on $\sfSigma$, and our goal is to construct the analytic continuation $\phi$ of its Borel transform $\hat{\phi}$.
To do so, we analyse the equation satisfied by $\phi$ which is obtained by applying the Borel transform in $\hbar$ to \eqref{240801211624}, resulting in a first-order nonlinear differential equation with convolution of the following form:
\begin{eqntag}
\label{240809160212}
	\big( \V - \del_t \big) \phi = \phi \ast \phi + w \phi + \omega
\fullstop{,}
\end{eqntag}
where $\omega$ is the Borel transform of the coefficient $\W$ in \eqref{240801211624}.
The significance of this equation may be understood as follows.
Although the Riccati equation \eqref{240801211624} is globally well-defined and equivalent to the Schrödinger equation, on its own it is underdetermined.
This is an inherent feature of the WKB analysis where we seek solutions to an $\hbar$-family of equations subject to a strong boundary condition; namely, asymptotic regularity as $\hbar \to 0$.
The Borel transform allows us to incorporate this boundary condition into the problem in an effective way.
Thus, by applying the Laplace transform, global solutions of \eqref{240809160212} yield solutions of \eqref{240801211624} which necessarily satisfy the required boundary conditions.

It is possible to construct solutions of \eqref{240809160212} using an iteration scheme where at each step we flow the vector field $\V - \del_t$.
Our third major observation is that the vector field $\V - \del_t$ is nothing but the unique lift of $\V$ to a left-invariant vector field $\tilde{\V}$ on the groupoid $\Pi_1 (\sfSigma_\sfDelta)$ tangent to the target fibres.
Consequently, we can use $\tilde{\V}$ to flow initial data from the identity bisection along the target foliation and obtain global holomorphic solutions on the groupoid $\Pi_1 (\sfSigma_\sfDelta)$, provided we avoid paths that exchange sheets of the spectral cover; i.e., the ramification locus of \eqref{240801220010}.
We deduce that the analytic continuation $\phi$ of the Borel transform $\hat{\phi}$ is a globally defined holomorphic function on the groupoid $\Pi_1 (\sfSigma_\sfDelta \smallsetminus R)$.
As a result, the resurgence properties of $\phi$ can be read off from the geometry of the spectral cover and in particular its fundamental groupoid.

\paragraph{Acknowledgements.}
The author wishes to thank Dylan Allegretti, Tom Bridgeland, Samuel Crew, Marco Gualtieri, Kohei Iwaki, Olivier Marchal, Marta Mazzocco, Andrew Neitzke, Gerg\HungarianUmlaut{o} Nemes, Nicolas Orantin, and Shinji Sasaki for helpful discussions during various stages of this project.
Special thanks go to Marco Gualtieri and Marta Mazzocco for their support and encouragement in pursuing this project as well as the numerous suggestions for improving the manuscript.

\paragraph{Funding.}
This research was funded by the European Union's Horizon 2020 Research and Innovation Programme under the Marie Skłodowska-Curie Grant Agreement No.\,101026083 (\href{https://cordis.europa.eu/project/id/101026083}{\texttt{AbQuantumSpec}}), and also partially funded by the Leverhulme Trust through the Leverhulme Research Project grant \textit{Extended Riemann-Hilbert Correspondence, Quantum Curves and Mirror Symmetry}.

\section{Elements of Resurgence}
\label{240730090001}
\setcounter{equation}{0}
\counterwithin{equation}{section}

In this section, we introduce some basic ingredients of resurgence and resurgent perturbation theory from a geometric point of view.
Beware that the literature is unfortunately replete with somewhat inconsistent and often imprecise definitions (see \autoref{240730092212}) and we hope the contents of this section will help infuse some clarity into the subject.

In its most basic form, a \textit{resurgent power series} is a formal power series that has a convergent Borel transform which (a) admits a mostly unobstructed analytic continuation and (b) has at most exponential growth at infinity in most directions.
In the literature, property (a) appears in many not necessarily equivalent forms but under the same name, whilst property (b) is often assumed in some form but only implicitly.
These are properties of the analytic continuation of a convergent power series, so resurgence is ultimately a subject that concerns the singularities and global properties of multivalued holomorphic functions.
It is therefore best approached using geometric techniques, and this is the point of view we take.

\paragraph{Baseline notation and terminology.}
We denote by $\CC \cbrac{t} \subset \CC \bbrac{t}$ the subalgebra of convergent power series.
If $u \in \CC$, we let $\cal{O}_u \coleq \cal{O}_{\CC,u} \cong \CC \set{ t - u } \cong \CC \set{t}$ denote the algebra of germs of holomorphic functions at $u$.
A \textit{domain} in a Riemann surface is any nonempty, open, and connected subset.

By a \textit{real curve} on a Riemann surface $X$ we mean a continuous map $\gamma : I \to X$ from an interval $I \subset \RR$, defined up to orientation-preserving reparameterisation.
A real curve is a \textit{path} if $I$ is closed.
The \textit{limit} of a real curve $\gamma : [0,1) \to X$ is the limit of $\gamma (t) \in X$ as $t \to 1$.
This limit may not exist or may consist of more than one point.

\paragraph{Arcs and rays.}
Elements of $\RR / 2\pi\ZZ \cong \SS^1$ will be called \dfn{phases}, \dfn{directions}, or \dfn{rays}.
An \dfn{arc} $A \subset \RR / 2\pi\ZZ$ is a nonempty open interval; i.e., a nonempty, connected, and simply connected open subset.
We usually write arcs using the interval notation $(\alpha_1, \alpha_2)$ for a pair of rays $\alpha_1, \alpha_2 \in \SS^1$ respecting the standard orientation of the circle.
It has a well-defined length, or angle, $|A| \coleq \alpha_2 - \alpha_1 \leq 2\pi$.
When $|A| = 2\pi$, $A$ is the complement of a single point in $\RR / 2\pi\ZZ$.
If $\alpha \in \RR / 2\pi\ZZ$, the corresponding \dfn{conormal arc} is the arc of directions that form an acute angle with the \dfn{normal direction} $\alpha$, and we write:
\begin{eqn}
	\sfop{Arc} (\alpha) \coleq \big( \alpha - \tfrac{\pi}{2}, \alpha + \tfrac{\pi}{2} \big)
\qqtext{and}
	\sfop{Arc} (\alpha)^\perp \coleq \alpha
\fullstop
\end{eqn}
Similarly, if $A = (\alpha_1, \alpha_2) \subset \RR / 2\pi\ZZ$ has length $|A| \leq \pi$, its conormal arc is
\begin{eqn}
	\sfop{Arc} (A) 
		\coleq \Cup_{\alpha \in A} \sfop{arc} (\alpha)
		= \big( \alpha_1 - \tfrac{\pi}{2}, \alpha_2 + \tfrac{\pi}{2} \big)
\qqtext{and}
	\sfop{Arc} (A)^\perp \coleq A
\fullstop	
\end{eqn}
If $A \subset \RR / 2\pi\ZZ$ has length $|A| \geq \pi$, then the corresponding \dfn{normal arc} is
\begin{eqn}
	A^\perp
		\coleq \Cap_{\alpha \in A} \sfop{Arc} (\alpha)
\fullstop
\end{eqn}
In particular, if $A = (\alpha_1, \alpha_2)$ has length $|A| > \pi$, then $A^\perp = \big( \alpha_2 - \tfrac{\pi}{2}, \alpha_1 + \tfrac{\pi}{2} \big)$.
If $|A| = \pi$, then $A = \sfop{Arc} (\alpha)$ for some $\alpha$.

\paragraph{Halflines and halfstrips.}
If $\alpha \in \RR / 2\pi\ZZ$ and $u \in \CC$, a \dfn{halfline} is the parametrised real curve
\begin{eqn}
	\RR_{\alpha, u} \coleq u + e^{i \alpha} \RR_+ \coleq \set{ u + r e^{i \alpha} ~\big|~ r > 0 } \subset \CC
\fullstop
\end{eqn}
We write $\bar{\RR}_{\alpha, u} \coleq \set{ u + r e^{i \alpha} ~\big|~ r \geq 0 }$ for its closure.
A \dfn{halfstrip} is any domain of the form 
\begin{eqn}
	\SS_{\alpha, u} \coleq \set{ t ~\big|~ \op{dist} (t, \RR_{\alpha, u}) < \epsilon }
\end{eqn}
for some $\epsilon > 0$, where $\op{dist}$ is the usual Euclidean distance.
If $A \subset \RR$ is an arc, then an \dfn{infinite sector} \textit{with aperture $A$} is a domain of the form 
\begin{eqn}
	\SS_{A, u} \coleq \set{ t \in \RR_{\alpha, u} ~\big|~ \alpha \in A }
\fullstop
\end{eqn}
When $u = 0$, we simply write $\RR_\alpha, \SS_\alpha, \SS_A$.

\subsection{The Real-Oriented Blowup}

In this subsection, we recall the \textit{real-oriented blowup} of a Riemann surface at a point which is a convenient geometric construction that allows to keep track of directions in coordinate-independent way.
It is, in other words, an invariant way to make sense of polar coordinates.

\paragraph{Tangent circle.}
Let $X$ be a Riemann surface and $p \in X$ a point.
Recall that the tangent space $T_p X$ carries a rescaling $\CC^\times$-action.
Let $\RR_+ \subset \CC^\times$ be the subgroup of positive real numbers, and define the \dfn{tangent circle} \textit{at} $p$ as the quotient space
\begin{eqn}
	S_p X \coleq (T_p X)^\times \big/ \RR_+
\fullstop
\end{eqn}
Elements $\alpha \in S_p X$ are called (\dfn{tangential}) \dfn{directions} or \dfn{rays} at $p$.
The tangent circle is an $\SS^1$-torsor, hence homeomorphic to $\RR / 2\pi\ZZ$ but not canonically.
Such a homeomorphism can be set up by choosing a local coordinate centred at $p$.
In particular, the tangent circle $S_u \CC$ at any point $u \in \CC$ is canonically the standard circle $\SS^1 \cong \RR / 2\pi\ZZ$.

A (\dfn{tangential}) \dfn{arc} at $p$ is any nonempty open interval $A = (\alpha_1, \alpha_2) \subset S_p X$ whose \dfn{length}, denoted by $|A| = \alpha_2 - \alpha_1$, is well-defined and equals to the length of its image arc under any homeomorphism $S_p X \cong \RR / 2\pi\ZZ$.

\begin{defn}[real-oriented blowup]
\label{240911211435}	
Let $X^\ast \coleq X \smallsetminus \set{p}$.
The \dfn{real-oriented blowup} \textit{of $X$ at $p$} is the bordered Riemann surface
\begin{eqntag}
\label{240606160647}	
	[X : p] \coleq X^\ast \sqcup S_p X
\qqtext{with boundary}
	\partial [X : p] = S_p X
\fullstop{,}
\end{eqntag}
together with a holomorphic surjective map $\beta : [X : p] \too X$ called the \dfn{blowdown map}, which is the identity map when restricted to $X^\ast$ and the constant map with value $p$ when restricted to $S_p X$.
See \autoref{240606115513}.
\end{defn}

\begin{figure}[t]
\centering
\begin{tikzpicture}[scale=0.45, every node/.style={scale=0.9}]
\begin{scope}
\fill [grey] (-3,-5) rectangle (7,5);
\draw [blue, densely dashed, fill = blue, fill opacity = 0.3] (110:50pt) to [out=110,in=180] (1.5,4) arc (90:-100:4) to [out=165,in=260] (-100:50pt);
\fill [white, draw = red, ultra thick] (0,0) circle (50pt);
\draw [ultra thick, blue] (-100:50pt) arc (-100:110:50pt);
\node at (-1.5,4) {$[X:p]$};
\node at (110:45pt) [below] {$\alpha_2$};
\node at (-100:50pt) [above] {$\alpha_1$};
\node at (0:45pt) [above left, blue] {$A$};
\node at (25:130pt) [blue] {$\tilde{V}$};
\node at (180:20pt) [red] {$S_p X$};
\node at (8.5,0) {$\xrightarrow{~~\beta~~}$};
\end{scope}
\begin{scope}[xshift=13cm]
\fill [grey] (-3,-5) rectangle (7,5);
\draw [blue, densely dashed, fill = blue, fill opacity = 0.3] (110:3pt) to [out=110,in=180] (1.5,4) arc (90:-100:4) to [out=165,in=260] (-100:3pt);
\fill [red, draw = red, ultra thick] (0,0) circle (3pt) node [left, red] {$p$};
\node at (-2,4) {$X$};
\node at (25:130pt) [blue] {$V$};
\end{scope}
\end{tikzpicture}
\caption{The real-oriented blowup of $X$ at $p$ and a sector at $p$ with opening $A = (\alpha_1, \alpha_2)$.}
\label{240606115513}
\end{figure}
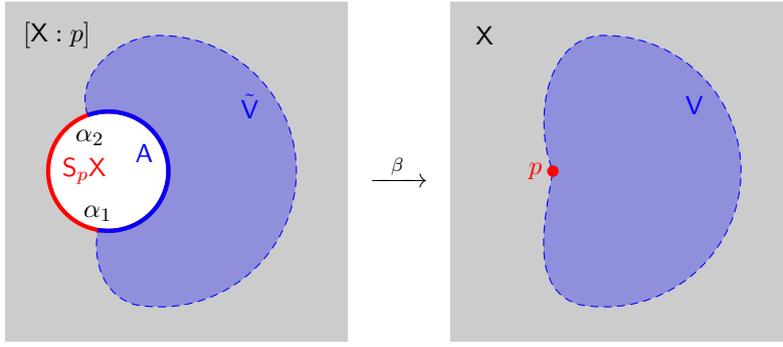

The real-oriented blowup lets us consider \textit{sectors} at points of a Riemann surface in a coordinate-independent way, defined as follows.

\begin{defn}[sectors]
\label{240606083834}
A \dfn{sector} at $p \in X$ is any simply connected domain
\begin{eqn}
	V \subset X^\ast = [X : p] \smallsetminus S_p X
\end{eqn}
with the property that there exists a (necessarily unique) open subset $\tilde{V} \subset [X : p]$ (in the sense of manifolds with boundary) which coincides with $V$ away from $p$ and intersects the tangent circle $S_p X$ in an arc $A$:
\begin{eqn}
	\tilde{V} \cap X^\ast = V
\qqtext{and}
	\tilde{V} \cap S_p X = A
\fullstop
\end{eqn}
\end{defn}
This arc $A$ is then called the \dfn{aperture} of $V$.
See \autoref{240606115513}.

\subsection{Endless Riemann Surfaces}

In this subsection, we introduce the basic object of study in resurgence: a Riemann surface spread over $\CC$.
It is a very general concept which arises in particular when constructing the analytic continuation of convergent power series.
It has a natural notion of singularities which captures the notion of singularities of analytic continuation.
We will be concerned with the case of \textit{endless Riemann surfaces} which are those that in a sense extend infinitely far away and their singularities are isolated.

\begin{defn}
\label{240913153305}
A \dfn{Riemann surface spread over} $\CC$ is a pair $(X, \pi)$ where $X$ is a connected Riemann surface and $\pi : X \to \CC$ is a nonconstant holomorphic map.
A morphism $(X, \pi) \to (X',\pi')$ of such pairs is itself a pair $(f,\T)$ consisting of a holomorphic map $f : X \to X'$ and a translation $\T : \CC \to \CC$ such that $\T \circ \pi = \pi' \circ f$.
\end{defn}

The map $\pi$ is sometimes called the \textit{projection} although we note that $\pi$ may not be surjective.
We also stress that $\pi$ is not even required to be a covering map (branched or not) over its image $\pi (X)$.
Note that it is important that we have restricted the automorphism group of $\CC$ to the subgroup of translations: it means that directions and straight line segments in $\CC$ have an invariant meaning on the Riemann surface $X$, but there is no invariantly distinguished fibre of $X$ that lies over the origin.

\paragraph{Analytic continuation.}
The most important source of Riemann surfaces spread over $\CC$ is the analytic continuation of holomorphic germs.
Recall that, if $Y$ is a Riemann surface, the \textit{étalé space} $|\cal{O}_Y|$ of the structure sheaf $\cal{O}_Y$ is the disjoint union of all the stalks $\cal{O}_{Y,q}$, $q \in Y$.
It is a topological space with a continuous surjective map $\pi : |\cal{O}_Y| \to Y$.
If $\hat{\phi} \in \cal{O}_{Y,q}$ is a holomorphic germ at $q$, the Riemann surface of $\hat{\phi}$ is by definition the connected component $X$ of the étalé space $|\cal{O}_Y|$ passing through the point $\hat{\phi} \in \cal{O}_{Y,q}$ which is now interpreted as a point $p \in X$ with projection $q = \pi (p)$.

The projection $\pi : |\cal{O}_Y| \to Y$ restricts to a holomorphic map $\pi : X \to Y$ which in general has ramification points, but the point $p$ is necessarily not one of them.
As a result, $\pi$ is a biholomorphism between a neighbourhood of $p$ and a neighbourhood of $q$ which yields an isomorphism $\pi^\ast : \cal{O}_{Y,q} \iso \cal{O}_{X,p}$.
Under this identification, the germ $\hat{\phi}$ is then the germ at $p$ of a global holomorphic function $\phi \in H^0 (X, \cal{O}_X)$.

When $Y = \CC$ and $\hat{\phi} (t) \in \CC \cbrac{t} \cong \cal{O}_0$ is a convergent power series, its Riemann surface $(X, \pi)$ is naturally spread over $\CC$ with a distinguished point $p \in X$ such that $\pi (p) = 0$.
In this case, we sometimes call the distinguished point $p$ the \textit{origin} in $X$.

\paragraph{Algebraic singularities and regular points.}
We wish to describe the \textit{singularities} of $(X,\pi)$ which are exceptional `points' in $X$, to be defined shortly, where the map $\pi$ exhibits special behaviour.
They are of two types --- \textit{algebraic} and \textit{nonalgebraic} --- which need to be treated separately.

The \dfn{algebraic singularities} (or \textit{algebraic ramification points}) are by definition the critical points of $\pi$.
That is, points $p \in X$ where the derivative $\pi_\ast : T_p X \to T_{\pi(p)} \CC$ vanishes.
Their projections to $\CC$ are called \dfn{algebraic singular values} (or \textit{algebraic branch points}).
All other points in $X$ are called \dfn{regular points}.

For any $p \in X$ with projection $u \in \CC$, the derivative $\pi_\ast : T_p X \to T_u \CC$ induces a natural map
\begin{eqntag}
\label{241010174155}
	S_p X \to S_u \CC = \SS^1
\fullstop
\end{eqntag}
If $p$ is a regular point, this map is an isomorphism $S_p X \iso \SS^1$, allowing us to canonically identify every tangent ray $\alpha \in S_p X$ with a phase $\alpha \in \RR / 2\pi\ZZ$.
On the other hand, if $p$ is an algebraic singularity, then this map is an $m$-fold covering map where $m$ is the ramification index of $\pi$ at $p$.
Nevertheless, every tangent ray $\tilde{\alpha} \in S_p X$ at an algebraic singularity has a well-defined phase $\alpha \in \RR / 2\pi\ZZ$.

\paragraph{Nonalgebraic singularities.}
Nonalgebraic singularities are more complicated to describe; one approach is as follows.
Let $\xi \in \CC$ be any point in the closure of $\pi (X)$.
For any domain $\UU \subset \CC$ around $\xi$, consider the collection of connected components of the preimage $\pi^\inv (\UU)$ in $X$.
Since $\pi$ is holomorphic and nonconstant, any such collection is necessarily countable.

Let $F_\xi$ denote the set of equivalence classes of domains $U \subset X$ such that $U$ is a connected component of $\pi^\inv (\UU)$ for some domain $\UU \subset \CC$ around $\xi$.
Here, two such $U,U' \subset X$ are considered equivalent if there is another domain $U'' \subset U \cap U'$ which is a connected component of $\pi^\inv (\UU'')$ for some domain $\UU'' \subset \CC$ around $\xi$.
We think of $F_\xi$ as the `fibre' of $\pi$ above $\xi$, although this set may be strictly larger than $\pi^\inv (\xi)$.
Indeed, it may happen that $\pi^\inv (\xi)$ is empty yet $F_\xi$ is countably infinite.
For example, if $(X, \pi) = (\CC, \exp)$, then $\xi = 0$ is in the closure of $\pi (X)$, but $\pi^\inv (0) = \emptyset$ yet $F_0 \cong \ZZ$.

Given an equivalence class $\frak{p} \in F_\xi$, we consider the total intersection
\begin{eqntag}
\label{240906135306}
	|\, \frak{p} \,| \coleq \Cap_{U \in \frak{p} } U \subset X
\fullstop
\end{eqntag}
Since $\pi$ is a nonconstant holomorphic map from a Riemann surface, this set is either empty or a single point.
If it is nonempty and equals $\set{p}$ for some $p \in X$, then $\xi = \pi (p)$ and $p$ is either a regular point or an algebraic singularity.
On the other hand, if it is empty, then we call $\frak{p}$ a \dfn{nonalgebraic singularity} and there are two possibilities: $\frak{p}$ is either a \textit{removable singularity} or a \textit{transcendental singularity}.

\begin{defn*}
\label{241016104312}
A \dfn{removable singularity} is characterised by the fact that there exist a domain $U^\ast \in \frak{p}$ and a simply connected domain $\UU \subset \CC$ around $\xi$ such that $\pi$ restricts to an isomorphism $\pi : U^\ast \iso \UU \smallsetminus \set{\xi}$.
On the other hand, $\frak{p}$ is a \dfn{transcendental singularity} \textit{over} $\xi \in \CC$ if the set \eqref{240906135306} is empty but $\frak{p}$ is not a removable singularity.
In this case, $\xi$ is called a \dfn{transcendental singular value}.
\end{defn*}

A removable singularity can always be `filled in' by adding a new abstract point $p$ to $X$ and extending $\pi$ to a holomorphic map $U \to \UU$ by defining $U \coleq U^\ast \sqcup \set{p}$ and $\pi (p) \coleq \xi$.
The point $p$ in the new Riemann surface $X \sqcup \set{p}$ is regular and $|\,\frak{p}\,| = \set{p}$.
At the same time, a transcendental singularity cannot in general be added to $X$ as a new point in a way that extends the Riemann surface structure.

\paragraf
All algebraic and nonalgebraic singularities are collectively called \dfn{singularities} of $(X, \pi)$ and their collection is denoted by $\sfop{Sing} (X,\pi)$.
We also denote by $\Xi \subset \CC$ the collection of all singular values.
Given any singularity $\frak{p} \in \sfop{Sing} (X,\pi)$, we call any representative $U \in \frak{p}$ a \dfn{domain around $\frak{p}$} even though, for example, a transcendental singularity is \textit{not} a point in the Riemann surface $X$.
If $p_0 \in X$ is a regular point over $u = \pi (p_0)$ and $\frak{p} \in \sfop{Sing} (X, \pi)$ is a singularity over $\xi \in \CC$, the \dfn{phase} \textit{of $\frak{p}$ relative to $p_0$} is the phase $\arg (\xi - u) \in \RR / 2\pi\ZZ$.

\paragraf
Algebraic singularities form a necessarily discrete subset of $X$.
In contrast, nonalgebraic singularities are not points in $X$, so it does not make sense to say that they form a discrete subset.
Instead, we formulate this notion as follows.

\begin{defn*}[isolated singularities]
\label{241005134703}
A singularity $\frak{p} \in \sfop{Sing} (X,\pi)$ is \dfn{isolated} if there is a domain $U \subset X$ around $\frak{p}$ with the following property: every singularity $\frak{q}$ different from $\frak{p}$ has a domain $U' \subset X$ around it such that $U \cap U' = \emptyset$.
\end{defn*}

\begin{defn*}[endless analytic continuation]
\label{240913161447}
An \dfn{endless Riemann surface} is a Riemann surface spread over $\CC$ with only isolated singularities.
A holomorphic germ $\hat{\phi} \in \CC \set{t}$ is said to admit an \dfn{endless analytic continuation} if it can be analytically continued to an endless Riemann surface.
\end{defn*}

Clearly, all removable singularities are isolated.
Also, note that even if $(X,\pi)$ is endless, the subset $\Xi \subset \CC$ of all singular values may in general have accumulation points or even be dense somewhere or everywhere.

\paragraf
Amongst the most well-behaved classes of endless Riemann surfaces are branched covers.
If $\pi$ is a branched cover, then $(X,\pi)$ is endless and $\Xi \subset \CC$ is necessarily a discrete subset; in this case, we say $(X, \pi)$ is of \dfn{discrete type}.
If, furthermore, $\Xi$ is finite, then we say $(X,\pi)$ is of \dfn{finite type}.
In this paper, we encounter another special class of endless Riemann surfaces which are universal covers of other endless Riemann surfaces that are themselves algebraic branched covers of $\CC$.

\begin{defn*}
\label{241011072440}
An endless Riemann surface $(X,\pi)$ is of
\begin{enumerate}
\item \dfn{algebraic type} if all singularities are algebraic;
\item \dfn{log-algebraic type} if there is an endless Riemann surface $(X',\pi')$ of algebraic type, a discrete subset $\Gamma' \subset X'$, and an isomorphism $(X,\pi) \iso (\tilde{X}',\pi')$ where $\tilde{X}'$ is a universal cover of $X' \smallsetminus \Gamma'$.
\end{enumerate}
\end{defn*}

\subsection{Visible Singularities and Stokes Rays}

The structure of being spread over $\CC$ equips a Riemann surface with a notion of straight lines, distances, and directions.
Heuristically, these notions can be used to distinguish some singularities (the \textit{visible singularities}) as being more directly accessible from the point of view of a particular point on the Riemann surface.
Other singularities cannot be accessed directly, but amongst them are those (the \textit{observable singularities}) which are almost directly accessible: they can be reached upon following a straight line while circumventing a few singularities that are encountered along the way.
These observations are central to the resurgence theory, and make them precise in this subsection.

\begin{defn}[asymptotic and infinite paths]
\label{241010104825}
Suppose $(X,\pi)$ is a Riemann surface spread over $\CC$.
An \dfn{asymptotic path} on $X$ is a real curve $\gamma : [0,1) \to X$ which has no limit in $X$ but its projection $\pi \circ \gamma : [0,1) \to \CC$ has a limit in $\CC$.
Similarly, an \dfn{infinite path} on $X$ is a real curve $\gamma : [0,1) \to X$ whose projection $\pi \circ \gamma$ limits to the point at infinity $\infty$ in the Riemann sphere $\bar{\CC}$.

If $\frak{p} \in \sfop{Sing} (X,\pi)$ is a singularity over $\xi \in \CC$, we say that a real curve $\gamma : [0,1) \to X$ \dfn{terminates at $\frak{p}$} if the projection $\pi \circ \gamma$ limits to $\xi$ and any domain $U \subset X$ around $\frak{p}$ contains the tail of $\gamma$; i.e., there is some $t_0 \in [0,1)$ such that $\gamma (t) \in U$ for all $t \geq t_0$.
A real curve $\gamma : (0,1] \to X$ that \dfn{starts at $\frak{p}$} is defined similarly.
\end{defn}

\begin{defn}[geodesics]
\label{241015170927}
A \dfn{geodesic} on $X$ is a smooth real curve $\gamma : I \to X$ that avoids algebraic singularities and whose projection $\pi \circ \gamma$ is a (possibly infinitely long) line segment in $\CC$.
The direction $\alpha \in \RR / 2\pi\ZZ$ of this line segment is the \dfn{phase} of $\gamma$.
Said differently, a smooth real curve $\gamma : I \to X$ is a geodesic if its tangent vector $\dot{\gamma}$ has constant phase $\alpha$ at every point of $\gamma$ in the sense that $\d{t} (\pi_\ast \dot{\gamma}) \in e^{i \alpha} \RR_+$ where $\d{t}$ is the canonical translation-invariant differential on $\CC$.

For any phase $\alpha$, any regular point $p \in X$ lies on a unique maximal geodesic with phase $\alpha$.
Likewise, if $p \in X$ is a regular point, then using the canonical identification $S_p X \cong \RR / 2\pi\ZZ$ given by the derivative of $\pi$, tangential rays at $p$ are in one-to-one correspondence with maximal geodesics starting at $p$.
\end{defn}

\begin{defn}[visible singularities; regular and Stokes rays]
\label{240926145223}
Let $(X,\pi)$ be a Riemann surface spread over $\CC$ and $p \in X$ a regular point.
A singularity $\frak{q} \in \sfop{Sing} (X,\pi)$ is called a \dfn{singularity visible from $p$} if there is a (necessarily unique) geodesic that starts at $p$ and terminates at $\frak{q}$.
The corresponding singular value is also called a \dfn{visible singular value}.
In this case, we call this geodesic a \dfn{Stokes geodesic}, and its phase $\alpha \in \RR / 2\pi\ZZ \cong S_p X$ a \dfn{Stokes ray} \textit{at $p$}.
All other rays $\alpha \in S_p X$ are called \dfn{regular rays}.
If $\alpha \in S_p X$ is a regular ray, the maximal geodesic starting at $p$ with phase $\alpha$ is an infinite path called a \dfn{regular geodesic} or \dfn{infinite geodesic}.
\end{defn}

We denote by $\sfop{Sing}^0_p (X,\pi) \subset \sfop{Sing} (X,\pi)$ the subset of singularities visible from $p$.
Denote by $\Xi^0_p \subset \Xi \subset \CC$ the subset of all singular values visible from $p$.
If $A \subset \RR / 2\pi\ZZ$ is an arc, put $\Xi^0_{p,A} \coleq \Xi^0_p \cap \SS_A$ and let $\sfop{Sing}^0_{p,A} (X, \pi)$ be subset of singularities visible from $p$ whose phase relative to $p$ belongs to $A$.

\paragraf
\autoref*{240926145223} can be phrased more concretely as follows.
By a \dfn{regular lift} of a real curve $\wp : I \to \CC$ to $X$ we mean a lift of $\wp$ to a real curve $\gamma : I \to X$ such that $\pi$ restricts to a local biholomorphism in a neighbourhood of $\gamma$.
Thus, for example, a lift of $\wp$ that passes through a ramification point is not regular.

Then, upon translating the projection $\pi (p) \in \CC$ to the origin, a singularity $\frak{q}$ over $\xi \in \CC$ is visible from $p$ if and only if the half-open line segment $[0, \xi)$ has a regular lift to $X$ starting at $p$, but the closed line segment $[0, \xi]$ does not.
In this case, the phase $\alpha = \arg \xi$ is a Stokes ray at $p$, and the lift of $[0, \xi)$ that starts at $p$ is the Stokes geodesic with phase $\alpha$.
The point $\xi$ is a visible singular value, and the closed halfline $\bar{\RR}_{\alpha, \xi}$ emanating from $\xi$ is called the \dfn{Stokes cut} at $\xi$; see \autoref{240820100638}.
On the other hand, if the halfline $\RR_{\alpha}$ has a regular lift to $X$ starting at $p$ then $\alpha$ is a regular ray at $p$, and this lift is an infinite geodesic path.

\begin{figure}[t]
\centering
\begin{adjustwidth}{-1.5cm}{-1.5cm}
\begin{subfigure}[t]{0.4\textwidth}
\centering
\usetikzlibrary{shapes.misc}

\tikzset{cross/.style={cross out, draw=black, ultra thick, minimum size=2*(#1-\pgflinewidth), inner sep=1pt, outer sep=1pt}}

\begin{tikzpicture}
\begin{scope}
\clip (-2,-2) rectangle (2,2);
\fill [grey] (-2,-2) rectangle (2,2);
\fill [black] (0,0) circle (1pt);
\node at (135:2.25) {$\CC$};
\node [cross, red] at (0:1) {};
\node [cross, red] at (60:1.1) {};
\node [cross, red] at (125:0.7) {};
\node [cross, red] at (170:1.2) {};
\node [cross, red] at (220:1.5) {};
\node [cross, red] at (260:1.3) {};
\node [cross, red] at (295:1.1) {};
\node [cross, red] at (310:2.1) {};
\draw [dashed, red, ultra thick] (0:1.1) -- (0:4);
\draw [dashed, red, ultra thick] (60:1.2) -- (60:4);
\draw [dashed, red, ultra thick] (125:0.8) -- (125:4);
\draw [dashed, red, ultra thick] (170:1.3) -- (170:4);
\draw [dashed, red, ultra thick] (220:1.6) -- (220:4);
\draw [dashed, red, ultra thick] (260:1.4) -- (260:4);
\draw [dashed, red, ultra thick] (295:1.2) -- (295:4);
\draw [dashed, red, ultra thick] (310:2.2) -- (310:4);
\end{scope}
\end{tikzpicture}
\caption{Visible singular values (red crosses) and Stokes cuts (red dashed rays).
The union of Stokes cuts $\Xi^+_p$ is the subset in red.
Its complement is the holomorphic star $\EE_p$.}
\label{240820100638}
\end{subfigure}
\quad
\begin{subfigure}[t]{0.4\textwidth}
\centering
\begin{tikzpicture}
\begin{scope}
\clip (-2,-2) rectangle (2,2);
\draw [black, ultra thick] (0,0) circle (1.5);
\fill [red] (0:1.5) circle (2pt);
\fill [red] (60:1.5) circle (2pt);
\fill [red] (125:1.5) circle (2pt);
\fill [red] (170:1.5) circle (2pt);
\fill [red] (220:1.5) circle (2pt);
\fill [red] (260:1.5) circle (2pt);
\fill [red] (295:1.5) circle (2pt);
\fill [red] (310:1.5) circle (2pt);
\draw [red, ->] (0,0) -- (0:1.4);
\draw [red, ->] (0,0) -- (60:1.4);
\draw [red, ->] (0,0) -- (125:1.4);
\draw [red, ->] (0,0) -- (170:1.4);
\draw [red, ->] (0,0) -- (220:1.4);
\draw [red, ->] (0,0) -- (260:1.4);
\draw [red, ->] (0,0) -- (295:1.4);
\draw [red, ->] (0,0) -- (310:1.4);
\end{scope}
\end{tikzpicture}
\caption{Stokes diagram.
Red points on the circle are Stokes rays; all other points are regular rays.}
\label{240820101437}
\end{subfigure}
\quad
\begin{subfigure}[t]{0.3\textwidth}
\centering
\begin{tikzpicture}
\draw [black, ultra thick] (0,0) circle (1.5);
\fill [red] (0:1.5) circle (1.5pt);
\fill [red] (30:1.5) circle (1.5pt);
\fill [red] (45:1.5) circle (1.5pt);
\fill [red] (55:1.5) circle (1.5pt);
\fill [red] (63:1.5) circle (1.5pt);
\fill [red] (70:1.5) circle (1.5pt);
\fill [red] (75:1.5) circle (1.5pt);
\fill [red] (79:1.5) circle (1.5pt);
\fill [red] (83:1.5) circle (1.5pt);
\fill [red] (85:1.5) circle (1.5pt);
\fill [red] (87:1.5) circle (1.5pt);
\begin{scope}[xscale=-1]
\fill [red] (0:1.5) circle (1.5pt);
\fill [red] (30:1.5) circle (1.5pt);
\fill [red] (45:1.5) circle (1.5pt);
\fill [red] (55:1.5) circle (1.5pt);
\fill [red] (63:1.5) circle (1.5pt);
\fill [red] (70:1.5) circle (1.5pt);
\fill [red] (75:1.5) circle (1.5pt);
\fill [red] (79:1.5) circle (1.5pt);
\fill [red] (83:1.5) circle (1.5pt);
\fill [red] (85:1.5) circle (1.5pt);
\fill [red] (87:1.5) circle (1.5pt);
\end{scope}
\fill [red] (-90:1.5) circle (2pt);

\begin{scope}[yscale=-1]
\fill [red] (30:1.5) circle (1.5pt);
\fill [red] (45:1.5) circle (1.5pt);
\fill [red] (55:1.5) circle (1.5pt);
\fill [red] (63:1.5) circle (1.5pt);
\fill [red] (70:1.5) circle (1.5pt);
\fill [red] (75:1.5) circle (1.5pt);
\fill [red] (79:1.5) circle (1.5pt);
\fill [red] (83:1.5) circle (1.5pt);
\fill [red] (85:1.5) circle (1.5pt);
\fill [red] (87:1.5) circle (1.5pt);
\end{scope}
\begin{scope}[xscale=-1,yscale=-1]
\fill [red] (30:1.5) circle (1.5pt);
\fill [red] (45:1.5) circle (1.5pt);
\fill [red] (55:1.5) circle (1.5pt);
\fill [red] (63:1.5) circle (1.5pt);
\fill [red] (70:1.5) circle (1.5pt);
\fill [red] (75:1.5) circle (1.5pt);
\fill [red] (79:1.5) circle (1.5pt);
\fill [red] (83:1.5) circle (1.5pt);
\fill [red] (85:1.5) circle (1.5pt);
\fill [red] (87:1.5) circle (1.5pt);
\end{scope}
\fill [black] (90:1.5) circle (1.5pt);

\draw [red, ->] (0,0) -- (0:1.4);
\draw [red, ->] (0,0) -- (30:1.4);
\draw [red, ->] (0,0) -- (45:1.4);
\draw [red, ->] (0,0) -- (55:1.4);
\draw [red, ->] (0,0) -- (63:1.4);
\draw [red, ->] (0,0) -- (70:1.4);
\draw [red, ->] (0,0) -- (75:1.4);
\fill [red] (80:1.25) circle (0.5pt);
\fill [red] (83:1.25) circle (0.5pt);
\fill [red] (86:1.25) circle (0.5pt);
\begin{scope}[xscale=-1]
\draw [red, ->] (0,0) -- (0:1.4);
\draw [red, ->] (0,0) -- (30:1.4);
\draw [red, ->] (0,0) -- (45:1.4);
\draw [red, ->] (0,0) -- (55:1.4);
\draw [red, ->] (0,0) -- (63:1.4);
\draw [red, ->] (0,0) -- (70:1.4);
\draw [red, ->] (0,0) -- (75:1.4);
\fill [red] (80:1.25) circle (0.5pt);
\fill [red] (83:1.25) circle (0.5pt);
\fill [red] (86:1.25) circle (0.5pt);
\end{scope}

\draw [red, ->] (0,0) -- (-90:1.4);

\begin{scope}[yscale=-1]
\draw [red, ->] (0,0) -- (0:1.4);
\draw [red, ->] (0,0) -- (30:1.4);
\draw [red, ->] (0,0) -- (45:1.4);
\draw [red, ->] (0,0) -- (55:1.4);
\draw [red, ->] (0,0) -- (63:1.4);
\draw [red, ->] (0,0) -- (70:1.4);
\draw [red, ->] (0,0) -- (75:1.4);
\fill [red] (80:1.25) circle (0.5pt);
\fill [red] (83:1.25) circle (0.5pt);
\fill [red] (86:1.25) circle (0.5pt);
\begin{scope}[xscale=-1]
\draw [red, ->] (0,0) -- (0:1.4);
\draw [red, ->] (0,0) -- (30:1.4);
\draw [red, ->] (0,0) -- (45:1.4);
\draw [red, ->] (0,0) -- (55:1.4);
\draw [red, ->] (0,0) -- (63:1.4);
\draw [red, ->] (0,0) -- (70:1.4);
\draw [red, ->] (0,0) -- (75:1.4);
\fill [red] (80:1.25) circle (0.5pt);
\fill [red] (83:1.25) circle (0.5pt);
\fill [red] (86:1.25) circle (0.5pt);
\end{scope}
\end{scope}
\end{tikzpicture}
\caption{Stokes diagram.
The rays in the directions $\tfrac{\pi}{2}$ and $-\tfrac{\pi}{2}$ are respectively a regular ray and a Stokes ray which are either semi-stable or unstable, but not stable; it is not possible to distinguish semi-stability from instability from the Stokes diagram.
All other rays are stable.}
\label{241022152328}
\end{subfigure}
\end{adjustwidth}
\caption{Visible singularities, Stokes cuts, and Stokes diagrams.}
\end{figure}
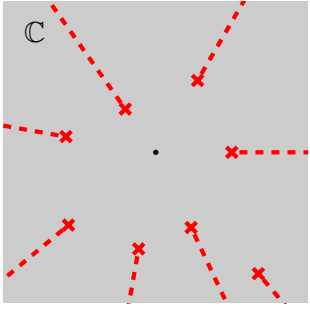
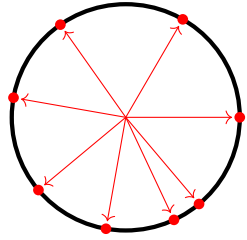
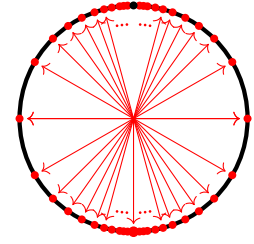

\begin{prop}
\label{241007110600}
Suppose $(X,\pi)$ is an endless Riemann surface and $p \in X$ is a regular point.
Then the number of Stokes rays at $p$ is at most countable.
Consequently, the union of all Stokes cuts is a closed subset of $\CC$.
Upon translating $\pi (p)$ to the origin for convenience, the complement of Stokes cuts $\EE_p \subset \CC$ is a star-shaped domain around the origin containing the (not necessarily open) cone consisting of all halflines $\RR_{\alpha}$ corresponding to regular rays $\alpha$ at $p$.
Moreover, $\EE_p$ lifts biholomorphically via $\pi$ to a simply connected neighbourhood $E_p \subset X$ of $p$.
\end{prop}

We call the star-shaped domain $E_p \cong \EE_p$ the \dfn{holomorphic star} at $p$.

\begin{proof*}
Since the singularities of $\pi$ are isolated, the subset of visible singular values $\Xi^0_p$ is a discrete and countable subset.
Therefore, the subset of all Stokes cuts is closed because it is a countable union of a locally finite collection of closed subsets.
All other assertions are just simple observations.
\end{proof*}

\paragraph{Resurgent Stokes diagrams.}
Thus, Stokes rays in an endless Riemann surface form at worst a countable subset of $\SS^1 \cong S_p X$ which we picture as in \autoref{240820101437} and call the \dfn{Stokes diagram} of $X$ at $p$, or \dfn{resurgent Stokes diagram} in order to distinguish it from a \textit{geometric Stokes diagram} which we will encounter later in this paper.
In general, a resurgent Stokes diagram can be rather complicated: it can have accumulation points or even be dense in a part or all of $\SS^1$; see \autoref{241022152328}.
To describe special kinds of better-behaved situations, we introduce the following terminology.

\begin{defn*}[stability in resurgent Stokes diagrams]
\label{240809103556}
Let $(X, \pi)$ be an endless Riemann surface and $p \in X$ a regular point with projection $\pi (p)$ translated to the origin.
We say that a ray $\alpha \in S_p X$ is
\begin{enumerate}
\item \dfn{stable} if $\alpha$ is contained in an arc $A = (\alpha_-, \alpha_+) \subset S_p X$ of consisting only of regular rays except possibly $\alpha$ itself, see \autoref{241022153737};
\item \dfn{semi-stable} if the halfline $\RR_{\alpha}$ has a halfstrip neighbourhood $\SS_\alpha$ which does not contain any Stokes cuts except possibly the one in the direction $\alpha$, see \autoref{241022161242};
\item \dfn{unstable} otherwise.	
\end{enumerate}
\end{defn*}

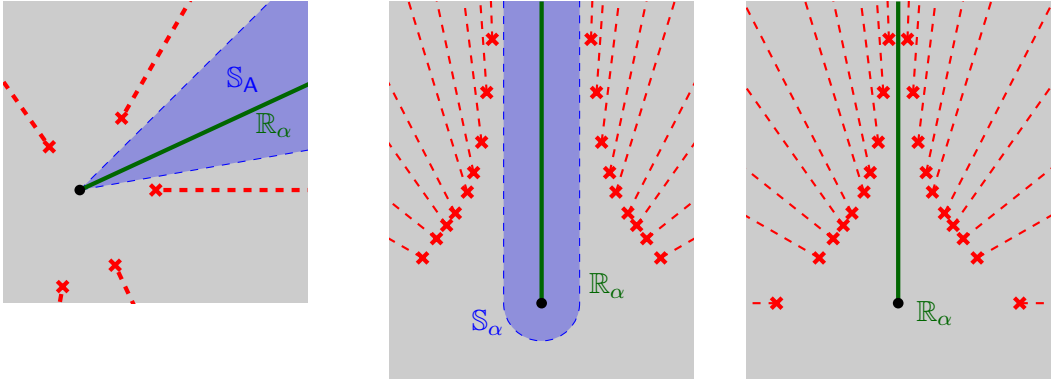
\begin{figure}[t]
\centering
\begin{subfigure}[t]{0.3\textwidth}
\centering
\usetikzlibrary{shapes.misc}

\tikzset{cross/.style={cross out, draw=black, ultra thick, minimum size=2*(#1-\pgflinewidth), inner sep=1pt, outer sep=1pt}}

\begin{tikzpicture}
\begin{scope}
\clip (-2,-2) rectangle (2,2);
\fill [grey] (-2,-2) rectangle (2,2);
\begin{scope}[shift={(-1,-0.5)}]
\node at (135:2.25) {$\CC$};
\node [cross, red] at (0:1) {};
\node [cross, red] at (60:1.1) {};
\node [cross, red] at (125:0.7) {};
\node [cross, red] at (170:1.2) {};
\node [cross, red] at (220:1.5) {};
\node [cross, red] at (260:1.3) {};
\node [cross, red] at (295:1.1) {};
\node [cross, red] at (310:2.1) {};
\draw [dashed, red, ultra thick] (0:1.1) -- (0:4);
\draw [dashed, red, ultra thick] (60:1.2) -- (60:4);
\draw [dashed, red, ultra thick] (125:0.8) -- (125:4);
\draw [dashed, red, ultra thick] (170:1.3) -- (170:4);
\draw [dashed, red, ultra thick] (220:1.6) -- (220:4);
\draw [dashed, red, ultra thick] (260:1.4) -- (260:4);
\draw [dashed, red, ultra thick] (295:1.2) -- (295:4);
\draw [dashed, red, ultra thick] (310:2.2) -- (310:4);
\draw [dashed, blue, fill=blue, fill opacity = 0.3] (0:0) -- (10:4) -- (4,4) -- (45:4) -- cycle;
\draw [darkgreen, ultra thick] (0:0) -- (25:4);
\fill [black] (0,0) circle (2pt);
\end{scope}
\end{scope}
\node [darkgreen] at (13:1.6) {$\RR_\alpha$};
\node [blue] at (40:1.5) {$\SS_A$};
\end{tikzpicture}
\caption{Stable regular ray $\RR_\alpha$ contained in an infinite sector $\SS_A$ which is contained inside the holomorphic star at $p$.}
\label{241022153737}
\end{subfigure}
\qquad
\begin{subfigure}[t]{0.3\textwidth}
\centering
\usetikzlibrary{shapes.misc}

\tikzset{cross/.style={cross out, draw=black, ultra thick, minimum size=2*(#1-\pgflinewidth), inner sep=1pt, outer sep=1pt}}

\begin{tikzpicture}
\begin{scope}
\clip (-2,-1) rectangle (2,4);
\fill [grey] (-2,-1) rectangle (2,4);
\begin{scope}[xshift=15]
\node [cross, red] at (0:1.6) {};
\node [cross, red] at (30:1.2) {};
\node [cross, red] at (45:1.21) {};
\node [cross, red] at (55:1.26) {};
\node [cross, red] at (63:1.34) {};
\node [cross, red] at (73:1.54) {};
\node [cross, red] at (78:1.77) {};
\node [cross, red] at (83:2.15) {};
\node [cross, red] at (86:2.8) {};
\node [cross, red] at (88:3.5) {};
\draw [dashed, red, thick] (0:1.6) -- (0:4);
\draw [dashed, red, thick] (30:1.2) -- (30:4);
\draw [dashed, red, thick] (45:1.21) -- (40:6);
\draw [dashed, red, thick] (55:1.26) -- (55:6);
\draw [dashed, red, thick] (63:1.34) -- (63:6);
\draw [dashed, red, thick] (73:1.54) -- (73:6);
\draw [dashed, red, thick] (78:1.77) -- (78:6);
\draw [dashed, red, thick] (83:2.15) -- (83:6);
\draw [dashed, red, thick] (86:2.8) -- (86:6);
\draw [dashed, red, thick] (88:3.5) -- (88:6);
\end{scope}
\begin{scope}[xscale=-1, xshift=15]
\node [cross, red] at (0:1.6) {};
\node [cross, red] at (30:1.2) {};
\node [cross, red] at (45:1.21) {};
\node [cross, red] at (55:1.26) {};
\node [cross, red] at (63:1.34) {};
\node [cross, red] at (73:1.54) {};
\node [cross, red] at (78:1.77) {};
\node [cross, red] at (83:2.15) {};
\node [cross, red] at (86:2.8) {};
\node [cross, red] at (88:3.5) {};
\draw [dashed, red, thick] (0:1.6) -- (0:4);
\draw [dashed, red, thick] (30:1.2) -- (30:4);
\draw [dashed, red, thick] (45:1.21) -- (40:6);
\draw [dashed, red, thick] (55:1.26) -- (55:6);
\draw [dashed, red, thick] (63:1.34) -- (63:6);
\draw [dashed, red, thick] (73:1.54) -- (73:6);
\draw [dashed, red, thick] (78:1.77) -- (78:6);
\draw [dashed, red, thick] (83:2.15) -- (83:6);
\draw [dashed, red, thick] (86:2.8) -- (86:6);
\draw [dashed, red, thick] (88:3.5) -- (88:6);
\end{scope}
\draw [dashed, blue, fill=blue, fill opacity = 0.3] (0:0.5) arc (0:-180:0.5) -- (-0.5,6) -- (0.5,6) -- cycle;
\draw [darkgreen, ultra thick] (0:0) -- (90:4);
\fill [black] (0,0) circle (2pt);
\node [blue] at (200:0.75) {$\SS_\alpha$};
\node [darkgreen] at (15:0.9) {$\RR_\alpha$};
\end{scope}
\end{tikzpicture}
\caption{Semi-stable regular ray $\RR_\alpha$ contained in a halfstrip $\SS_\alpha$}
\label{241022161242}
\end{subfigure}
\quad
\begin{subfigure}[t]{0.3\textwidth}
\centering
\usetikzlibrary{shapes.misc}

\tikzset{cross/.style={cross out, draw=black, ultra thick, minimum size=2*(#1-\pgflinewidth), inner sep=1pt, outer sep=1pt}}

\begin{tikzpicture}
\begin{scope}
\clip (-2,-1) rectangle (2,4);
\fill [grey] (-2,-1) rectangle (2,4);
\node [cross, red] at (0:1.6) {};
\node [cross, red] at (30:1.2) {};
\node [cross, red] at (45:1.21) {};
\node [cross, red] at (55:1.26) {};
\node [cross, red] at (63:1.34) {};
\node [cross, red] at (73:1.54) {};
\node [cross, red] at (78:1.77) {};
\node [cross, red] at (83:2.15) {};
\node [cross, red] at (86:2.8) {};
\node [cross, red] at (88:3.5) {};
\draw [dashed, red, thick] (0:1.6) -- (0:4);
\draw [dashed, red, thick] (30:1.2) -- (30:4);
\draw [dashed, red, thick] (45:1.21) -- (40:6);
\draw [dashed, red, thick] (55:1.26) -- (55:6);
\draw [dashed, red, thick] (63:1.34) -- (63:6);
\draw [dashed, red, thick] (73:1.54) -- (73:6);
\draw [dashed, red, thick] (78:1.77) -- (78:6);
\draw [dashed, red, thick] (83:2.15) -- (83:6);
\draw [dashed, red, thick] (86:2.8) -- (86:6);
\draw [dashed, red, thick] (88:3.5) -- (88:6);
\begin{scope}[xscale=-1]
\node [cross, red] at (0:1.6) {};
\node [cross, red] at (30:1.2) {};
\node [cross, red] at (45:1.21) {};
\node [cross, red] at (55:1.26) {};
\node [cross, red] at (63:1.34) {};
\node [cross, red] at (73:1.54) {};
\node [cross, red] at (78:1.77) {};
\node [cross, red] at (83:2.15) {};
\node [cross, red] at (86:2.8) {};
\node [cross, red] at (88:3.5) {};
\draw [dashed, red, thick] (0:1.6) -- (0:4);
\draw [dashed, red, thick] (30:1.2) -- (30:4);
\draw [dashed, red, thick] (45:1.21) -- (40:6);
\draw [dashed, red, thick] (55:1.26) -- (55:6);
\draw [dashed, red, thick] (63:1.34) -- (63:6);
\draw [dashed, red, thick] (73:1.54) -- (73:6);
\draw [dashed, red, thick] (78:1.77) -- (78:6);
\draw [dashed, red, thick] (83:2.15) -- (83:6);
\draw [dashed, red, thick] (86:2.8) -- (86:6);
\draw [dashed, red, thick] (88:3.5) -- (88:6);
\end{scope}
\draw [darkgreen, ultra thick] (0:0) -- (90:4);
\fill [black] (0,0) circle (2pt);
\node [darkgreen] at (-15:0.5) {$\RR_\alpha$};
\end{scope}
\end{tikzpicture}
\caption{Unstable regular ray.}
\label{241022161330}
\end{subfigure}
\caption{Stability of regular rays. Stability of Stokes rays is similar.}
\label{241015181238}
\end{figure}

\paragraph{Visible directions.}
Suppose $(X, \pi)$ is endless, $p \in X$ is a regular point with projection $\pi (p)$ translated to the origin, and $\alpha \in S_p X$ is a Stokes ray with visible singularity $\frak{q} \in \sfop{Sing}_{p,\alpha}^0 (X, \pi)$ over the visible singular value $\xi \in \RR_\alpha$.
The halfline $\RR_\alpha$ does not admit a regular lift to an infinite path on $X$ starting at $p$.
However, the Stokes cut $\bar{\RR}_{\alpha, \xi}$ is isolated and contained in an open subset $\SS \subset \CC$ which is separated from all other Stokes cuts.
We can therefore deform the halfline $\RR_\alpha$ by moving the part contained in $\SS$ to an infinite path in the direction $\alpha$ contained in $\SS \smallsetminus \bar{\RR}_{\alpha, \xi}$ that runs either along the left or along the right side of the Stokes cut.
These infinite paths in $\CC$, which we can denote by $\RR^\rm{L}_\alpha, \RR^\rm{R}_\alpha$, have regular lifts to infinite paths $\gamma^\rm{L}, \gamma^\rm{R}$ on $X$ starting at $p$.
Furthermore, these infinite paths are uniquely determined up to homotopy inside the holomorphic star $\EE_p \cong E_p$ that fixes the starting point and the asymptotic direction at infinity.

\begin{defn*}[visible direction]
\label{241014173703}
Suppose $(X, \pi)$ is an endless Riemann surface and $p \in X$ is a regular point.
A \dfn{visible direction} at $p$ is either a regular ray $\alpha \in S_p X$ or a Stokes ray $\alpha \in S_p X$ together with a choice of a letter $\varepsilon \in \set{\rm{L}, \rm{R}}$.
In other words, a \textit{visible direction} at $p$ is a pair $(\alpha, \varepsilon)$ where $\alpha \in S_p X$ is a ray and $\varepsilon$ is either the empty word (if $\alpha$ is a regular ray) or the $1$-letter word $\rm{L}$ or $\rm{R}$ (if $\alpha$ is a Stokes ray).
\end{defn*}

\paragraph{Observable singularities.}
Suppose $(X, \pi)$ is endless, $p \in X$ is a regular point with projection $\pi (p)$ translated to the origin, and $\alpha \in S_p X$ is a Stokes ray with visible singularity $\frak{q} \in \sfop{Sing}_{p,\alpha}^0 (X, \pi)$ over the visible singular value $\xi \in \RR_\alpha$.
Since $\frak{q}$ is isolated, we can find a domain $U \in \frak{q}$ around it, separated from all other singularities, such that $\pi$ restricts to a covering map $U \to \UU^\ast$ over a punctured neighbourhood of $\xi \in \CC$.
Now, a line segment contained in $\UU$ that passes through $\xi$ does not admit a regular lift to $U$, but it can always be distorted either to the left or to the right by adding a small semi-circular diversion around $\xi$.
These distorted line segments have regular lifts to $U$ which start at the same point and terminate at two distinct preimages $p_\rm{L}, p_\rm{R}$ of the same point in $\CC$.

Pick one of these preimages, say $p_\rm{L} \in X$, and repeat the analysis by replacing $p$ with $p_\rm{L}$.
The same phase $\alpha \in \RR /2\pi\ZZ$ now determines a ray $\alpha \in S_{p_\rm{L}} X$ at $p_\rm{L}$.
This ray is either regular or a Stokes ray.
If it is a Stokes ray, let $\frak{q}_\rm{L} \in \sfop{Sing}_{p_\rm{L}, \alpha} (X, \pi)$ be the corresponding visible singularity.
Note that $\frak{q}_\rm{L}$ is not visible from $p$ because it is obscured by the singularity $\frak{q}$ (i.e., the geodesic starting at $p$ with phase $\alpha$ does not limit to $\frak{q}_\rm{L}$).
Nevertheless, we can reach $\frak{q}_\rm{L}$ from $p$ along a real curve that is almost a geodesic with phase $\alpha$; namely, we can travel from $p$ along the geodesic with phase $\alpha$ until we reach the neighbourhood of the singular point $\frak{q}$, then travelling along a leftward detour which is the lift of a left or right semi-circular diversion around $\xi$, and then proceeding to travel along the geodesic with phase $\alpha$ starting at $p_1$.
We may therefore refer to the singularity $\frak{q}_\rm{L}$ as being \textit{observable} from $p$, and the real curve that takes us from $p$ to $\frak{q}_\rm{L}$ is fully determined (up to homotopy) by the phase $\alpha$ and a letter $\rm{L}$ or $\rm{R}$ corresponding to the purely combinatorial information of whether the singularity $\frak{q}$ was avoided on the left or the right.

Now, $\frak{q}_\rm{L}$ is again an isolated singularity, so it can be circumvented on the left or on the right.
This yields another pair of distinct points in $X$ that we may call $p_{\rm{LL}}$ and $p_{\rm{LR}}$ lying over the same point in $\CC$, and the analysis repeats.
After a finite number of iterations of this analysis, we obtain a word $\varepsilon$ on letters $\set{\rm{L}, \rm{R}}$ of finite length, written from left to right, and a regular point $p_\varepsilon$ which has phase $\alpha$ relative to $p$.
These observations motivate the following definitions.

\begin{defn*}[observable singularity]
\label{241014140525}
Suppose $(X, \pi)$ is an endless Riemann surface, $p \in X$ is a regular point, and $\alpha \in S_p X$ is a ray.
A singularity $\frak{q} \in \sfop{Sing} (X,\pi)$ is called a \dfn{singularity observable from $p$} if there is a \dfn{distorted geodesic} with phase $\alpha$ starting at $p$ that limits to $\frak{q}$; i.e., a real curve $\gamma : I \to X$ obtained as the regular lift of a line segment $\wp : I \to \CC$ which has been distorted to the left or to the right at a finite subset of points.
Similarly, an \dfn{infinite distorted geodesic} is an infinite path $\gamma : I \to X$ obtained as the regular lift of an infinite path $\ell_\alpha : [0,1) \to \CC$ which runs along the halfline $\RR_\alpha$ while circumventing (on the left or the right) a discrete subset of points on $\RR_\alpha$.
\end{defn*}

Note that the phase of $\frak{q}$ relative to $p$ is necessarily equal to $\alpha$, so $\frak{q} \in \sfop{Sing}_{p,\alpha} (X,\pi)$.
The subset of all singularities observable from $p$ in the direction $\alpha$ is denoted by $\sfop{Sing}_{p,\alpha}^+ (X,\pi) \subset \sfop{Sing}_{p,\alpha} (X,\pi)$.
Given a (finite or infinite) distorted geodesic $\gamma : I \to X$, the sequence of left and right diversions forms a (possibly infinite) binary word $\varepsilon$, which we call the \dfn{signature} of $\gamma$, comprised of letters $\set{\rm{\rm{L}, \rm{R}}}$.
The homotopy class of $\gamma$ is completely determined by the phase $\alpha$ and the signature $\varepsilon$.

\begin{defn*}[observable direction]
\label{241014141058}
An \dfn{observable direction} is then a pair $(\alpha, \varepsilon)$ consisting of a phase $\alpha$ and a binary word $\varepsilon$ on letters $\set{\rm{L}, \rm{R}}$.
If $\varepsilon$ has infinite length, we also require that it stabilises at either $\rm{L}$ or $\rm{R}$.
An observable direction $(\alpha, \varepsilon)$ at $p \in X$ is \dfn{stable}, \dfn{semi-stable}, or \dfn{unstable} if the ray $\alpha \in S_{p_\varepsilon} X$ at $p_\varepsilon$ is respectively, stable, semi-stable, or unstable.
\end{defn*}

\subsection{Borel Resummation and Resurgent Series}

In this subsection, we give a quick review of Borel resummation and introduce the notion of resurgent series.

\begin{defn}[Borel transform and factorial type]
Given a formal power series,
\begin{eqntag}
\label{240806083555}
	\hat{f} (\hbar) = \sum_{k = 0}^\infty a_k \hbar^k \in \CC \bbrac{\hbar}
\fullstop{,}
\end{eqntag}
its \dfn{Borel transform} is by definition the following formal power series in $t$:
\begin{eqntag}
\label{240806222257}
	\hat{\phi} (t) 
		\coleq \Borel \big[ \ \hat{f} \ \big] (t)
		\coleq \sum_{k=0}^\infty \frac{a_{k+1}}{(k+1)!} t^k
		\in \CC \bbrac{t}
\fullstop
\end{eqntag}
The series $\hat{\phi}$ is convergent if and only if the series $\hat{f}$ is of \dfn{factorial type} (or is a \dfn{Gevrey-1 series}), which means its coefficients have at most factorial growth: there are real constants $\C, \M > 0$ such that $|a_k| \leq \C \M^k k!$ for all $k$.
\end{defn}

\begin{defn}[Borel summability]
\label{241011190046}
A formal power series $\hat{f} \in \CC \bbrac{\hbar}$ of factorial type is called \dfn{Borel summable} \textit{with phase $\alpha \in \RR / 2\pi\ZZ$} (or \textit{in the direction $\alpha$}) if its Borel transform $\hat{\phi} (t) \in \CC \cbrac{t}$ admits analytic continuation $\phi \in \cal{O} (\SS)$ to a neighbourhood $\SS \subset \CC$ of the halfline $\RR_\alpha$ with \dfn{exponential growth at $\infty$}; i.e., there are constants $\C, \K > 0$ such that
\begin{eqntag}
\label{240806222248}
	\big| \phi (t) \big| \leq \C e^{\K |t|}
\qqquad
	\forall t \in \SS
\fullstop
\end{eqntag}
This exponential bound \eqref{240806222248} implies that the \dfn{Laplace transform} of $\phi$ \textit{with phase $\alpha$},
\begin{eqntag}
\label{240807083858}
	\Laplace_\alpha [ \phi ] (\hbar) \coleq \int_{\RR_\alpha} \phi (t) e^{-t/\hbar} \d{t}
\fullstop{,}
\end{eqntag}
defines a holomorphic function in a sector at the origin with aperture $\sfop{Arc} (\alpha)$,
\begin{eqntag}
\label{241020130805}
	V_\alpha \coleq \set{ \hbar ~\big|~ \Re ( e^{i\alpha} / \hbar ) > \K }
\end{eqntag}
In this case, the \dfn{Borel resummation} of $\hat{f}$ \textit{with phase $\alpha$} is the holomorphic function $f_\alpha \in \cal{O} (V_\alpha)$ defined by
\begin{eqntag}
\label{240806222612}
	f_\alpha (\hbar) = s_\alpha [ \ \hat{f} \ ] (\hbar) \coleq a_0 + \Laplace_\alpha [ \phi ] (\hbar)
\fullstop
\end{eqntag}
Furthermore, $\hat{f}$ is the asymptotic expansion of $f_\alpha$ in the sector $V_\alpha$:
\begin{eqntag}
\label{241008092632}
	f_\alpha (\hbar) \sim \hat{f} (\hbar)
\qquad
\text{as $\hbar \to 0$ along $\sfop{Arc} (\alpha)$}
\fullstop
\end{eqntag}
\end{defn}

\begin{defn}[lateral Borel summability]
\label{241020124230}
Suppose $\RR_\alpha^\rm{L/R}$ is an infinite path in $\CC$ which runs along the halfline $\RR_\alpha$ while circumventing, always on the left or the right respectively, a discrete subset of points on $\RR_\alpha$.
Then a formal series $\hat{f} \in \CC \bbrac{\hbar}$ of factorial type is called respectively \dfn{left} or \dfn{right laterally Borel summable} \textit{in the direction $\alpha$} if its Borel transform $\hat{\phi} (t) \in \CC \cbrac{t}$ admits analytic continuation $\phi \in \cal{O} (\SS)$ to a neighboourhood $\SS \subset \CC$ of the laterally distorted halfline $\RR_\alpha^\rm{L/R}$ with exponential growth at $\infty$; i.e., satisfying \eqref{240806222248}.
The \dfn{left} or \dfn{right lateral Laplace transform} \textit{of $\phi$ with phase $\alpha$},
\begin{eqntag}
\label{241020135022}
	\Laplace_\alpha^\rm{L/R} [ \phi ] (\hbar) \coleq \int_{\RR_\alpha^\rm{L/R}} \phi (t) e^{-t/\hbar} \d{t}
\fullstop{,}
\end{eqntag}
also defines a holomorphic function in a sector $V_\alpha$ of the form \eqref{241020130805} but with possibly a bigger constant $\K$ than in the exponential bound \eqref{240806222248}.
In this case, the \dfn{left} or \dfn{right lateral Borel resummation} of $\hat{f}$ \textit{with phase $\alpha$} is respectively the holomorphic function $f_\alpha^\rm{L/R} \in \cal{O} (V_\alpha)$ defined by
\begin{eqntag}
\label{241020135027}
	f_\alpha^\rm{L/R} (\hbar) = s_\alpha^\rm{L/R} [ \ \hat{f} \ ] (\hbar) \coleq a_0 + \Laplace_\alpha^\rm{L/R} [ \phi ] (\hbar)
\fullstop
\end{eqntag}
Furthermore, $\hat{f}$ is the asymptotic expansion of $f_\alpha$ in the sector $V_\alpha$:
\begin{eqntag}
\label{241020135030}
	f_\alpha^\rm{L/R} (\hbar) \sim \hat{f} (\hbar)
\qquad
\text{as $\hbar \to 0$ along $\sfop{Arc} (\alpha)$}
\fullstop
\end{eqntag}
\end{defn}

\paragraf
It turns out that the shape of the domain $\SS$ plays an important role for Borel resummation (cf. \autoref{241011174325}), so we introduce the following terminology.

\begin{defn*}[stability of Borel resummation]
\label{241011190023}
If $\hat{f} \in \CC \bbrac{\hbar}$ is Borel summable with phase $\alpha$, let $\SS \subset \CC$ be a neighbourhood of the halfline $\RR_\alpha$ where the exponential bound \eqref{240806222248} holds.
Then we say that the Borel resummation $f_\alpha$ of $\hat{f}$ is
\begin{enumerate}
\item \dfn{stable} if $\SS$ can be chosen to be an infinite sector $\SS = \SS_A$;
\item \dfn{semi-stable} if $\SS$ can be chosen to be a halfstrip $\SS = \SS_\alpha$;
\item \dfn{unstable} otherwise.
\end{enumerate}
In these instances, we say that $\hat{f}$ is \dfn{stably}, \dfn{semi-stably}, or \dfn{unstably Borel summable} with phase $\alpha$, respectively.
We also say that $\hat{\phi} (t)$ has \dfn{stable}, \dfn{semi-stable}, or \dfn{unstable exponential growth at $\infty$}, respectively.
Clearly, stable implies semi-stable but not vice versa, whereas unstable and semi-stable are mutually exclusive; we often write ``(semi-)stable'' to mean ``at least semi-stable''.
\end{defn*}

\paragraf
The significance of stability is the ability to glue Borel resummations for all nearby phases into a single holomorphic function, as explained next.

\begin{defn*}[Borel resummation along an arc]
\label{241011191045}
Suppose $A \subset \RR / 2\pi\ZZ$ is an arc of length $|A| \leq \pi$.
A formal power series $\hat{f} \in \CC \bbrac{\hbar}$ of factorial type is called \dfn{Borel summable} \textit{along $A$} if it is stably Borel summable in every direction $\alpha \in A$.
\end{defn*}

In this case, all the Laplace transforms $\Laplace_\alpha [\phi]$ for different phases $\alpha \in A$ glue together into a single holomorphic function $\Laplace_{A} [\phi]$ on a sector $V_A$ at $\hbar = 0$ with aperture $\sfop{Arc} (A)$; namely,
\begin{eqntag}
\label{241009115236}
	\Laplace_{A} [\phi] (\hbar) \coleq \Laplace_\alpha [\phi] (\hbar)
\qqtext{for all $\hbar \in V_A$ with $\arg (\hbar) = \alpha \in A$.}
\end{eqntag}
Consequently, the Borel resummations of $\hat{f}$ for different phases $\alpha \in A$ glue into a single holomorphic function $f_A \in \cal{O} (V_A)$:
\begin{eqntag}
\label{241009120208}
	f_A (\hbar)
		\coleq s_A [ \hat{f} ] (\hbar)
		\coleq a_0 + \Laplace_A [ \hat{f} ] (\hbar)
\fullstop
\end{eqntag}
Furthermore, $\hat{f}$ is the asymptotic expansion of $f_A$ in the sector $V_A$:
\begin{eqntag}
\label{241008092633}
	f_A (\hbar) \sim \hat{f}_A (\hbar)
\qquad
\text{as $\hbar \to 0$ along $\sfop{Arc} (A)$}
\fullstop
\end{eqntag}

\begin{rem}[significance of semi-stability]
\label{241011174325}
A fundamentally important Theorem of Watson-Nevanlinna-Sokal \cite{zbMATH02629428,nevanlinna1918theorie,MR558468}\footnote{See also \cite[Théorème 1.4.1.1]{zbMATH00797135} or \cite[Theorem B.15]{MY2008.06492} for a more contemporary exposition.} asserts that if the Borel resummation $f_\alpha$ is semi-stable, then $f_\alpha$ is the unique sectorial germ along $\sfop{Arc} (\alpha)$ such that the asymptotic expansion \eqref{241008092632} holds with uniform factorial growth.
In other words, $f_\alpha$ is the canonical (i.e., coordinate-invariant) lift of the formal power series $\hat{f}$ to a sectorial germ along the semicircular arc $\sfop{Arc} (\alpha)$ bisected by $\alpha$.

This canonical relationship between $f_\alpha$ and $\hat{f}$ breaks down if the Borel resummation $f_\alpha$ is unstable; i.e., if we ask that the exponential bound \eqref{240806222248} holds in an arbitrary open neighbourhood $\SS$ of the halfline $\RR_\alpha$ that is not necessarily a halfstrip.
Namely, although the Laplace integral \eqref{240807083858} still defines a holomorphic function in the sector $V_\alpha$, and although the asymptotic equivalence \eqref{241008092632} still holds, it is no longer possible to conclude in general that $f_\alpha$ is the canonical lift of $\hat{f}$.
Indeed, in this case there are (uncountably) infinitely many functions which satisfy \eqref{241008092632}.
\end{rem}

\begin{defn}[exponential type]
\label{241008162601}
Suppose a holomorphic germ $\hat{\phi} \in \CC \cbrac{t}$ admits endless analytic continuation.
We say that $\hat{\phi}$ has \dfn{exponential type at $\infty$} (resp. \textit{along an arc $A$}) if it has exponential growth at $\infty$ in every observable direction (resp. with phase in $A$) with matching stability.
Concretely, we ask for the following properties to hold:
\begin{enumerate}
\item If $\alpha$ is a stable, semi-stable, or unstable regular ray, then $\hat{\phi}$ has respectively a stable, semi-stable, or unstable exponential growth along $\RR_\alpha$.
\item If $\alpha$ is a stable, semi-stable, or unstable Stokes ray, then $\hat{\phi}$ has respectively a stable, semi-stable, or unstable exponential growth along the diverted halflines $\RR_\alpha^\rm{L}$ and $\RR_\alpha^\rm{R}$.
\item If $(\alpha, \varepsilon)$ is an observable direction at $p$, and $\gamma_\varepsilon$ is the path from $p$ to $p_\varepsilon$, then the analytically continued germ $\gamma_\varepsilon \cdot \hat{\phi}$ at $p_\varepsilon$ has exponential type at $\infty$ in the direction $\alpha$.
\end{enumerate}
\end{defn}

\begin{defn}
\label{241005184045}
A \dfn{resurgent series of log-algebraic type} is a formal power series $\hat{f} (\hbar) \in \CC \bbrac{\hbar}$ of factorial type whose Borel transform $\hat{\phi} (t) \coleq \Borel [ \, \hat{f} \, ] (t) \in \CC \set{t}$ admits an endless analytic continuation of log-algebraic type and exponential type at $\infty$ in all directions observable from the origin.

If the directions are restricted to an arc $A$, then we say $\hat{f} (\hbar)$ is \dfn{resurgent along $A$}.
If the word ``observable'' is replaced with ``visible'', then we refer to $\hat{f} (\hbar)$ as \dfn{weakly resurgent}.
On the other hand, if ``observable from the origin'' is replaced with ``visible from every regular point of $X$'', where $X$ is the endless Riemann surface to which the Borel transform is analytically continued, then we refer to $\hat{f} (\hbar)$ as \dfn{strongly resurgent}.
\end{defn}

\begin{rem}[terminology]
\label{240730092212}
Let us address some of the naming discrepancies found in the literature.
\autoref{240913161447} of \textit{endless analytic continuation} is our interpretation of the notion by the same name (or ``\textit{prolongeabilité sans fin}'' in French) introduced by Candelpergher-Nosmas-Pham in \cite[Définition (iii), p. 204]{CNP1993}; it is possible to show that our formulation is equivalent to theirs.
Elsewhere, the name \textit{endless analytic continuation} is used for any of the special cases that we named \textit{discrete type} and \textit{finite type}; see, e.g., \cite[Chapter 6]{MR3526111} and \cite[p.177]{MR3495546}.

Our \autoref{241005184045} of \textit{resurgent series} (log-algebraic type property aside) is more restrictive than usual definitions found in the literature that do not explicitly assume exponential growth at $\infty$, with notable exceptions such as \cite{MR2116813,MR2427514}.
However, we believe that our definition actually captures more the sense in which resurgence is used in modern applications (especially in the physics literature).
There, at least some form of exponential growth is assumed, often implicitly; see, for example, \cite{MR3955133,Alim2024} and references therein.
Our notions of \textit{weakly} and \textit{strongly resurgent series} are closer to (but not the same as) what is sometimes called \textit{summable and resurgent series} and \textit{summable-resurgent series}, respectively; see \cite[Definition 2.5]{MR3495546} (see also \cite[§5.3.4]{MR2737219}) but note that they also assume having only finitely many Borel singularities.

In contrast, if one explicitly makes no assumptions on the exponential growth (as in, e.g., \cite{MR0819354,MR1250603,MR1397029}), then the Laplace transform is not well-defined, and so the Borel resummation as well as the associated Stokes phenomenon have no direct meaning in this generality.
The standard approach seems to be to consider the \textit{truncated} Laplace transform by integrating along a ray only up to a finite cutoff value, and speak therefore of Borel \textit{pre}summation.
The difficulty with this approach is that the Borel presummation of a series depends on the cutoff and is therefore only well-defined up to the set of all possible subexponential corrections.
Correspondingly, endlessly continuable germs in the $t$-plane are considered only up to all possible entire functions and one speaks of \textit{microfunctions}.
In this language, a \textit{resurgent function} is essentially an endlessly continuable microfunction, which is closer to the original point of view of Écalle \cite{zbMATH03971144}.
Yet Écalle's more recent texts \textit{do} assume exponential growth explicitly; cf. \cite[§2.1]{MR4417677}.
\end{rem}

\subsection{Resurgent Perturbation Theory}

In this subsection, we formulate the ideas of the previous subsections in the setting of perturbation theory.

\paragraf
The fundamental goal of the WKB analysis is to solve the Schrödinger equation \eqref{240711180437} using perturbation theory in the parameter $\hbar$.
Specifically, WKB analysis is a method for obtaining linearly independent solutions with prescribed asymptotics as $\hbar \to 0$.
The suitable notion of asymptotics for this purpose is sometimes known as \textit{exponential asymptotics}.
Namely, an \dfn{exponential perturbative series} on a domain $U \subset \CC$ is an expression of the form
\begin{eqntag}
\label{240813085456}
	\hat{\Psi} (x, \hbar) = e^{-\S(x) / \hbar} \hat{\A} (x, \hbar)
\qtext{with}
	\hat{\A} (x, \hbar) = \sum_{k=0}^\infty \A_k (x) \hbar^k
\fullstop{,}
\end{eqntag}
where $\S \in \cal{O} (U)$ is a holomorphic function on $U$ (sometimes called \textit{phase function} or \textit{prepotential}), and $\hat{\A} \in \cal{O} (U) \bbrac{\hbar}$ is a \dfn{perturbative series} on $U$ (sometimes called \textit{amplitude}).
A \dfn{perturbative transseries} on $U$ is a formal (possibly countable) $\CC$-linear combination of exponential perturbative series.

Correspondingly, we must extend the notions from the previous sections to allow them to vary with respect to $x$, a point on some Riemann surface.
One way to do this is as follows.

\begin{defn}
\label{241015104954}
Let $S$ be a Riemann surface.
An \dfn{$S$-family of Riemann surfaces spread over $\CC$} is the data $(X, \pi, s)$ consisting of a $2$-dimensional holomorphic manifold $X$ equipped with a holomorphic submersion $s : X \to S$ and a holomorphic map $\pi : X \to \CC$ such that, for every $p \in S$, the pair $(X_p, \pi_p)$, consisting of the fibre $X_p \coleq s^\inv (p)$ and the restriction $\pi_p \coleq \pi |_{X_p} : X_p \to \CC$, defines a Riemann surface spread over $\CC$.
Moreover, $(X, \pi, s)$ is \dfn{endless of log-algebraic type} if each $(X_p, \pi_p)$ is endless of log-algebraic type.
\end{defn}

\begin{defn}
\label{241015104211}
Suppose $S$ is a Riemann surface and $\hat{f} \in \cal{O}_S \bbrac{t}$ is a formal power series with coefficients which are local holomorphic functions on $S$.
We say that $\hat{f}$ is a \dfn{uniformly resurgent series} \textit{of log-algebraic type} on a domain $U \subset S$ if the following properties hold uniformly for all $p \in U$.
The formal power series $\hat{f} (p, \hbar) \in \CC \bbrac{\hbar}$ is of factorial type and its Borel transform $\hat{\phi} (p, t) \in \CC \set{t}$ admits an analytic continuation to a $U$-family of endless Riemann surfaces of log-algebraic type and exponential type at $\infty$ in all directions observable from each origin $(p,0) \in U \times \CC$.
\end{defn}

\section{The Geometry of Spectral Curves}
\label{240919123411}

In this section, describe the relevant geometry of the spectral curve of a meromorphic quadratic differential.
A Schrödinger operator, as we will define in \autoref{241022163033}, has a naturally associated quadratic differential, and the geometry it gives rise to is central to the resurgence of WKB solutions.

\subsection{Quadratic Differentials}
\label{240730094616}

\paragraf
Let $(X,D)$ be a \dfn{marked curve}; i.e., a compact Riemann surface $X$ equipped with an effective divisor $D$.
Let $\omega_X$ be the holomorphic cotangent bundle of $X$ and set
\begin{eqn}
	\omega_{X,D} \coleq \omega_X (D)	
\fullstop
\end{eqn}
This is the line bundle of meromorphic differential forms on $X$ with poles bounded by $D$, which we think of as the \dfn{cotangent bundle of the marked curve} $(X,D)$.
We denote its total space by $T^\ast_{X,D} \coleq \sf{tot} (\omega_{X,D})$.
Note that we allow $D$ to have higher multiplicity or be empty, but we always assume that $\deg (D) \geq 2 - 2g$ where $g$ is the genus of $X$.

Concretely, suppose $p \in D$ has multiplicity $k$ and let $x$ be a local coordinate centred at $p$.
The main difference between the cotangent bundles of $X$ and $(X,D)$ is that $\omega_X$ is generated near $p$ by $\d{x}$ whilst $\omega_{X,D}$ is generated by $x^{-k} \d{x}$.
Any local holomorphic section $f(x) \d{x}$ of $\omega_X$ can be viewed as a local holomorphic section of $\omega_{X,D}$ by writing $f(x) \d{x} = f(x) x^{k} (x^{-k} \d{x})$, but not necessarily vice versa.

\begin{defn}
\label{240813120707}
A \dfn{quadratic differential} $\phi$ on a marked curve $(X,D)$ is by definition a global holomorphic section of the line bundle $\omega^{2}_{X,D} \coleq \omega_{X,D} \otimes \omega_{X,D} = \omega^{2}_X (2D)$:
\begin{eqntag}
\label{241016144548}
	\phi \in H^0 (X, \omega_{X,D}^{2})
\fullstop
\end{eqntag}
\end{defn}

Note that according to this definition a holomorphic quadratic differential on $X$ is in particular a quadratic differential on $(X,D)$ for any $D$.
We say the divisor $D$ is \dfn{minimal} \textit{relative to} $\phi$ if $\phi$ is a holomorphic section of $\omega_{X,D}^{2}$ but not a holomorphic section of $\omega_{X,D'}^{2}$ for any proper subdivisor $D' < D$.
It follows that if $D$ is minimal for $\phi$, then $\phi$ has at most a simple zero at every point of $D$.
We also say that the divisor $D$ is \dfn{full} relative to $\phi$ if it is minimal and $\phi$ is nonvanishing along $D$.
If $D$ is minimal for $\phi$, we denote the vanishing divisor of $\phi$ by
\begin{eqn}
	B \coleq \Z ( \phi )
\fullstop
\end{eqn}
Thus, $D$ is minimal if $B \cap D$ has no higher multiplicity, and full if $B \cap D = \emptyset$.
If $D$ is minimal, points in the support of $D$ and $B \smallsetminus D$ are called \dfn{poles} and \dfn{zeros} respectively.

\begin{defn}[simple quadratic differential]
\label{241017182312}
A quadratic differential $\phi$ on $(X,D)$ is called \dfn{simple} if $D$ is minimal relative to $\phi$ and all zeros of $\phi$ are simple.
A simple quadratic differential is called \dfn{complete} if it has no simple poles.
\end{defn}

\paragraf
Concretely, if $p \in D$ has multiplicity $k$ and $x$ is a local coordinate centred at $p$, then the line bundle $\omega^{2}_{X,D}$ is generated near $p$ by $(x^{-k} \d{x})^2 = x^{-2k} \d{x}^2$, so locally near $p$ any quadratic differential on $(X,D)$ is of the form $\phi \loceq f(x) x^{-2k} \d{x}^2$.
The divisor $D$ is minimal if $f$ vanishes at $p$ to order at most $1$, and full if $f$ is nonvanishing at $p$.

\paragraf
Equivalently, a quadratic differential $\phi$ on $(X,D)$ can be viewed as a meromorphic section of the line bundle $\omega^2_X$ given locally by
\begin{eqntag}
\label{240726122413}
	\phi \loceq \Q(x) \d{x}^{2}
\fullstop{,}
\end{eqntag}
where $q$ is some meromorphic function with poles at the points of $D$ and zeros at the points of the divisor $B \smallsetminus D \coleq B - (B \cap D)$.
Although, a quadratic differential $\phi$, as a section of $\omega_{X,D}^{2}$, is everywhere regular, we still refer to points in $D$ as \dfn{poles} of $\phi$.
The \dfn{order} of a pole at $p$ refers to the pole order of $q$; i.e., the pole order as a section of $\omega^2_X$.
Similarly, although $\phi \in \omega_{X,D}^{2}$ has a zero at every point in $B$, including in the intersection $B \cap D$, we distinguish between \dfn{true zeros} $B \smallsetminus D$ and \dfn{virtual zeros} $B \cap D$.
Virtual zeros correspond to odd-order poles.
If $D$ is full, there are no poles of odd order and hence no virtual zeros.

If $D$ is minimal for $\phi$, all virtual zeros of $\phi$ are simple.
In this case, if $\phi$, viewed as a section of $\omega^{2}_X$, has a pole at $p \in X$ of odd order, then when viewed instead as a section of $\omega_{X,D}^{2}$ it has a simple zero at $p$.
Concretely, if $p \in \sfop{Pol} (\phi)$ is a pole of order $m$ and $x$ is a local coordinate centred at $p$, then $\phi \in \omega^{2}_X$ has a local expression $f(x) x^{-m} \d{x}^2$ (for some function $f (x)$ holomorphic and nonzero at $x = 0$), whilst $\phi \in \omega^{2}_{X,D}$ has the local expression $f(x) \big(x^{-k} \d{x} \big)^2$ if $m = 2k$ or $x f(x) \big(x^{-k} \d{x} \big)^2$ if $m = 2k - 1$.
Thus, when $m$ is odd, the virtual zero at $p$ is only visible when passing to the bundle $\omega^2_{X,D}$.

\paragraf
If $\phi$ is a quadratic differential on $(X,D)$ with minimal $D$, then we denote by $\sf{Pol} (\phi)$ and $\sf{Zer} (\phi)$ respectively the subsets of poles and true zeros of $\phi$.
It is sometimes important to distinguish between poles of even and odd orders, so let $\sf{Pol}_\text{ev} (\phi)$ and $\sf{Pol}_\text{od} (\phi)$ denote their respective subsets.
We also write $\sf{Pol}_1 (\phi)$ for the subset of simple poles: they require special treatment.
The relationship between these subsets and the divisors $B$ and $D$ is summarised next:
\begin{eqn}
\begin{gathered}
	\sfop{Supp} (D) = \sfop{Pol} (\phi),
\qquad
	\sfop{Supp} (B) = \sfop{Zer} (\phi) \cup \sfop{Pol}_\textup{od} (\phi)
\fullstop{,}
\\
	\sfop{Supp} (B \smallsetminus D) = \sfop{Zer} (\phi),
\quad
	\sfop{Supp} (B \cap D) = \sfop{Pol}_\textup{od} (\phi)
\fullstop{,}
\\
	\sfop{Supp} (D \smallsetminus B) 
		= \sfop{Pol} (\phi) \smallsetminus \sfop{Pol}_1 (\phi)
		\eqcol \sfop{Pol}_{\geq 2} (\phi)
\fullstop
\end{gathered}
\end{eqn}
If $m_p > 0$ is the pole order at $p \in \sf{Pol} (\phi)$ or the vanishing order at $p \in \sf{Zer} (\phi)$, then
\begin{eqntag}
\label{240728162258}
	D = \sum_{p \in \sf{Pol}_\text{ev} (\phi)} \tfrac{m_p}{2} \cdot p + \sum_{p \in \sf{Pol}_\text{od} (\phi)} \tfrac{m_p + 1}{2} \cdot p 
\fullstop{,}
\end{eqntag}
\begin{eqntag}
\label{240728171830}
	B = \sum_{p \in \sfop{Zer} (\phi)} m_p \cdot p + \sum_{p \in \sfop{Pol}_\text{od} (\phi)} 1 \cdot p
\fullstop
\end{eqntag}
The zeros and poles of $\phi$ are collectively called \dfn{critical points}, and all other points are called \dfn{regular points}.
Poles of order $2$ or greater are often called \dfn{infinite critical points}, whilst zeros and simple poles are called \dfn{finite critical points}.
In WKB analysis, true zeros and simple poles of $\phi$ are usually called \dfn{turning points}; a \dfn{simple turning point} is a true simple zero.
We denote the subsets of turning points and simple turning points by
\begin{eqntag}
\begin{aligned}
\label{241016182934}
	\sfop{Turn} (\phi) 
		&\coleq \sfop{Zer} (\phi) \cup \sfop{Pol}_1 (\phi)
		= \sfop{Supp} \big( B \smallsetminus (D \smallsetminus B) \big)
\fullstop{,}
\\
	\sfop{Turn}_0 (\phi) 
		&\coleq \sfop{Zer} (\phi)
		= \sfop{Supp} \big( B \smallsetminus D \big)
\fullstop
\end{aligned}
\end{eqntag}
To see that these equalities are true, recall how to work with divisors: locally around a single point, if $B = m \cdot p$ and $D = n \cdot p$, then $D - B = (n-m)\cdot p$, $B \cap D = \min \set{n,m} \cdot p$, and so $D \smallsetminus B = (n - \min \set{n,m})\cdot p$, and therefore $B \smallsetminus (D \smallsetminus B) = k \cdot p$ where $k = m - \min \set{m, n - \min \set{n,m}}$.
If $p \in \sfop{Zer} (\phi)$, then $m \geq 1$ and $n = 0$, and so $k = m$.
If $p \in \sfop{Pol}_\textup{od} (\phi)$ of order $3$ or greater, then $m = 1$ and $n \geq 2$, and so $k = 0$.
If $p \in \sfop{Pol}_1 (\phi)$, then $m = 1$ and $n = 1$, and so $k = 1$.

\paragraf
Explicitly, we record the following local normal form result for quadratic differential which is straightforward to verify.
In each case, we express the local normal form for $\phi$ both thought of as a meromorphic section of $\omega^2_X$ and as a holomorphic section of $\omega^2_{X,D}$.

\begin{prop*}[local normal forms]
\label{241017104429}
Let $\phi$ be a quadratic differential on $(X,D)$ with minimal divisor $D$.
Pick $p \in X$ and let $x$ be a local coordinate centred at $p$.
If $p$ is a regular point, i.e., $p \in X \smallsetminus (B \cup D)$, then $x$ can be chosen such that $\phi (x) = \d{x}^2$.
If $p \in \sfop{Zer} (\phi)$ is a zero of order $m \geq 1$, then $x$ can be chosen such that 
\begin{eqntag}
\label{241017110704}
	\phi (x) = x^m \d{x}^2
\fullstop
\end{eqntag}
If $p \in \sfop{Pol}_\textup{od} (\phi)$ has order $m = 2k-1 \geq 1$, then $x$ can be chosen such that 
\begin{eqntag}
\label{241017105138}
	\phi (x) = x^{-(2k-1)} \d{x}^2 = x \left( \tfrac{\d{x}}{x^k} \right)^2
\fullstop
\end{eqntag}
If $p \in \sfop{Pol}_\textup{ev} (\phi)$ has order $m = 2k \geq 2$, then $x$ can be chosen such that
\begin{eqntag}
\label{241017110836}
	\phi (x) =
\begin{cases}
	r x^{-2} \d{x}^2 = r \left( \tfrac{\d{x}}{x^k} \right)^2
		&\qquad \text{if $m = 2$;}
\\	(x^{-k} + b x^{-1})^2 \d{x}^2 = \big( 1 + b x^{k-1} \big)^2 \left( \tfrac{\d{x}}{x^k} \right)^2
		&\qquad \text{if $m = 2k \geq 4$,}
\end{cases}
\end{eqntag}
for some $r \in \CC^\times$ and $b \in \CC$.
In fact, $r$ and $b^2$ are coordinate-invariant quantities called the \dfn{quadratic residue} of $\phi$ at $p$, denoted by $\Res^2_p (\phi)$.
Thus $\Res^2_p (\phi) = r$ if $m = 2$; $\Res^2_p (\phi) = b^2$ if $m = 2k \geq 4$; and $\Res^2_p (\phi) = 0$ if $m$ is odd.
Finally, we note that formula \eqref{241017110836} can be written a little more uniformly: if $m$ is even, write $m = 2k$ for $k \geq 1$, then $x$ can be chosen such that
\begin{eqntag}
\label{241017105802}
	\phi (x) = \big( 1 - k + \sqrt{r} x^{k-1} \big)^2 \left( \tfrac{\d{x}}{x^k} \right)^2
\fullstop{,}
\qqquad
	r = \Res^2_p (\phi)
\fullstop
\end{eqntag}
\end{prop*}

\subsection{The Spectral Curve}
\label{240813110632}

Suppose $\phi$ is a simple quadratic differential on a marked curve $(X,D)$.

\begin{defn}
\label{241017182538}
The \dfn{spectral curve} of $\phi$ is the branched double cover $\pi : \sfSigma \to X$ defined as the following subvariety of the cotangent bundle $T^\ast_{X,D}$:
\begin{eqntag}
\label{241016144535}
	\sfSigma \coleq \set{ \big(p, \lambda (p) \big) \in T^\ast_{X,D} ~\Big|~ \lambda (p) \otimes \lambda (p) = \phi (p) }
\fullstop
\end{eqntag}
\end{defn}

Here, $\phi (p)$ is the value at $p \in X$ of the holomorphic section $\phi \in \omega_{X,D}^{2}$ and $\lambda (p)$ is a point in the fibre $T^\ast_{X,D} |_p$.
The map $\pi : \sfSigma \to X$ sending $\big(p, \lambda (p) \big) \mapsto p$ is just the restriction of the bundle projection $\pi : T^\ast_{X,D} \to X$.
As any double cover, $\sfSigma$ carries a covering involution $\iota : \sfSigma \to \sfSigma$ which sends $\big(p, \lambda (p) \big) \mapsto \big(p, -\lambda (p) \big)$.

\paragraf
Since all zeros of $\phi$ are simple, it follows that $\sfSigma$ is a compact Riemann surface which is branched over the vanishing locus of $\phi \in \omega_{X,D}^{2}$ which is the subset of all zeros and odd-order poles of $\phi$.
In other words, the branch locus of $\pi$ is the divisor $B$ introduced earlier, and we highlight the (sometimes misunderstood) discrepancy between branch points of the spectral curve and turning points: every turning point is a branch point, but not every branch point is necessarily a turning point:
\begin{eqn}
	B = \sfop{Bra} (\pi)
\qqtext{and}
	\sfop{Turn} (\phi) \subset \sfop{Supp} (B)
\fullstop
\end{eqn}
This inclusion is an equality only if there are no poles of odd order $3$ or higher.

\paragraf
Concretely, pick $p \in X$ and let $x$ be a local coordinate centred at $p$.
If $p \in X \smallsetminus D$, then $\phi (x) = \Q(x) \d{x}^2$ and the local equation for the spectral cover $\sfSigma$ is $y^2 = \Q(x)$.
On the other hand, if $p \in D$ with multiplicity $k$, then $\phi (x) = \Q(x) (x^{-k} \d{x})^2$ and the local equation is the equation $y^2 = \Q(x)$.

We can use local normal forms from \autoref{241017104429} to make this even more explicit.
If $p$ is a regular point, choose $x$ such that $\phi (x) = \d{x}^2$.
Then $\sfSigma \loceq \set{y^2 = 1} \subset \CC^2$.

If $p \in \sfop{Zer} (\phi)$ and $x$ is such that $\phi (x) = x^m \d{x}^2$, then $\sfSigma \loceq \set{y^2 = x^m}$.
Notice that the algebraic curve $y^2 = x^m$ has a singularity at $(0,0)$ unless $m = 1$.
This demonstrates explicitly why we consider quadratic differentials with simple zeros only.
The case of higher-order zeros can nevertheless be treated using slightly more sophisticated geometric techniques, which we will not discuss here.

If $p \in \sfop{Pol}_\textup{ev} (\phi)$ is a pole of even order $m = 2k$ and $x$ is such that $\phi$ has local normal form \eqref{241017105802}, then $\sfSigma \loceq \set{y^2 = (1 - k + \sqrt{r} x^{k-1})^2}$.
In particular, $\sfSigma \loceq \set{y^2 = r}$ when $m = 2$.
Finally, if $p \in \sfop{Pol}_\textup{od} (\phi)$ is a pole of even order $m = 2k-1$ and $x$ is such that $\phi (x) = x \left( \frac{\d{x}}{x^k} \right)^2$, then $\sfSigma \loceq \set{y^2 = x}$.

\paragraf
Next, we introduce the ramification divisor of $\pi$ as well as a modification of the pullback of the divisor $D$:
\begin{eqntag}
\label{241017184647}
	R \coleq \sfop{Ram} (\pi)
\qqtext{and}
	\sfDelta \coleq \pi^\ast D \smallsetminus 2R = \pi^\ast (D \smallsetminus B)
\fullstop
\end{eqntag}
We also introduce the following shorthand notation:
\begin{eqntag}
\label{240922122415}
	\sfSigma_R \coleq \sfSigma \smallsetminus R
~,\qquad
	\sfSigma_\sfDelta \coleq \sfSigma \smallsetminus \sfDelta
~,\qquad
	\sfSigma_{R, \sfDelta} \coleq \sfSigma \smallsetminus (R \cup \sfDelta)
\fullstop
\end{eqntag}
Notice that $R$ is supported along the preimages of zeros and odd-order poles of $\phi$, whilst $\sfDelta$ is supported only along the preimages of poles of order $2$ or greater.
In particular, $2R \smallsetminus \pi^\ast D$ is supported only at the preimages of zeros of $\phi$, whilst $\pi^\ast D - 2R$ is supported at all critical points except the preimages of simple poles of $\phi$.
In summary:
\begin{eqntag}
\label{241017184949}
\begin{gathered}
	\sfop{Supp} (R) = \pi^\inv \sfop{Zer} (\phi) \cup \pi^\inv \sfop{Pol}_\textup{od} (\phi)
\fullstop{,}\qquad
	\sfop{Supp} (\sfDelta) = \pi^\inv \sfop{Pol}_{\geq 2} (\phi)
\fullstop{,}
\\
	\sfop{Supp} (2R \smallsetminus \pi^\ast D)
		= \pi^\inv \sfop{Zer} (\phi)
\fullstop{,}
\\
	\sfop{Supp} (\pi^\ast D - 2R)
		= \pi^\inv \sfop{Zer} (\phi) \cup \pi^\inv \sfop{Pol}_{\geq 2} (\phi)
\fullstop
\end{gathered}
\end{eqntag}
Explicitly, the divisors so far are:
\begin{eqnstag}
\label{241017173345}
	\sfDelta &= \sum_{p \in \sf{Pol}_\text{ev} (\phi)} \tfrac{m_p}{2} \cdot ( p^+ + p^-) + \sum_{p \in \sf{Pol}_\text{od} (\phi)} (m_p - 1) \cdot \tilde{p} 
\fullstop{,}
\\
\label{241017173557}
	R &= \sum_{p \in \sf{Zer} (\phi)} 1 \cdot \tilde{p} + \sum_{p \in \sf{Pol}_\text{od} (\phi)} 1 \cdot \tilde{p}
\fullstop{,}
\\
\label{241017221251}
	R_\sfDelta &= \sum_{p \in \sf{Zer} (\phi)} 1 \cdot \tilde{p} + \sum_{p \in \sf{Pol}_1 (\phi)} 1 \cdot \tilde{p}
\fullstop{,}
\\
\label{241018131921}
	R_{\pi^\ast D} &= \sum_{p \in \sf{Zer} (\phi)} 1 \cdot \tilde{p}
\fullstop{,}
\end{eqnstag}
where $m_p > 0$ is the pole order at $p \in \sfop{Pol} (\phi)$,~ $p^+,p^-$ are the two preimages of an even-order pole $p$, and $\tilde{p}$ is the preimage of an odd-order pole or a zero $p$.

In the special case where the divisor $D$ is full relative to the quadratic differential $\phi \in \omega^2_{X,D}$ (i.e., $\phi$ only has even-order poles), we have that $\sfDelta = \pi^\ast D$.

\paragraf
We call (\dfn{simple}) \dfn{transition points} the preimages of all (simple) turning points:
\begin{eqntag}
\label{241018131432}
\begin{aligned}
	\sfop{Tran} (\phi) 
		&\coleq \pi^\inv \sfop{Turn} (\phi)
		= \pi^\inv \sfop{Zer} (\phi) \cup \pi^\inv \sfop{Pol}_1 (\phi)
\fullstop{,}
\\
	\sfop{Tran}_0 (\phi) 
		&\coleq \pi^\inv \sfop{Turn}_0 (\phi)
		= \pi^\inv \sfop{Zer} (\phi)
\fullstop
\end{aligned}
\end{eqntag}
We highlight the discrepancy between transition points and ramification points: every transition point is a ramification point of the spectral curve, but not every ramification point is a transition point: $\sfop{Tran} (\phi) \subset \sfop{Supp} (R)$.
In fact, transition points and simple transition points are nothing but the support of the divisors $R \smallsetminus \sfDelta$ and $R \smallsetminus \pi^\ast D$, respectively:
\begin{eqntag}
\label{241017181721}
\begin{aligned}
	\sfop{Tran} (\phi)
		&= \sfop{Supp} ( R_\sfDelta )
\qqtext{~~where}
	R_\sfDelta \coleq R \smallsetminus \sfDelta
\fullstop
\\
	\sfop{Tran}_0 (\phi)
		&= \sfop{Supp} ( R_{\pi^\ast D} )
\qqtext{where}
	R_{\pi^\ast D} \coleq R \smallsetminus \pi^\ast D
\fullstop
\end{aligned}
\end{eqntag}
To see this, compute locally around a single point $p \in \sfSigma$ with projection $q = \pi (p)$.
If, locally, $R = 1 \cdot p$ and $D = n \cdot q$, then $\pi^\ast D = 2n \cdot p$ and $\pi^\ast D \smallsetminus 2R = (2n - \min \set{2n,2}) \cdot p = 2(n-1) \cdot p$, so $R_\sfDelta = k \cdot p$ where $k = 1 - \min \set{ 1 , 2(n - 1)}$.
If $p \in \pi^\inv \sfop{Zer} (\phi)$, then $m = 1$ and $n = 0$, so $k = 1$.
If $p \in \sfop{Pol}_\textup{od} (\phi)$ of order $3$ or greater, then $m = 1$ and $n \geq 2$, and so $k = 0$.
But if $p \in \sfop{Pol}_1 (\phi)$, then $m = 1$ and $n = 1$, and so $k = 1$.

\paragraf
The spectral cover $\sfSigma$ carries a canonical differential often called the \dfn{Liouville form}:
\begin{eqntag}
\label{241018114416}
	\lambda \in H^0 (\sfSigma, \omega_{\sfSigma, \sfDelta})
\fullstop
\end{eqntag}
It comes via $\sfDelta \subset \pi^\ast D$ from the canonical meromorphic differential form
\begin{eqntag}
\label{241018110723}
	\lambda \in H^0 \big(\sfSigma, \omega_{\sfSigma, \pi^\ast D} \big)
\qqtext{such that}
	\lambda \otimes \lambda = \pi^\ast \phi
\qtext{and}
	\iota^\ast \lambda = -\lambda
\fullstop{,}
\end{eqntag}
which is obtained by restricting to $\sfSigma$ the tautological section $\lambda_\text{taut}$ of the line bundle $\pi^\ast \omega_{X,D}$ on $T^\ast_{X,D}$.
Namely, there is a canonical isomorphism $\pi^\ast \omega_X \iso \omega_\sfSigma (-R)$ induced by the derivative of $\pi$, and $\lambda$ is the image of $\lambda_\text{taut}$ under the composition of maps
\begin{eqn}
	\pi^\ast \omega_{X,D}
		= \pi^\ast \omega_X (\pi^\ast D)
		\iso \omega_{\sfSigma} (\pi^\ast D - R)
		= \omega_{\sfSigma, \pi^\ast D} (-R)
		\to \omega_{\sfSigma, \pi^\ast D}
\fullstop	
\end{eqn}
Notice that this map vanishes along the ramification locus $R$.
At the same time, a square root of $\pi^\ast \phi$ is a holomorphic section of $\pi^\ast \omega_{X,D}$ which has a simple zero at every ramification point.
Consequently, as a section of the line bundle $\omega_{\sfSigma, \pi^\ast D}$, we see that $\lambda$ has a double zero at every ramification point and nowhere else; i.e., $\lambda$ is a global nonvanishing holomorphic section
\begin{eqntag}
\label{241018113803}
	\lambda 
		\in H^0 \big( \sfSigma, \omega_{\sfSigma,\pi^\ast D} (- 2R) \big)
		= H^0 \big( \sfSigma, \omega_{\sfSigma} (\pi^\ast D - 2R) \big)
\fullstop
\end{eqntag}
More concretely, as a meromorphic section of the cotangent bundle $\omega_{\sfSigma}$, we see that $\lambda$ has a double zero at the preimages of all zeros of $\phi$ and it has poles at the preimages of all the poles of $\phi$ of order $2$ or greater, but it is holomorphic and nonvanishing at the preimages of all simple poles $\sfop{Pol}_1 (\phi)$.
We denote the subsets of zeros and poles of $\lambda \in \omega_{\sfSigma}$ by
\begin{eqn}
	\sf{Zer} (\lambda) 
		\coleq \sfop{Supp} (2R \smallsetminus \pi^\ast D)
		= \pi^\inv \sf{Zer} (\phi)
\qtext{and}
	\sf{Pol} (\lambda)
		\coleq \sfop{Supp} (\sfDelta)
		= \pi^\inv \sf{Pol}_{\geq 2} (\phi)
\fullstop
\end{eqn}
On the other hand, $\pi^\ast D - 2R = \sfDelta - 2R + \pi^\ast D \cap 2R = \sfDelta - (2R \smallsetminus \pi^\ast D)$, so we find that $\lambda$ is a global nonvanishing section
\begin{eqntag}
\label{241018113808}
	\lambda 
		\in H^0 \big( \sfSigma, \omega_{\sfSigma, \sfDelta} (- 2R \smallsetminus \pi^\ast D)\big)
		\subset H^0 ( \sfSigma, \omega_{\sfSigma, \sfDelta} )
\fullstop
\end{eqntag}

\paragraf
Concretely, if $x$ is a local coordinate centred at a point $p \in X$, then the tautological one-form is $\lambda_\rm{taut} = y\d{x}$ or $y \d{x}/x^k$ depending on whether or not $p$ is away from the divisor $D$, where $y$ is the linear coordinate in the fibres of $T^\ast_{X,D}$.
Then $\lambda$ is the restriction of $\lambda_\rm{taut}$ restriction to $\sfSigma$.
Specifically, if $\phi (x) = \Q(x) \d{x}^2$, then roughly speaking $\lambda (x) = \pm \sqrt{ \Q(x) } \d{x}$, but this expression is not well-defined at all $p$.
Instead, we can use local normal forms from \autoref{241017104429} to give an accurate coordinate expression as follows.

If $p$ is a regular point, choose $x$ such that $\phi (x) = \d{x}^2$.
In this case, $\sfSigma$ is unramified over $p$ so we can use $x$ as a coordinate near either of the two preimages of $p$, and in this coordinate, $\lambda (x) = \pm \d{x}$.

If $p \in \sfop{Pol}_\textup{ev} (\phi)$ has order $m = 2k$, choose $x$ in which $\phi$ has local normal form \eqref{241017105802}.
In this case, $\sfSigma$ is again unramified so we use $x$ as a coordinate near either of the two preimages of $p$, and in this coordinate,
\begin{eqntag}
\label{241018111318}
	\lambda (x) 
		= \pm \big( 1 - k + \sqrt{r} x^{k-1} \big) \tfrac{\d{x}}{x^k}
		= \pm (1-k) \tfrac{\d{x}}{x^k} \pm \sqrt{r} \tfrac{\d{x}}{x}
\fullstop
\end{eqntag}
Note that we should think that the first formula expresses $\lambda (x)$ as a section of $\omega_{\sfSigma, \sfDelta}$ (with generator $x^{-k} \d{x}$), whilst the second formula expresses it as a section of $\omega_{\sfSigma}$ (with generator $\d{x}$).

On the other hand, if $p \in \sfop{Zer} (\phi)$ is a simple zero, choose $x$ such that $\phi (x) = x\d{x}^2$.
In this case, $\sfSigma$ is ramified over $p$, so $x$ does not give a well-defined coordinate near the preimage of $p$.
Instead, we use the linear fibre coordinate $y$ on the cotangent bundle $T^\ast_{X,D}$ near $p$ to parameterise this neighbourhood in $\sfSigma$.
Using the relation $y^2 = x$ which locally cuts out $\sfSigma$, we have $\pi^\ast \phi (y) = y^2 \d{(y^2)}^2 = (2y^2 \d{y})^2$.
Thus, 
\begin{eqntag}
\label{241019071032}
	\lambda (y) = 2y^2 \d{y}
\fullstop
\end{eqntag}
Note that $\lambda$ has a double zero at this transition point.

Similarly, if $p \in \sfop{Pol}_\textup{od} (\phi)$ is an odd pole of order $m = 2k - 1$, then $\sfSigma$ has a simple ramification point over $p$.
Choose $x$ such that $\phi (x) = x \left( \tfrac{\d{x}}{x^k} \right)^2$, and use the linear fibre coordinate $y$ on $T^\ast_{X,D}$ near $p$ to parameterise the neighbourhood of the preimage of $p$ on $\sfSigma$.
Using the local equation $y^2 = x$ for $\sfSigma$, we have $\pi^\ast \phi (y) = y^2 \left( y^{-2k} 2y\d{y} \right)^2 = \left( 2y^{-2(k-1)}\d{y} \right)^2$.
Thus, 
\begin{eqntag}
\label{241019071039}
	\lambda (y) = 2y^{-2(k-1)}\d{y}
\fullstop
\end{eqntag}
Thus, $\lambda$ defines a meromorphic differential form near the preimage of $p$ which has a pole of order $2(k-1)$.
In particular, if $p \in \sfop{Pol}_1 (\phi)$ is a simple pole of $\phi$, then $\lambda (y) = 2y \d{y}$ is a holomorphic and nonvanishing differential form.
On the other hand, to view $\lambda$ as a section of $\omega_{\sfSigma, \pi^\ast D}$ as in \eqref{241018110723}, we must write in terms of the generator $y^{-2k} \d{y}$.
Then we can clearly see that $\lambda \in \omega_{\sfSigma, \pi^\ast D}$ has a double zero at the preimage of $p$ for any $k$:
\begin{eqntag}
\label{241019071043}
	\lambda (y) = 2y^2 (y^{-2k} \d{y})
\fullstop
\end{eqntag}

\subsection{Critical Paths, Period, and Central Charge}

Suppose $\phi$ is a simple quadratic differential on a marked curve $(X,D)$, and let $\sfSigma$ be the associated spectral cover.

\paragraf
\label{241021121815}
Recall that the \dfn{fundamental groupoid} $\Pi_1 ({\sfSigma_\sfDelta})$ is the groupoid parameterising paths on ${\sfSigma_\sfDelta} = \sfSigma \smallsetminus \sfDelta$ up to homotopy with fixed endpoints.
It is a $2$-dimensional holomorphic manifold equipped with the \dfn{source} and \dfn{target maps}, which are holomorphic surjective submersions that extract the endpoints of a path $\gamma : [0,1] \to \sfSigma_\sfDelta$:
\begin{eqn}
	\rm{s}, \rm{t} : \Pi_1 ({\sfSigma_\sfDelta}) \to {\sfSigma_\sfDelta}
\qqtext{given by}
	\rm{s} : \gamma \mapsto \gamma (0)
\qtext{and}
	\rm{t} : \gamma \mapsto \gamma (1)
\fullstop
\end{eqn}
Viewing each point $p$ in ${\sfSigma_\sfDelta}$ as a constant path $1_p$ at $p$ defines an embedding ${\sfSigma_\sfDelta}$ into $\Pi_1 ({\sfSigma_\sfDelta})$ as a holomorphic submanifold called the \dfn{identity bisection}.
It is called so because it is a global holomorphic section of both $\rm{s}$ and $\rm{t}$.

\begin{defn}[Borel surface]
\label{241021115044}
We denote the source fibre of any $p \in {\sfSigma_\sfDelta}$ by
\begin{eqntag}
\label{240916160311}
	\tilde{\sfSigma}_{\sfDelta,p} 
		\coleq \rm{s}^\inv (p) 
		= \set{ \gamma \in \Pi_1 ({\sfSigma_\sfDelta}) ~\big|~ \rm{s} (\gamma) = p }
\fullstop{,}
\end{eqntag}
and call it the \dfn{Borel surface} associated with $p$.
\end{defn}

The Borel surface is a simply connected Riemann surface, and the restriction of the target map $\rm{t} : \tilde{\sfSigma}_{\sfDelta,p} \to {\sfSigma_\sfDelta}$ is nothing but the universal cover of ${\sfSigma_\sfDelta}$ based at $p$.
So any path $\gamma : I \to \sfSigma_\sfDelta$ starting at $p$ determines a point in the Borel surface $\tilde{\sfSigma}_{\sfDelta,p}$.
Furthermore, the Borel surface associated with $p$ has a distinguished point represented by the constant path $1_p \in \tilde{\sfSigma}_{\sfDelta, p}$, which we refer to as the \dfn{origin} in $\tilde{\sfSigma}_{\sfDelta, p}$.
The point $p \in \sfSigma_\sfDelta$ itself has a neighbourhood which can be canonically identified with a neighbourhood of the constant path $1_p \in \tilde{\sfSigma}_{\sfDelta,p}$.

Recall that a \dfn{groupoid $1$-cocycle} on $\Pi_1 ({\sfSigma_\sfDelta})$ is a holomorphic map $\Z : \Pi_1 ({\sfSigma_\sfDelta}) \to \CC$ that respects the groupoid law: i.e., $\Z (\gamma_2 \circ \gamma_1) = \Z (\gamma_1) + \Z (\gamma_2)$ for all pairs of composable paths $\gamma_1, \gamma_2 \in \Pi_1 ({\sfSigma_\sfDelta})$.

\begin{defn}[period and central charge]
\label{240808175800}
Integrating the Liouville form $\lambda$ along paths on $\sfSigma_\sfDelta$ defines a groupoid 1-cocycle on the fundamental groupoid $\Pi_1 (\sfSigma_\sfDelta)$ which we call the \dfn{period}:
\begin{eqntag}
\label{241005091803}
	\S : \Pi_1 ({\sfSigma_\sfDelta}) \too \CC
\qtext{sending}
	\gamma \mapstoo \S (\gamma) \coleq \int_\gamma \lambda
\fullstop
\end{eqntag}
Similarly, using the anti-invariance property of $\lambda$, we define another meromorphic differential on $\sfSigma$ by $\sigma \coleq \lambda - \iota^\ast \lambda = 2 \lambda \in H^0 (\sfSigma, \omega_{\sfSigma,\sfDelta})$.
Integrating the one-form $\sigma$ along paths on ${\sfSigma_\sfDelta}$ defines a groupoid $1$-cocycle which we call the \dfn{central charge}:
\begin{eqntag}
\label{240620104057}
	\Z : \Pi_1 ({\sfSigma_\sfDelta}) \too \CC
\qtext{sending}
	\gamma \mapstoo \Z (\gamma) \coleq \int_\gamma \sigma
\fullstop
\end{eqntag}
\end{defn}

\begin{defn}
\label{240725152318}
A \dfn{critical path} is any $\gamma \in \Pi_1 ({\sfSigma_\sfDelta})$ which terminates at a transition point; i.e., $\rm{t} (\gamma) \in \sfop{Tran} (\phi)$.
We call it a \dfn{simple critical path} if it terminates at a simple transition point; i.e., $\rm{t} (\gamma) \in \sfop{Tran}_0 (\phi)$; otherwise, we call it \dfn{non-simple}.
We denote the set of all critical paths and simple critical paths on ${\sfSigma_\sfDelta}$ respectively by
\begin{eqn}
	\Gamma \coleq \rm{t}^\inv (\sfop{Tran} (\phi)) \subset \Pi_1 ({\sfSigma_\sfDelta})
\qtext{and}
	\Gamma_0 \coleq \rm{t}^\inv (\sfop{Tran}_0 (\phi)) \subset \Pi_1 ({\sfSigma_\sfDelta})
\fullstop
\end{eqn}
Correspondingly, denote the set of all critical paths starting at $p \in {\sfSigma_\sfDelta}$ by 
\begin{eqn}
	\Gamma_p \coleq \Gamma \cap \tilde{\sfSigma}_{\sfDelta,p}
\qtext{and}
	\Gamma_{p,0} \coleq \Gamma_0 \cap \tilde{\sfSigma}_{\sfDelta,p}
\fullstop
\end{eqn}
Finally, denote the corresponding sets of central charges by
\begin{eqn}
	\Xi_p \coleq \Z (\Gamma_p) \subset \CC
\qtext{and}
	\Xi_{p,0} \coleq \Z (\Gamma_{p,0}) \subset \CC
\fullstop
\end{eqn}
\end{defn}

\begin{lem}
\label{241018134656}
The set $\Gamma_p$ is a discrete subset of the Borel surface $\tilde{\sfSigma}_{\sfDelta,p}$.
\end{lem}

\begin{proof*}
If $\gamma_1, \gamma_2$ are two critical paths starting at $p$ and terminating at the same transition point $p_\ast \in R_\sfDelta$, then the composition $\gamma_2^\inv \circ \gamma_1$ is a loop based at $p$.
So $\Gamma_p \cap \rm{t}^\inv (p_\ast)$ for each $p_\ast$ is a torsor for the fundamental group $\pi_1 ({\sfSigma_\sfDelta}, p)$.
In other words, there is a non-canonical isomorphism
\begin{eqn}
	\Gamma_p \cong \pi_1 ({\sfSigma_\sfDelta}, p)^{\#R_\sfDelta}
\end{eqn}
where $\#R_\sfDelta$ is the total number of transition points.
\end{proof*}

\paragraf
An alternative way to see that $\Gamma_p$ is discrete is to analyse the critical points of the central charge $\Z$ as follows.

\begin{prop*}
\label{240328174841}
For any $p \in {\sfSigma_\sfDelta}$, the central charge $\Z : \Pi_1 ({\sfSigma_\sfDelta}) \to \CC$ restricts to the Borel surface $\tilde{\sfSigma}_{\sfDelta,p}$ to yield a branched covering map
\begin{eqntag}
\label{240307121533}
	\Z_p : \tilde{\sfSigma}_{\sfDelta,p} \to \CC
\end{eqntag}
with only algebraic singularities.
Its ramification locus is the subset of simple critical paths $\Gamma_{p,0} \subset \Gamma_p \subset \tilde{\sfSigma}_{\sfDelta,p}$ starting at $p$.
In other words, the pair $(\tilde{\sfSigma}_{\sfDelta, p}, \Z_p)$ is an endless Riemann surface of algebraic type and $\sfop{Sing} (\tilde{\sfSigma}_{\sfDelta, p}, \Z_p) = \Gamma_{p,0}$.
In particular, if $p$ is regular, then it has a neighbourhood which is identified by $\Z_p$ with a neighbourhood of the origin in $\CC$.
\end{prop*}

\begin{proof*}
A ramification point of $\Z_p$ is by definition any path $\gamma \in \tilde{\sfSigma}_{\sfDelta,p}$ whose terminal point $q = \rm{t} (\gamma)$ is such that the differential pushforward map $\d \Z_p = \sigma$ drops rank at $q$.
In other words, $q$ is a zero of $\sigma$, hence a simple transition point.

To be completely explicit, if $\gamma \in \Gamma_p$, let $\xi \coleq \Z (\gamma) \in \CC$ be its central charge.
Choose a local coordinate $x$ on ${\sfSigma_\sfDelta}$ centred at the terminal point $\rm{t} (\gamma) \in \sfop{Tran} (\phi)$ with respect to which $\sigma$ is $\d{(x^r)}$ for some integer $r \geq 1$.
We can treat $x$ as a local coordinate on $\tilde{\sfSigma}_{\sfDelta,p}$ centred at $\gamma$.
Then with respect to this coordinate, the map $\Z_p$ is given by
\begin{eqn}
	\Z_p : x \mapstoo t
		= \xi + x^r
		= \xi + \int_0^x \d{(x^r)}
\fullstop
\end{eqn}
When the turning point is simple, the Liouville form $\lambda$ has double zero at the corresponding ramification point, which means $r = 3$.
Therefore, in the neighbourhood of the point in $\tilde{\sfSigma}_{\sfDelta, p}$ representing a simple critical path $\gamma$, the map $\Z_p$ is a cyclic $3$:$1$ ramified cover, branched over the central charge $\xi$ of $\gamma$.
On the other hand, when the turning point is a simple pole, the Liouville form $\lambda$ is nonvanishing at the corresponding ramification point, which means $r = 1$.
Thus, in the neighbourhood of the point in $\tilde{\sfSigma}_{\sfDelta, p}$ representing a non-simple critical path $\gamma$, the map $\Z_p$ is $1$:$1$.
\end{proof*}

We can now immediately conclude that the fundamental groupoid $\Pi_1 (\sfSigma_\sfDelta)$ has the following resurgent-geometric structure.

\begin{cor*}
\label{241015112906}
The triple $\big( \Pi_1 (\sfSigma_\sfDelta), \rm{s}, \Z \big)$ --- consisting of the fundamental groupoid of the punctured spectral curve $\sfSigma_\sfDelta$, the groupoid source map $\rm{s} : \Pi_1 (\sfSigma_\sfDelta) \to \sfSigma_\sfDelta$, and the central charge $\Z : \Pi_1 (\sfSigma_\sfDelta) \to \CC$ --- is a $\sfSigma_\sfDelta$-family of endless Riemann surfaces of algebraic type.
In particular, every Borel surface $(\tilde{\sfSigma}_{\sfDelta, p}, \Z_p)$, $p \in \sfSigma_\sfDelta$, is an endless Riemann surface of algebraic type, and $\sfop{Sing} (\tilde{\sfSigma}_{\sfDelta, p}, \Z_p) = \Gamma_{p,0}$.
\end{cor*}

\paragraf\label{240801145442}%
Given a critical path $\gamma \in \Gamma_p$, we can consider its reversed reflection $\check{\gamma} \coleq \iota (\gamma)^\inv$.
They are composable paths because $\rm{t} (\gamma) = \rm{s} (\check{\gamma})$.
The composition $\check{\gamma} \circ \gamma$ terminates at the point $\check{p} \coleq \iota (p)$ which lies in the same fibre of $\pi$ as the point $p$.
Thus, any critical path canonically gives rise to a path that switches sheets of the spectral curve.
This defines a surjective, but never injective, map 
\begin{eqn}
	\Gamma_p \to \rm{s}^\inv (p) \cap \rm{t}^\inv (\check{p})
\quad:\quad
	\gamma \mapstoo \check{\gamma} \circ \gamma
\fullstop
\end{eqn}
Moreover, from the anti-invariance of $\lambda$, we conclude that the central charge of a critical path $\gamma$ is the period of $\check{\gamma} \circ \gamma$:
\begin{eqn}
	\Z (\gamma)
	= \int_\gamma \big( \lambda - \iota^\ast \lambda \big)
	= \int_\gamma \lambda - \int_{\iota (\gamma)} \lambda
	= \S (\check{\gamma} \circ \gamma)
\fullstop
\end{eqn}

\paragraf
The inclusion ${\sfSigma_{R,\sfDelta}} \inj {\sfSigma_\sfDelta}$ gives rise to a holomorphic map $\Pi_1 ({\sfSigma_{R,\sfDelta}}) \to \Pi_1 ({\sfSigma_\sfDelta})$ which is neither injective nor surjective.
However, if we restrict the groupoid $\Pi_1 ({\sfSigma_\sfDelta})$ to ${\sfSigma_{R,\sfDelta}}$, then we do get a surjective holomorphic map
\begin{eqn}
	\Pi_1 ({\sfSigma_{R,\sfDelta}}) \too \Pi_1 ({\sfSigma_\sfDelta}) \big|_{{\sfSigma_{R,\sfDelta}}} \coleq \rm{s}^\inv ({\sfSigma_{R,\sfDelta}}) \cap \rm{t}^\inv ({\sfSigma_{R,\sfDelta}})
\fullstop	
\end{eqn}
In particular, for any $p \in {\sfSigma_{R,\sfDelta}}$, the restriction of this map to the source fibre
\begin{eqn}
	\tilde{\sfSigma}_{R,\sfDelta, p} \coleq \rm{s}^\inv (p) \subset \Pi_1 ({\sfSigma_{R,\sfDelta}})
\end{eqn}
is nothing but the universal cover of the Borel surface $\tilde{\sfSigma}_{\sfDelta, p} \subset \Pi_1 (\sfSigma_\sfDelta)$ punctured along the locus of critical paths $\Gamma_p$ and based at the point $1_p$:
\begin{eqn}
	\tilde{\sfSigma}_{R, \sfDelta,p} = \widetilde{ \tilde{\sfSigma}_{\sfDelta,p} \!\smallsetminus\! \Gamma_p } \too \tilde{\sfSigma}_{\sfDelta,p} \smallsetminus \Gamma_p
\fullstop
\end{eqn}
Furthermore, any simply-connected neighbourhood of the point $p$ in $\tilde{\sfSigma}_{\sfDelta,p}$, viewed as the constant path $1_p$ at $p$, is canonically identified with a neighbourhood of $p$ in $\tilde{\sfSigma}_{R, \sfDelta,p}$.
Any path $\gamma : I \to \sfSigma_\sfDelta$ starting at $p$ which avoids the transition points $\sfop{Tran} (\phi)$ is canonically a point in the source fibre $\tilde{\sfSigma}_{R, \sfDelta, p}$.
At the same time, the source fibre $\tilde{\sfSigma}_{R,\sfDelta,p}$ is simply connected, so $\gamma$ determines a path $\tilde{\gamma} : I \to \tilde{\sfSigma}_{R,\sfDelta,p}$, unique up to fixed-endpoint homotopy, which starts at $1_p$ and terminates at $\gamma$ in $\tilde{\sfSigma}_{R,\sfDelta,p}$.

\begin{cor*}
\label{241015113613}
The triple $\big( \Pi_1 (\sfSigma_{R,\sfDelta}), \rm{s}, \Z \big)$ --- consisting of the fundamental groupoid of the subset $\sfSigma_{R,\sfDelta}$ of regular points in the spectral curve, the groupoid source map $\rm{s} : \Pi_1 (\sfSigma_{R,\sfDelta}) \to \sfSigma_{R,\sfDelta}$, and the central charge $\Z : \Pi_1 (\sfSigma_{R,\sfDelta}) \to \CC$ --- is a $\sfSigma_{R,\sfDelta}$-family of endless Riemann surfaces of log-algebraic type.
In particular, every source fibre $(\tilde{\sfSigma}_{R,\sfDelta, p}, \Z_p)$, $p \in \sfSigma_{R,\sfDelta}$, is an endless Riemann surface of log-algebraic type.
\end{cor*}

\begin{proof*}
For each $p \in \sfSigma_{R,\sfDelta}$, the pair $(\tilde{\sfSigma}_{\sfDelta, p}, \Z_p)$ is an endless Riemann surface of algebraic type by \autoref{241015112906}, and $\Gamma_p$ is a discrete subset of $\tilde{\sfSigma}_{\sfDelta, p}$ by \autoref{241018134656}.
\end{proof*}

\subsection{Critical Trajectories and Stokes Rays}

The spectral curve $\sfSigma$ comes equipped with a natural $\SS^1$-family of singular foliations, called the \textit{trajectory structure}.
It is obtained by integrating the singular distribution inside the real tangent bundle of $\sfSigma$ determined by the real part of the rotated canonical differential $\Re (e^{i \alpha} \sigma)$.
We now briefly recall the properties of this foliation.

\begin{defn}[trajectories]
\label{240916210413x}
A \dfn{trajectory} with \dfn{phase} $\alpha \in \RR / 2\pi\ZZ$ is a smooth real curve $\gamma : I \to \sfSigma$ whose tangent vector $\dot{\gamma}$ satisfies $\sigma (\dot{\gamma}) \in \RR_\alpha$ everywhere along $\gamma$.
Equivalently, the projection of a trajectory to $\CC$ via the central charge $\Z$ is a line segment (possibly unbounded) with phase $\alpha$.

A \dfn{trajectory starting at} $p \in \sfSigma$ (or a \dfn{trajectory ray}) is a maximal trajectory $\gamma$ with fixed starting point $p$.
A \dfn{critical trajectory} is a trajectory ray whose homotopy class in $\sfSigma_\sfDelta$ is a critical path; otherwise we call it a \dfn{regular trajectory}.
The phase of a regular or a critical trajectory starting at $p \in \sfSigma$ is called respectively a \dfn{regular phase} or a \dfn{critical phase} at $p$.
A \dfn{saddle trajectory} is a trajectory that both starts and terminates at transition points, not necessarily distinct ones.
\end{defn}

For every fixed phase $\alpha$, maximal trajectories with phase $\alpha$ form a real-$1$-dimensional singular foliation (sometimes called the \dfn{trajectory foliation} or the \dfn{geodesic foliation}) of the spectral curve $\sfSigma$ with singularities at the simple transition points $\sfop{Tran}_0 (\phi_0)$ and the preimages of non-simple poles $\pi^\inv \sfop{Pol}_{\geq 2} (\phi_0)$.
Note that this foliation is smooth at all non-simple transition points $\sfop{Tran}_1 (\phi_0) = \pi^\inv \sfop{Pol}_1 (\phi_0)$.

\paragraph{Trajectory-based Stokes diagrams.}
For every phase $\alpha \in \RR / 2\pi\ZZ$, any regular point $p \in \sfSigma_{R,\sfDelta}$ lies on a unique maximal trajectory $\gamma$ with phase $\alpha$.
This means that points on the tangent circle $S_p \sfSigma$ at $p$ are in one-to-one correspondence with trajectories starting at $p$.
Therefore, the points of the tangent circle $S_p \sfSigma$ get labelled as regular or critical depending on the type of maximal trajectory passing through $p$ in that tangential direction.
In analogy with the stability notions in a resurgent diagram (\autoref{240809103556}), we introduced the following terminology.

\begin{defn}[stability in trajectory-based Stokes diagrams]
\label{240924153408}
Let $p \in \sfSigma$ be a regular point.
We say that a regular trajectory with phase $\alpha$ starting at $p$ is 
\begin{enumerate}
\item \dfn{stable} if $\alpha$ is contained in an open arc $A = (\alpha_-, \alpha_+)$ such that every maximal trajectory with phase $\alpha' \in A$ starting at $p$ is regular;
\item \dfn{semi-stable} if $p$ has a neighbourhood $U \subset \sfSigma$ of regular points such that the maximal trajectory with phase $\alpha$ starting at any $p' \in U$ is also regular;
\item \dfn{unstable} otherwise.
\end{enumerate}
We also say that a critical trajectory with phase $\alpha$ starting at $p$ is
\begin{enumerate}
\item \dfn{stable} if $\alpha$ is contained in an open arc $A = (\alpha_-,\alpha_+)$ such that all maximal trajectories starting at $p$ with phase $\alpha' \in A \smallsetminus \set{\alpha}$ are regular;
\item \dfn{semi-stable} if $p$ has a neighbourhood $U \subset X$ with the following property: for every $p' \in U$, the maximal trajectory with phase $\alpha$ starting at $p'$ is critical if and only if $p'$ and $p$ belong to the same trajectory with phase $\alpha$;
\item \dfn{unstable} otherwise.
\end{enumerate}
\end{defn}

\section{The WKB Method for Schrödinger Equations}
\label{241022163033}

\subsection{The Schrödinger Equation}

In this subsection, we introduce a coordinate-invariant notion of Schrödinger equations and impose some mild assumptions on its singular points.

\begin{defn}
\label{240818084025}
A \dfn{spin structure} on $X$ is a choice of square root of $\omega_X$; i.e., a line bundle $\cal{L}$ on $X$ together with an isomorphism $\cal{L} \otimes \cal{L} \iso \omega_X$.
Once a spin structure has been chosen, we use the more suggestive notation for the \dfn{spin bundle} $\omega^\sqroot_X \coleq \cal{L}$, as well as $\omega^\invsqroot_X \coleq \cal{L}^\inv$ and $\omega^\cubesqroot_X \coleq \omega^\invsqroot_X \otimes \omega^2_X$.
\end{defn}

A spin structure (aka \textit{theta characteristic}) on a genus $g$ compact curve always exists but is not unique: there are $2^{2g}$ non-isomorphic spin structures.
Explicitly, if $(U,x)$ is any local coordinate chart on $X$, the cotangent bundle $\omega_X$ acquires a canonical generator $\d{x}$ over $U$.
This generator in turn induces a generator on the spin bundle which is often denoted by $\d{x}^\sqroot$.
It is characterised by the property that $\d{x}^\sqroot \otimes \d{x}^\sqroot = \d{x}$ via the isomorphism $\omega_X^\sqroot \otimes \omega_X^\sqroot \iso \omega_X$.
So this generator $\d{x}^\sqroot$ is uniquely determined up to sign.

\begin{defn}[Schrödinger operator]
\label{240730095844}
Let $(X,D)$ be a marked curve with a chosen spin structure $\omega^\sqroot_X$.
A \dfn{Schrödinger operator} on $(X,D)$ is a family of $\CC$-linear maps of sheaves
\begin{eqntag}
\label{240730101105}
	\P : \omega^\invsqroot_X \too \omega^\invsqroot_X \otimes \omega^2_{X,D} = \omega^\cubesqroot_X (2D)
\end{eqntag}
parameterised by $\hbar \in \CC$ and characterised by the following property.
In any local coordinate chart $(U,x)$ on $X$, the map $\P$ must have the following expression in terms of the induced local generators $\d{x}^\invsqroot$ and $\d{x}^\cubesqroot$:
\begin{eqntag}
\label{240730103038}
	\P
		\loceq \hbar^2 \del_x^2 - \Q (x, \hbar)
\fullstop{,}
\end{eqntag}
where $\Q (x, \hbar)$ is a polynomial in $\hbar$ whose coefficients are meromorphic functions with poles bounded by the divisor $2D$:
\begin{equation}
\label{240730105621}
	\Q (x, \hbar)
		= \sum_{k \geq 0}^{\textup{\tiny finite}} \Q_{k} (x) \hbar^k
\qqtext{where}
	\Q_k \in \cal{O}_U (2D)
\fullstop
\end{equation}
Furthermore, we shall always assume that the divisor $D$ is minimal in the sense that $\P$ does not restrict to $(X,D')$ for any proper subdivisor $D' < D$.
\end{defn}

\paragraf
A \dfn{Schrödinger equation} on $(X,D)$ is then the equation
\begin{eqntag}
\label{240730161252}
	\P \psi = 0
\end{eqntag}
for an $\hbar$-family of local sections $\psi$ of the dual spin bundle $\omega^\invsqroot_X$.
Such families will be called \dfn{wave functions}.
In a local coordinate chart $(U,x)$, any section $\psi \in \omega^\invsqroot_X$ has a canonical (up to sign) representation 
\begin{eqn}
	\psi \loceq \Psi (x, \hbar) \d{x}^\invsqroot
\end{eqn}
If $\psi$ is a wave function, we call its local representation $\Psi$ a \dfn{local wave function}.
Then local wave functions are solutions of a \dfn{local Schrödinger equation}:
\begin{eqntag}
\label{240730163844}
	\hbar^2 \del_x^2 \Psi (x, \hbar) = \Q (x, \hbar) \Psi (x, \hbar)	
\fullstop
\end{eqntag}
Insisting that this Schrödinger form is retained in \textit{every} local coordinate imposes the following nontrivial transformation law on the \dfn{local potential} $\Q$ whose verification is a straightforward calculation.

\begin{lem}
\label{240728093124}
If $\Psi (x, \hbar)$ is a solution of the local Schrödinger equation \eqref{240730163844} and $x = x (\tilde{x})$ is a holomorphic coordinate transformation, then
\begin{eqntag}
\label{240728092147}
	\tilde{\Psi} (\tilde{x}, \hbar) \coleq \Psi (x, \hbar) \dot{x}^\invsqroot
\qqtext{where}
	\dot{x} \coleq \d{x} / \d{\tilde{x}}
\fullstop{,}
\end{eqntag}
is a solution of $\hbar^2 \del_{\tilde{x}}^2 \tilde{\Psi} (\tilde{x}, \hbar) =  \tilde{\Q} (\tilde{x}, \hbar) \tilde{\Psi} (\tilde{x}, \hbar)$ with potential
\begin{eqntag}
\label{240731112139}
	\tilde{\Q} (\tilde{x}, \hbar) \coleq \Q (x , \hbar) \dot{x}^2 - \tfrac{\hbar^2}{2}  \cbrac{ x ; \tilde{x}}
\fullstop{,}
\end{eqntag}
where $\cbrac{ x ; \tilde{x}}$ is the Schwarzian derivative:
\begin{eqn}
	\cbrac{ x ; \tilde{x}}
		\coleq \frac{\d}{\d{\tilde{x}}} \left( \frac{ \ddot{x} }{ \dot{x} } \right)
			- \frac{1}{2} \left( \frac{ \ddot{x} }{ \dot{x} } \right)^2
		= \frac{ \dddot{x} }{ \dot{x} } - \frac{3}{2} \left( \frac{ \ddot{x} }{ \dot{x} } \right)^2
\fullstop
\end{eqn}
\end{lem}

\paragraf
It follows from the transformation law \eqref{240731112139} that the local potentials $\Q$ would determine a well-defined quadratic differential $\Q (x, \hbar) \d{x}^2$ were it not for the presence of the Schwarzian derivative term.
However, notice that the failure to transform like a quadratic differential is concentrated in order $2$ in $\hbar$.
Therefore, a Schrödinger operator $\P$ naturally determines quadratic differentials
\begin{eqntag}
\label{240918154226}
	\phi_k \loceq \Q_k (x) \d{x}^2 \in H^0 \big(X, \omega^2_{X,D} \big)
\qquad
	\text{for all $k \neq 2$.}
\end{eqntag}
The most important of these is the quadratic differential $\phi_0$, which we call the \dfn{classical quadratic differential} associated with $\P$.

\begin{assumption}
\label{240924155219}
If $\P$ is a Schrödinger operator on a marked curve $(X,D)$, let $\phi_0$ be the associated classical quadratic differential.
We make the following additional assumptions:
\begin{enumerate}
\item All zeros of $\phi_0$ are simple;
\item In any local coordinate $x$, if $x_\ast$ is a pole of some $\Q_k$, then $x_\ast$ is also a pole of $\Q_0$ of order no less; i.e., $\op{ord}_{x_\ast} (\Q_k) \leq \op{ord}_{x_\ast} (\Q_0)$.\end{enumerate}
\end{assumption}

\paragraf
It follows from the transformation law \eqref{240731112139} for the local potential $\Q$ that the poles of its coefficients $\Q_k$ as well as the zeros of $\Q_0$ are coordinate-independent data.
More precisely, if $p \in X$ is any point and $x, \tilde{x}$ are any two local coordinates near $p$, let $\Q, \tilde{\Q}$ be the corresponding local potentials.
Then, for every $k$, the coefficient $\Q_k$ of $\Q$ has a pole of order $m$ at $x_\ast = x(p)$ if and only if the coefficient $\tilde{\Q}_k$ of $\tilde{\Q}$ has a pole of order $m$ at $\tilde{x}_\ast = \tilde{x}(p)$.
We can therefore say that a point $p \in X$ is a pole of order $m$ of a given Schrödinger operator if its potential $\Q$ in any (and hence every) local coordinate near $p$ has a coefficient $\Q_k$ that has a pole at $p$ of order $m$.

\subsection{WKB Solutions}

Let $\P$ be a given Schrödinger operator on a marked curve $(X,D)$ with a chosen spin structure.
In this subsection, we describe in a coordinate-invariant way the WKB method for obtaining a special basis of formal solutions of the Schrödinger equation.

\begin{defn}
\label{240813114000}
A local \dfn{WKB solution} is any holomorphic family of local sections $\psi$ of $\omega^\invsqroot_X$ over an open set $U \subset X$, parameterised by $\hbar$ in a sector $V \subset \CC$ at the origin, which satisfies the Schrödinger equation, and which admits a uniform exponential perturbative series expansion, which is by definition a formal expression of the form
\begin{eqntag}
\label{240915142831}
	\hat{\psi} = e^{-\S/\hbar} \hat{a}
\qqtext{where}
	\hat{a} \in \omega^\invsqroot_X (U) \bbrac{\hbar}
\qtext{and}
	\S \in \cal{O}_X (U)
\fullstop
\end{eqntag}
\end{defn}

\begin{defn}
\label{240915142819}
A local \dfn{formal WKB solution} of the Schrödinger equation over $U \subset X$ is any exponential perturbative series \eqref{240915142831} which formally satisfies the Schrödinger equation.
A local \dfn{formal WKB transseries solution} is a formal (possibly countable) $\CC$-linear combination of formal WKB solutions.
\end{defn}

\begin{prop}
\label{240915131722}
Let $p_0 \in X$ be a regular point with preimages $p_0^\pm \in \sfSigma$.
Then there is a pair of canonical (up to common constant scale) formal WKB solutions, respectively called the \dfn{formal WKB solution based at $p_0^\pm$}:
\begin{eqntag}
\label{240915183003}
	\hat{\psi}_\pm = e^{-\S_\pm/\hbar} \hat{a}_\pm
\qtext{where}
	\hat{a}_\pm \in \omega^\invsqroot_{X,p_0} \bbrac{\hbar}
\qtext{and}
	\S_\pm (p) \coleq \S (p_\pm) = \int_{p_0^\pm}^{p_\pm} \lambda
\fullstop{,}
\end{eqntag}
for all $p_\pm$ near $p_0^\pm$ with projection $p \in X$.
In fact, there is a canonical (up to constant scale) generator $a_0$ of $\omega_{X,p_0}^\invsqroot$, called \dfn{canonical generator based at $p_0$}, such that $\hat{\psi}_\pm$ is the unique formal WKB solution near $p_0$ satisfying
\begin{eqntag}
\label{240915182858}
	\hat{\psi}_\pm (p_0) = a_0 (p_0)
\qqtext{and}
	\d{\S}_\pm (p_0) = \lambda (p_0^\pm)
\fullstop
\end{eqntag}
In terms of this generator,
\begin{eqntag}
\label{240916121540}
	\hat{\psi}_\pm = a_0 e^{-\S_\pm / \hbar} \hat{\A}_\pm
\qqtext{where}
	\hat{\A}_\pm \in \cal{O}_{X, p_0} \bbrac{\hbar}
\qtext{with}
	\hat{\A}_\pm (p_0) = 1
\fullstop
\end{eqntag}
\end{prop}

We will prove this Proposition in a few steps below, but let us make a couple of comments first.
\autoref*{240915131722} is the coordinate-invariant formulation of a well-known Formal Existence and Uniqueness Theorem for WKB solutions of Schrödinger equations; for example, see \cite[p.486]{MR1721985} or \cite[Proposition 3.1]{MY210623112236}.
Condition \eqref{240915182858} should be thought of as the initial condition at $p_0$.
In particular, the second equality in \eqref{240915182858} can be interpreted as the logarithmic derivative of $\hat{\psi}_\pm$ at $p_0$.
The canonical generator $a_0$ will be fully characterised below in \autoref{240815162128}.

\enlargethispage{10pt}

\begin{rem}
\label{241020142351}
Normalisations different from \eqref{240915182858} appear in the literature and will be described elsewhere.
We also remark that the `approximate solutions' 
\begin{eqn}
	a_0 e^{-\S_\pm / \hbar}	
\end{eqn}
are the coordinate-invariant form of the \dfn{WKB approximation} or the \dfn{semiclassical approximation} as traditionally known in quantum mechanics.
\end{rem}

\begin{rem}
\label{241015094508}
If $p_0 \in X$ is a regular point, the map $\pi : \sfSigma \to X$ is a canonical isomorphism between the infinitesimal neighbourhoods of $p_0$ and each of $p_0^+$ and $p_0^-$.
Using this isomorphism, we can regard each formal WKB solution $\hat{\psi}_\pm$ as defined near $p_0^\pm$ by identifying $\hat{\psi}_\pm$ with its pullback $\pi^\ast \hat{\psi}_\pm$.
From this point of view, we can say that for every regular point $p_1 \in \sfSigma$ on the spectral curve, there is just one canonical (up to constant scale) formal WKB solution: $\hat{\psi} = e^{-\S/\hbar} \hat{a}$ where $\hat{a} \in \big(\pi^\ast \omega^\invsqroot_X\big)_{p_1} \bbrac{\hbar}$ and $\S \in \cal{O}_{\sfSigma,p_1}$ is given by $\S (p) = \int_{p_1}^p \lambda$ for all $p$ near $p_1$.
Namely, once a generator $a_0$ of $\pi^\ast \omega_X^\invsqroot$ near $p_1$ has been chosen, $\hat{\psi}$ is the unique formal WKB solution which satisfies $\hat{\psi} (p_1) = a_0 (p_1)$ and $\d{\S} (p_1) = \lambda (p_1)$.
\end{rem}

\paragraf
To prove \autoref{240915131722}, let us begin by examining the given Schrödinger operator in a local coordinate chart $x$ near $p_0$:
\begin{eqntag}
\label{240810202521} 
	\hbar^2 \del_x^2 \Psi (x, \hbar) = \Q (x, \hbar) \Psi (x, \hbar)
\fullstop
\end{eqntag}
Then we make a change of the unknown variable known as the \dfn{WKB ansatz}:
\begin{eqntag}
\label{240815150637} 
	\Psi (x, \hbar) = \exp \left( -\frac{1}{\hbar} \int_{x_0}^x \Y (x', \hbar) \d{x'} \right)
\fullstop{,}
\end{eqntag}
where $x_0 \coleq x (p_0)$.
The point is that this expression for $\Psi$ satisfies the local Schrödinger equation \eqref{240810202521} if and only if $\Y$ satisfies the \dfn{local Riccati equation}
\begin{eqntag}
\label{240815150746}
	\hbar \del_x \Y = \Y^2 - \Q
\fullstop
\end{eqntag}
The advantage of this Riccati equation over the Schrödinger equation is that it is easily solved in perturbative series.
We therefore have the following explicit coordinate description of the formal WKB solutions $\hat{\psi}_\pm$ from \autoref*{240915131722}.

\begin{lem*}
\label{240915175932}
Let $x$ be any local coordinate near $p_0 \in X$, put $x_0 \coleq x (p_0)$, and fix the sign of the local generator $\d{x}^\invsqroot$.
Then, up to constant scale, the formal WKB solutions based at $p_0$ are given explicitly by the following formula:
\begin{eqntag}
\label{240915180049}
	\hat{\psi}_\pm (x) 
		= \hat{\Psi}_\pm (x, \hbar) \d{x}^\invsqroot
		= \exp \left( -\frac{1}{\hbar} \int_{x_0}^x \hat{\Y}_\pm (x', \hbar) \d{x'} \right) \d{x}^\invsqroot
\fullstop
\end{eqntag}
Here, the integral is the term-wise integration of the perturbative series
\begin{eqntag}
\label{240916120032}
	\hat{\Y}_\pm (x, \hbar) \coleq \sum_{k=0}^\infty y_k^\pm (x) \hbar^k \in \CC \cbrac{x} \bbrac{\hbar} \cong \cal{O}_{X,p_0} \bbrac{\hbar}
\fullstop{,}
\end{eqntag}
whose leading-order part is the local coordinate representation $\lambda (x) = y_0^\pm (x) \d{x}$ of the Liouville form $\lambda$ near the two preimages $p_0^\pm$ of $p_0$.
Moreover, the series $\hat{\Y}_\pm$ is the unique formal perturbative local solution with leading-order $y_0^\pm$ of the local Riccati equation \eqref{240815150746}, and all of its higher-order coefficients are uniquely determined from the leading-order coefficient $y_0^\pm$ by the \dfn{WKB recursion}:
\begin{eqntag}
\label{240915180751}
	y^2_0 = \Q_0,
\quad
	y_1 = \tfrac{1}{2y_0} \big(\del_x y_0 - \Q_1 \big),
\quad
	y_k = \frac{1}{2y_0} \left( \del_x y_{k-1} - \sum_{i + j = k}^{i,j\neq0} y_i y_j - \Q_k \right)
\fullstop
\end{eqntag} 
\end{lem*}

The proof of this Lemma involves a straightforward calculation.
Namely, the fact that the perturbative series $\hat{\Y}_\pm (x, \hbar)$ formally satisfies the Riccati equation \eqref{240815150746} can be verified using the WKB recursion \eqref{240915180751}.
The uniqueness of $\hat{\Y}_\pm$ is also evident from the WKB recursion.
The only assertion that is not straightforward is that the formal local section $\hat{\psi}_\pm$ defined by formula \eqref{240915180049} is actually independent of the chosen coordinate $x$.
This is not a trivial statement: just as the local Schrödinger equation \eqref{240810202521} depends on the coordinate $x$, so do the WKB ansatz \eqref{240815150637} and the local Riccati equation \eqref{240815150746}.
The following Lemma describes this dependence.

\begin{lem}
\label{240728090056}
If $\Y (x, \hbar)$ is a solution of the Riccati equation \eqref{240815150746}, and $x = x(\tilde{x})$ is a holomorphic coordinate transformation, then
\begin{eqntag}
\label{240728091232}
	\tilde{\Y} (\tilde{x}, \hbar) \coleq \Y (x, \hbar) \dot{x} + \frac{\hbar}{2} \frac{ \ddot{x} }{ \dot{x} }
\fullstop{,}
\end{eqntag}
is a solution of the Riccati equation $\hbar \del_{\tilde{x}} \tilde{\Y} = \tilde{\Y}^2 - \tilde{\Q}$, where $\tilde{\Q}$ is related to $\Q$ by \eqref{240731112139}.
\end{lem}

\begin{proof*}
Calculate using \eqref{240728091232} and \eqref{240731112139} as follows:
\begin{eqns}
	\hbar \tfrac{\d}{\d{\tilde{x}}} \tilde{\Y}
		&= \hbar \tfrac{\d}{\d{\tilde{x}}} \left( \dot{x} \Y 
			+ \tfrac{\hbar}{2} \tfrac{ \ddot{x} }{ \dot{x} } \right)
\\		&= \hbar \ddot{x} \Y 
			+ (\dot{x})^2 \hbar \tfrac{\d}{\d{x}} \Y 
			+ \tfrac{\hbar^2}{2} \tfrac{\d}{\d{\tilde{x}}} \left( \tfrac{\ddot{x}}{\dot{x}} \right)
\\		&= \hbar \ddot{x} \Y 
			+ (\dot{x})^2 \big( \Y^2 - \Q \big)
			+ \tfrac{\hbar^2}{2} \tfrac{\d}{\d{\tilde{x}}} \left( \tfrac{\ddot{x}}{\dot{x}} \right)
\\		&= \hbar \tfrac{\ddot{x}}{\dot{x}} \big( \tilde{\Y} - \tfrac{\hbar}{2} \tfrac{\ddot{x}}{\dot{x}} \big)
			+ \big( \tilde{\Y} - \tfrac{\hbar}{2} \tfrac{\ddot{x}}{\dot{x}} \big)^2 - \dot{x}^2 \Q
			+ \tfrac{\hbar^2}{2} \tfrac{\d}{\d{\tilde{x}}} \left( \tfrac{\ddot{x}}{\dot{x}} \right)
\\		&= \tilde{\Y}^2 
			- \Big( \dot{x}^2 \Q 
				- \tfrac{\hbar^2}{2} \tfrac{\d}{\d{\tilde{x}}} \big( \tfrac{\ddot{x}}{\dot{x}} \big) 
				+ \tfrac{\hbar^2}{4} \big( \tfrac{\ddot{x}}{\dot{x}} \big)^2 \Big)
\\		&= \tilde{\Y}^2 - \tilde{\Q}
\fullstop
\tag*{\qedhere}
\end{eqns}	
\end{proof*}

\paragraf
It follows from this Lemma that the integrand $\Y \d{x}$ in the WKB ansatz \eqref{240815150637} does \textit{not} determine a well-defined differential form near $p_0 \in X$ even when the square-root ambiguity in $y_0$ has been fixed.
The problem is that $\Y$ does not obey the correct transformation law under coordinate changes.
That is to say, if we take another local coordinate $\tilde{x}$ near $p_0$, put $\tilde{x}_0 \coleq \tilde{x} (p_0)$, and write down the new WKB ansatz
\begin{eqn}
	\tilde{\Psi} (\tilde{x}, \hbar) = \exp \left( -\frac{1}{\hbar} \int_{\tilde{x}_0}^{\tilde{x}} \tilde{\Y} (\tilde{x}', \hbar) \d{\tilde{x}'} \right)
\fullstop{,}
\end{eqn}
then the differential forms $\Y \d{x}$ and $\tilde{\Y} \d{\tilde{x}}$ are in general \textit{not} the same.
Observe, however, that the failure of $\Y$ to transform like a differential form is concentrated at exactly the first order in $\hbar$.
We can take advantage of this observation as follows.

From the WKB recursion \eqref{240915180751}, the Riccati equation \eqref{240815150746} has the following approximate solutions valid to first order in $\hbar$:
\begin{eqntag}
\label{240815161025}
	y_0 + y_1 \hbar
\qqqtext{where}
	y^2_0 = \Q_0
\qtext{and}
	y_1 = \tfrac{1}{4} \del_x \log \Q_0 - \frac{\Q_1}{2y_0}
\fullstop
\end{eqntag}
Here, we have used the fact that $\del_x \log y_0 = \frac{1}{2} \del_x \log \Q_0$ for any of the two possible determinations $y_0 = \pm \sqrt{\Q_0}$.
It follows from \eqref{240728091232} that $y_0$ and $y_1$ transform under coordinate changes like so:
\begin{eqntag}
\label{240815183519}
	\tilde{y}_0 = y_0 \dot{x}
\qqtext{and}
	\tilde{y}_1 = y_1 \dot{x} + \tfrac{1}{2} \ddot{x}/\dot{x}
\fullstop
\end{eqntag}
At the same time, from the potential's coordinate transformation law \eqref{240731112139} we see that $\Q_1$ transforms like $\tilde{\Q}_1 = \Q_1 \dot{x}^2$.
This means the second term in \eqref{240815161025} transforms like $\tilde{\Q}_1/2\tilde{y}_0 = \dot{x} \Q_1/2y_0$.
Consequently, the first term in the equation for $y_1$ in \eqref{240815161025} transforms like so:
\begin{eqntag}
\label{240917104315}
	\tfrac{1}{4} \del_{\tilde{x}} \log \tilde{\Q}_0
		= \dot{x} \tfrac{1}{4} \del_{x} \log \Q_0 + \tfrac{1}{2} \ddot{x}/\dot{x}
\qqtext{i.e.}
	\dlog \tilde{\Q}_0 = \dlog \Q_0 + 2 \dlog \dot{x}
\fullstop
\end{eqntag}
Of course, this identity could have also been deduced directly from $\tilde{\Q}_0 = \Q_0 \dot{x}^2$.
In any case, this identity implies the following Lemma.

\begin{lem}
\label{240815162128}
Given a quadratic differential $\phi_0$ on a marked curve $(X,D)$ and a regular point $p_0 \in X$, there is a canonical local generator $a_0$ of $\omega^\invsqroot_{X}$ near $p_0$, well-defined up to constant scale, given in any local coordinate $x$ near $p_0$ by the following formula:
\begin{eqntag}
\label{240927114929}
	a_0 (x) \coleq \exp \left( - \frac{1}{4} \int^x \dlog \Q_0 \right) \d{x}^\invsqroot
\fullstop
\end{eqntag}
In other words, \eqref{240927114929} is a canonical $\CC^\times$-family of local generators, any of which is called a \dfn{canonical generator} \textit{of $\omega^\invsqroot_X$ near $p_0$}.
The constant scale can be fixed, uniquely up to sign, by choosing the starting point of integration $x_0$ in any coordinate chart:
\begin{eqntag}
\label{240815195953}
	a_0 (x) = \exp \left( - \frac{1}{4} \int_{x_0}^x \dlog \Q_0 \right) \d{x}^\invsqroot
\fullstop
\end{eqntag}
\end{lem}

\begin{proof*}
We need to check that formula \eqref{240927114929} determines a local generator of $\omega^\invsqroot_X$ up to constant scale.
Suppose $\tilde{x}$ is another local coordinate near $p_0$.
Then, using \eqref{240917104315}, we find:
\begin{eqn}
	\exp \left( - \frac{1}{4} \int^{\tilde{x}} \dlog \tilde{\Q}_0 \right) \d{\tilde{x}}^\invsqroot
	= \exp \left( - \frac{1}{4} \int^x \dlog \Q_0 - \tfrac{1}{2} \log \dot{x} \right) \dot{x}^\sqroot \d{x}^\invsqroot
\fullstop
\end{eqn}
This expression is proportional to $a_0 (x)$ because the factor $\exp \big( - \tfrac{1}{2} \log \dot{x} \big) = \dot{x}^\invsqroot$ cancels the line bundle transition function $\dot{x}^\sqroot$.
\end{proof*}

\begin{proof}[Proof of \autoref{240915131722} and \autoref{240915175932}.]
It follows from \autoref{240728090056} and \ref{240815162128} that the local section $\hat{\psi}_\pm$ defined by formula \eqref{240915180049} does not depend on the chosen coordinate $x$ except for an overall multiplicative factor.
Fixing this choice in any local coordinate, fixes it uniformly in all local coordinates.
Therefore, \eqref{240915180049} are well-defined formal WKB solutions.
\end{proof}

\newpage
\subsection{The WKB Differential}

We continue to analyse a Schrödinger operator $\P$ on a marked curve $(X,D)$ satisfying \autoref{240924155219}.

\paragraf
Recall that the perturbative series $\hat{\Y}_\pm (x, \hbar)$ from \autoref{240915175932} define the formal WKB solutions $\hat{\psi}_\pm$ based at $p_0$.
Although, as we explained in \autoref{240728090056}, the coordinate transformation law \eqref{240728091232} means $\hat{\Y}_\pm (x, \hbar) \d{x}$ is not a well-defined differential form, \autoref{240728090056} also tells us that removing the leading- and next-to-leading-order parts of $\hat{\Y}_\pm$ results in a pair of well-defined differential forms near $p_0$.
In fact, they come from two different local determinations of a single global differential form on the spectral curve $\sfSigma$ as follows.
Recall that $R_\sfDelta = R \smallsetminus \sfDelta$ and $\sfDelta = \pi^\ast D \smallsetminus 2R$.

\begin{prop*}
\label{240916114621}
There is a canonical formal meromorphic differential form
\begin{eqntag}
\label{240916090225}
	\hat{\Lambda} (\hbar) = \sum_{k=0}^\infty \Lambda_k \hbar^k \in H^0 \Big(\omega_{\sfSigma,\sfDelta} (\ast R_\sfDelta) \bbrac{\hbar} \Big)
\fullstop{,}
\end{eqntag}
which we call the \dfn{formal WKB differential} on $\sfSigma$ associated with $\P$.
It is the unique formal differential with the property that the two formal WKB solutions based at any regular point $p_0 \in X$ can be expressed as follows:
\begin{eqntag}
	\hat{\psi}_\pm = a_0 e^{-\S_\pm / \hbar} \hat{\A}_\pm
\qtext{where}
	\hat{\A}_\pm (p, \hbar) \coleq
		\exp \left( - \int_{p_0^\pm}^{p_\pm} \hat{\Lambda} (\hbar) \right)
\fullstop{,}
\end{eqntag}
for all $p$ near $p_0$, where $a_0$ is a canonical generator of $\omega^\invsqroot_X$ near $p_0$, and $p_\pm, p_0^\pm \in \sfSigma$ are the appropriately labelled preimages of $p,p_0 \in X$.
Moreover, the $\hbar$-leading-order part of $\hat{\Lambda} (\hbar)$ is the differential form $\Lambda_0 = - \pi^\ast \phi_1 / \sigma \loceq - \Q_1 / 2y_0 \d{x}$ where $\phi_1$ is the quadratic differential from \eqref{240918154226}.
\end{prop*}

The construction of $\hat{\Lambda} (\hbar)$ goes through its coordinate representation as follows.

\begin{lem*}
\label{241005100313}
Explicitly, in a local coordinate $x$ near a regular point $p_0 \in \sfSigma$, the formal WKB differential $\hat{\Lambda} (\hbar)$ can be expressed as follows:
\begin{eqntag}
\label{240918131242}
	\Lambda_0 = y_1 \d{x} - \tfrac{1}{4} \dlog \Q_0
\qqtext{and}
	\Lambda_k = y_{k+1} \d{x}
\quad\text{$(k \geq 1)$}
\end{eqntag}
where the functions $y_k$ are given by the WKB recursion \eqref{240915180751}.
In fact, near $p_0$,
\begin{eqntag}
\label{240918145601}
	\hat{\Y} \d{x} = \lambda + \hbar \tfrac{1}{4} \dlog \Q_0 + \hbar \hat{\Lambda} (\hbar)
\fullstop
\end{eqntag}
\end{lem*}

\begin{proof*}[Proof of \autoref{240916114621} and \autoref{241005100313}.]
By \autoref{240915175932}, the formal WKB solutions based in a local coordinate $x$ are expressed using the uniquely determined (with respect to this coordinate) pair of formal series $\hat{\Y}_\pm (x, \hbar)$.
These series together give different branches of the same formal series $\hat{\Y} (x, \hbar)$ on the spectral curve, though $\hat{\Y}$ still depends on the coordinate.
At the same time, $\lambda \loceq y_0 \d{x}$ and $y_1 = \tfrac{1}{2} \del_x \log y_0 - \Q_1 / 2y_0 = \frac{1}{4} \del_x \log \Q_0 - \Q_1 / 2y_0$.
So identity \eqref{240918145601} defines the series $\hat{\Lambda} (\hbar)$ in the local coordinate $x$ with coefficients given by \eqref{240918131242}.
By \autoref{240728090056}, $\hat{\Lambda} (\hbar)$ is well-defined independent of coordinates.
Thus, we obtain a uniquely determined formal differential form on $\sfSigma_{R,\sfDelta} = \sfSigma \smallsetminus (R \cup \sfDelta)$.

To see that it extends to a global section as in \eqref{240916114621}, we analyse the behaviour of the coefficients $y_k$ at the points of $R \cup \sfDelta$ using the WKB recursion \eqref{240915180751}.
Thus, let $p \in X$ be any point and let $x$ be a coordinate centred at $p$.
Recall that $\sfop{Supp} (\sfDelta) = \pi^\inv \sfop{Pol}_{\geq 2} (\phi_0)$ and $\sfop{Supp} (R_\sfDelta) = \pi^\inv \sfop{Zer} (\phi_0) \cup \pi^\inv \sfop{Pol}_1 (\phi_0)$.

First, we treat the case of simple transition points.
Suppose $p \in \sfop{Zer} (\phi_0)$ is a simple zero of $\phi_0$, and let $\tilde{p} \in R$ be its preimage on $\sfSigma$.
Then we can prove by induction on $k$ using the WKB recursion that for each $k \geq 0$, there is a nonzero constant $c_k$ such that\footnote{Recall that the notation ``$f(x) \sim g(x)$ as $x \to 0$'' means ``$|f(x)|/|g(x)| \to 1$ as $x \to 0$'', whilst ``$f(x) \preceq g(x)$ as $x \to 0$'' means ``$|f(x)|/|g(x)| \to c$ as $x \to 0$ for some $c \geq 0$''.}
\begin{eqntag}
\label{240821172249x}
	y_k (x) \sim c_k x^{(1-3k)/2}
\qquad
	\text{as $x \to 0$}
\fullstop
\end{eqntag}
The spectral curve is locally cut out by the equation $y^2 = x$, so writing the pullback of the differential form $x^{(1-3k)/2} \d{x}$ in the coordinate $y$, we get $2y^{2-3k} \d{y}$.
Therefore, $\Lambda_k$ extends to a meromorphic section of the line bundle $\omega_{\sfSigma, \sfDelta}$ near $\tilde{p}$ with a pole at $\tilde{p}$ of order $2-3k$.
In particular, each $\Lambda_k$ extends to a local section of $\omega_{\sfSigma, \sfDelta} (\ast \tilde{p})$ near $\tilde{p}$.

Next, suppose $p \in \sfop{Pol}_{\geq 2} (\phi_0)$ is a pole of order $m$.
Recall that, if $m = 2n$ is even, then both preimages $p^+,p^-$ of $p$ appear in the divisor $\sfDelta$ with multiplicity $m/2 = n$ each; cf. \eqref{241017173345}.
On the other hand, if $m = 2n - 1$ is odd, the preimage $\tilde{p}$ of $p$ is a ramification point which belongs to the divisor $\sfDelta$ with multiplicity $2n-2$.

If $m = 2n$ is even, then we can prove by induction on $k$, using the WKB recursion and \autoref{240924155219}, that the coefficients $y_k$ have the following behaviour at $p$:
\begin{eqntag}
\label{240821183918}
	y_0 (x) \sim c x^{-n}
\qtext{and}
	y_k (x) \preceq x^{-n}
\qquad
	\text{as $x \to 0$}
\fullstop{,}
\end{eqntag}
for some nonzero complex constant $c$.
It follows that each $\Lambda_k$ extends to a local section of $\omega_{\sfSigma, \sfDelta}$ near $p^\pm$ because it can be written in a local coordinate $x$ near $p^\pm$ as $\Lambda_k (x) = y_{k+1} (x) \d{x} = g_k(x) x^{-n} \d{x}$ for some holomorphic function $g_k$ near $p^\pm$.

If $m = 2n - 1 \geq 3$ is odd, then we can prove again by induction that
\begin{eqntag}
\label{240821185032}
	y_k (x) \preceq x^{-n+1/2}
\qquad
	\text{as $x \to 0$}
\fullstop
\end{eqntag}
The spectral curve is locally cut out by the equation $y^2 = x$, so writing the pullback of the differential form $x^{-n+1/2} \d{x}$ in the coordinate $y$ yields $(y^2)^{-n+1/2} \d{(y^2)} = 2y^{2(n-1)} \d{y}$.
It follows that each $\Lambda_k$ again extends to a local section of $\omega_{\sfSigma, \sfDelta}$ near $\tilde{p}$ because it can be written in a local coordinate $y$ near $\tilde{p}$ as $\Lambda_k (y) = y_{k+1} (x) \d{x} = g_k (y) 2y^{2(n-1)} \d{y}$ for some holomorphic function $g_k$ near $\tilde{p}$.

Finally, if $m = 1$, so that $p$ is a simple pole and hence $\tilde{p}$ is a non-simple transition point, then we can prove by induction that
\begin{eqntag}
\label{241018180606}
	y_0 (x) \sim c x^{-1/2}
\qtext{and}
	y_k (x) \preceq x^{-(k+1)/2}
\qquad
	\text{as $x \to 0$}
\fullstop{,}
\end{eqntag}
for some nonzero complex constant $c$.
The spectral curve is again locally cut out by the equation $y^2 = x$, so if we write the pullback of the differential form $x^{-(k+1)/2} \d{x}$ in the coordinate $y$, we get $2y^{-k} \d{y}$.
So $\Lambda_k$ extends to a meromorphic section of the line bundle $\omega_{\sfSigma, \sfDelta}$ near $\tilde{p}$ with a pole at $\tilde{p}$ of order $k$.
In particular, each $\Lambda_k$ extends to a local section of $\omega_{\sfSigma, \sfDelta} (\ast \tilde{p})$ near $\tilde{p}$.
\end{proof*}

\paragraf
Recall the canonical differential form $\sigma = \lambda - \iota^\ast \lambda \in H^0 (\sfSigma, \omega_{\sfSigma, \sfDelta})$.
Let $\V$ be the unique meromorphic vector field on the spectral cover $\sfSigma$ which is dual to $\sigma$ in the sense that $\sigma (\V) = 1$.
Explicitly in a local coordinate $x$ near a regular point on $\sfSigma$, we have
\begin{eqn}
\label{240126122452}
	\sigma \loceq 2y_0 (x) \d{x} \qquad \xleftrightarrow{~~~ \sigma (\V) = 1 ~~~} \qquad \V \loceq \frac{1}{2y_0 (x)} \del_x
\fullstop
\end{eqn}

\begin{prop*}
\label{240915182325}
In terms of the canonical differential $\sigma$, the formal WKB differential $\hat{\Lambda}$ is given by the identity
\begin{eqntag}
\label{240916140840}
	\hat{\Lambda} = \Lambda_0 + \hat{f} \sigma
\end{eqntag}
where $\hat{f} \in H^0 \big(\cal{O}_{\sfSigma} (\ast R_\sfDelta) \bbrac{\hbar} \big)$ is the unique power series solution of the Riccati equation
\begin{eqntag}
\label{240916141300}
	\hbar \V (f) - f = \hbar (f^2 + w f + \W)
\end{eqntag}
whose coefficients $w$ and $\W$ are explicitly given by
\begin{eqntag}
\label{240918160044}
	w \coleq \frac{\Lambda_0}{\lambda} = -2 \frac{\pi^\ast \phi_1}{\sigma^2}
\qqtext{and}
	\W \coleq \frac{\Lambda_1}{\sigma} - \sum_{k \geq 1} \frac{\pi^\ast \phi_{k+2} }{\sigma^2} \hbar^{k}
\fullstop{,}
\end{eqntag}
where $\phi_k \in H^0 (X, \omega_{X,D}^2)$ are the quadratic differentials associated with the given Schrödinger operator $\P$.
\end{prop*}

We emphasise that, as opposed to a local Riccati equation \eqref{240815150746}, the Riccati equation \eqref{240916141300} is an equation for a function $f$ on $\sfSigma$ which makes global coordinate-invariant sense.
Note also that $\W = \W_0 + \W_1 \hbar + \cdots$ is a polynomial in $\hbar$, and its coefficients $\W_k$ as well as the function $w$ are global meromorphic functions on $\sfSigma$ with simple poles along the transition locus:
\begin{eqntag}
\label{241018183300}
	w, \W_k \in H^0 \big(\sfSigma, \cal{O}_\sfSigma (R_\sfDelta) \big)
\fullstop
\end{eqntag}
Explicitly, in any local coordinate $x$ near a regular point on $\sfSigma$,
\begin{eqntag}
\label{240918173507}
	w (x) = - \frac{\Q_1}{2 y_0^2}
\qtext{and}
	\W (x, \hbar) = - \frac{y_2}{2y_0} - \sum_{k \geq 1} \frac{\Q_{k+2}}{(2y_0)^2}
\fullstop
\end{eqntag}

\begin{proof*}[Proof of \autoref*{240915182325}]
The Riccati equation \eqref{240916141300} can be derived from a local Riccati equation \eqref{240815150746} by combining identities \eqref{240918145601} and \eqref{240916140840} into a single change of the unknown variable.
Namely, let $x$ be a local coordinate near any regular point on $\sfSigma$, and make the following change of the unknown variable in the Riccati equation \eqref{240815150746}:
\begin{eqntag}
\label{240728214711}
	\Y (x, \hbar) ~\to~ f (x, \hbar)
\qtext{given by}
	\Y = y_0 + y_1 \hbar + 2y_0 \hbar f
\fullstop
\end{eqntag}
This leads to a Riccati equation for $f(x, \hbar)$ which, after some simplification, reads as follows:
\begin{eqntag}
\label{241018183948}
	\frac{\hbar}{2y_0} \del_x f - f = \hbar \left( f^2 + w(x) f - \frac{y_2}{2y_0} + \W (x,\hbar) \right)
\fullstop{,}
\end{eqntag}
where $w(x)$ and $\W(x,\hbar)$ are given precisely by \eqref{240918173507}.
We can also recognise the vector field on the lefthand side as the coordinate representation of the vector field $\V$.
Therefore, we can canonically extend this Riccati equation to a global equation on the spectral curve $\sfSigma$ because both the vector field $\V$ and the coefficients $w, \W$ make sense on $\sfSigma$.
\end{proof*}

\section{Resurgence of WKB Solutions}
\label{240801153058}

We can now state and prove the main results of this paper.
We begin by briefly recalling our setup, notation, and some facts that will facilitate the description of resurgence properties of formal WKB solutions.

\paragraph{Recollection.}
Let $(X,D)$ be a marked curve with a chosen spin structure $\omega^\invsqroot_X$.
We consider a Schrödinger operator $\P$ on $(X,D)$ in the sense of \autoref{240730095844}.
The locus of turning points $\sfop{Turn} (\phi_0)$ is the set of zeros and simple poles of the  quadratic differential $\phi_0 \in \omega^2_{X,D}$ associated with $\P$.
Let $\sfSigma \subset T^\ast_{X,D}$ be the associated spectral curve with Liouville form $\lambda \in \omega_{\sfSigma, \sfDelta}$, where $\sfDelta = \pi^\ast D \smallsetminus 2R$ and $R \subset \sfSigma$ is the ramification locus of $\pi : \sfSigma \to X$.
The divisor $\sfDelta$ is supported at the preimages of poles of $\phi_0$ of order $2$ or greater.
The preimage of turning points $\sfop{Turn} (\phi_0)$ is the locus of transition points $\sfop{Tran} (\phi_0)$, which is the support of the divisor $R_\sfDelta = R \smallsetminus \sfDelta$.

Now, consider the \hyperref[241021121815]{fundamental groupoid} $\Pi_1 ({\sfSigma_\sfDelta})$ of the punctured spectral curve ${\sfSigma_\sfDelta} = \sfSigma \smallsetminus \sfDelta$.
It is a two-dimensional holomorphic manifold equipped with surjective submersions $\rm{s},\rm{t} : \Pi_1 ({\sfSigma_\sfDelta}) \to \sfSigma_\sfDelta$.
If $p \in \sfSigma$ is a regular point, we consider the \hyperref[241021115044]{Borel surface} associated with $p$, which is the groupoid source fibre
\begin{eqn}
	\tilde{\sfSigma}_{\sfDelta,p} \coleq \rm{s}^\inv (p)
\fullstop
\end{eqn}
It is the universal cover of $\sfSigma_\sfDelta$ based at $p$, hence a simply connected Riemann surface with a distinguished basepoint given by the constant path $1_p \in \tilde{\sfSigma}_{\sfDelta,p}$.

Let $\Gamma_{p} \subset \tilde{\sfSigma}_{\sfDelta,p}$ be the countable discrete subset of \hyperref[240725152318]{critical paths} starting at $p$.
These are paths in ${\sfSigma_\sfDelta}$ that start at $p$ and terminate at a transition point.
If we puncture the Borel surface $\tilde{\sfSigma}_{\sfDelta,p}$ at the points of $\Gamma_p$, its universal cover based at the origin $1_p$ is canonically isomorphic to the universal cover $\tilde{\sfSigma}_{R, \sfDelta,p}$ of $\sfSigma_{R,\sfDelta} = \sfSigma \smallsetminus (R \cup \sfDelta)$, the spectral curve with basepoint $p$ punctured along both $\sfDelta$ and $R$:
\begin{eqn}
	\tilde{\sfSigma}_{R, \sfDelta,p} \cong \widetilde{\tilde{\sfSigma}_{\sfDelta,p} \!\! \smallsetminus \! \Gamma_{p}} \to \tilde{\sfSigma}_{\sfDelta,p} \smallsetminus \Gamma_{p}
\fullstop
\end{eqn}

We also consider the one-form $\sigma \coleq \lambda - \iota^\ast \lambda$ on $\sfSigma$, where $\iota : \sfSigma \to \sfSigma$ is the involution, and define the central charge $\Z : \Pi_1 ({\sfSigma_\sfDelta}) \to \CC$ by $\Z (\gamma) = \int_\gamma \sigma$.
Critical trajectories are critical paths which are projected under $\Z$ to straight lines in $\CC$.

If $q_0 \in X$ is a regular point, pick one of its preimages $p_0 \in \sfSigma$, and let $\hat{\psi}$ be a formal WKB solution based at $p_0$:
\begin{eqntag}
\label{240918184406}
	\hat{\psi} = a_0 e^{-\S / \hbar} \hat{\A}
\qtext{where}
	\S (q) = \int_{p_0}^{p} \lambda \in \cal{O}_{X,q_0}
\fullstop{,}\quad
	\hat{\A} \in \cal{O}_{X,p_0} \bbrac{\hbar}
\fullstop{,}
\end{eqntag}
for all $p \in \sfSigma$ near $p_0$ with projection $q = \pi (p)$.
Here, $a_0$ is a canonical (uniquely determined up to constant scale) local generator of $\omega^\invsqroot_X$ near $q_0$; cf. \autoref{240815162128}.

\subsection{Main Results}
\label{240918202409}

The main result of this paper can be formulated as follows.

\begin{theorem}[resurgence of WKB solutions]
\label{240801161941}
Consider a Schrödinger operator $\P$ on a marked curve $(X,D)$ satisfying \autoref{240924155219}.
Fix a regular point $q_0 \in X$ and pick one of its preimages $p_0 \in \sfSigma$.
Then the formal WKB solution $\hat{\psi}$ based at $p_0 \in \sfSigma$ is resurgent of log-algebraic type, uniformly for all $p$ near $p_0$, in all but the unstable directions at $p_0$.
\end{theorem}

We now break this Theorem down into several parts, each presented as a separate Proposition that elaborates on a specific aspect of resurgence of WKB solutions.
Namely, \autoref{241019111435} demonstrates the convergence of the Borel transform, \autoref{241019113314} describes the Borel transform's endless analytic continuation, \autoref{241019115130} clarifies the relationship between the resurgent and the trajectory-based Stokes diagrams, whilst Propositions \ref{240923104232}, \ref{241019122801}, and \ref{241022040824} spell out the Borel summability properties.
The proof of \autoref{240801161941} amounts to proving these Propositions.

The proof of these Propositions involves several steps, organised into the subsections below.
The first step is to reformulate the problem in terms of the resurgence properties of the formal WKB differential $\hat{\Lambda} (\hbar)$ which are of interest in their own right.
The advantage of this point of view is offered by \autoref{240915182325} which allows us to convert the given Schrödinger equation into a global coordinate-invariant Riccati equation on the spectral curve $\sfSigma$.
So, in \autoref{240815131241}, we state an analogous \autoref{240918164353} that asserts the resurgence properties of the formal WKB differential.

\autoref*{240918164353} is then similarly broken down into several Propositions.
Namely, the list of implications between these Propositions is as follows: \ref{241019125911} $\Rightarrow$ \ref{241019111435}, \ref{241019132652} $\Rightarrow$ \ref{241019113314}, \ref{241019134004} $\Rightarrow$ \ref{241019115130}, \ref{241019135458} $\Rightarrow$ \ref{240923104232}, \ref{241019144520} $\Rightarrow$ \ref{241019122801}, and finally \ref{241019145013} $\Rightarrow$ \ref{241022040824}.
As a result and with the help of \autoref{240918203006}, the proof of \autoref{240801161941} gets reduced to proving Propositions \ref{241019125911}-\ref{241019145013}, which we will begin in \autoref{240815131241}.

\paragraf
We start by spelling out the convergence of the Borel transform of $\hat{\A}$, and how it can be canonically lifted to the groupoid.

\begin{prop*}[Borel transform]
\label{241019111435}
Assume the hypothesis of \autoref{240801161941}.
Then the formal power series $\hat{\A} (q, \hbar) \in \CC \bbrac{\hbar}$ is of factorial type, uniformly for all $q$ near $q_0$.
In other words, the Borel transform of $\hat{\A} (q, \hbar)$ in the variable $\hbar$ is a convergent power series in $t$, uniformly for all $q$ near $q_0$:
\begin{eqntag}
\label{240916182534}
	\text{$\hat{\Phi} (t) \in \mathcal{O}_{X,q_0} \cbrac{t}$}
\qtext{where}
	\hat{\Phi} (q, t) \coleq \Borel \big[ \, \hat{\A} \, \big] (q, t) \in \CC \cbrac{t}
\fullstop
\end{eqntag}
More concretely, $p_0$ has a neighbourhood $U \subset \sfSigma$ such that $\hat{\Phi} (t) \in \mathcal{O}_\sfSigma (U) \cbrac{t}$.
Furthermore, for every $p \in U$, the convergent power series $\hat{\Phi} (p, t)$ can be canonically regarded, using the central charge $\Z_{p}$ and the projection $\pi$, as a holomorphic germ $\hat{\Phi} (p)$ at the point $1_{p}$ of the Borel surface $\tilde{\sfSigma}_{\sfDelta,p}$, which is the source fibre $\rm{s}^\inv (p)$ in the fundamental groupoid $\Pi_1 (\sfSigma_\sfDelta)$.
Equivalently, $\hat{\Phi} (p)$ is a holomorphic germ at the point $1_{p}$ of the source fibre $\tilde{\sfSigma}_{R,\sfDelta,p}$ in the fundamental groupoid $\Pi_1 (\sfSigma_{R,\sfDelta})$.
Consequently, the Borel transform of $\hat{\A}$ is canonically a holomorphic germ $\hat{\Phi}$ at the point $1_{p_0} \in \Pi_1 (\sfSigma_\sfDelta)$ or equivalently at the point $1_{p_0} \in \Pi_1 (\sfSigma_{R,\sfDelta})$.
\end{prop*}

The precise relationship between the holomorphic germs
\begin{eqntag}
\label{241019072223}
	\hat{\Phi} (q,t) \in \CC \cbrac{t}
\qqtext{and}
	\hat{\Phi} (p) \in \cal{O}_{\tilde{\sfSigma}_{\sfDelta, p}, 1_{p}}	
\end{eqntag}
is as follows.
By \autoref{240328174841}, the central charge $\Z_{p} : \tilde{\sfSigma}_{\sfDelta,p} \to \CC$ restricts to a biholomorphism from a neighbourhood of the point $1_{p}$ to a neighbourhood of the origin in $\CC$.
It induces an isomorphism
\begin{eqntag}
\label{241019074127}
	\Z_{p} : \CC \cbrac{t} \too \cal{O}_{\tilde{\sfSigma}_{\sfDelta,p}, 1_{p}}
\qtext{and so}
	\hat{\Phi} (p) = \Z_{p}^\ast \hat{\Phi} (q,t)
\fullstop
\end{eqntag}
This identity is true uniformly at all points $p$ near $p_0$, which may be rephrased as follows.
First, $\hat{\Phi} (q,t)$ is a convergent power series in $t$ uniformly for all $q$ near $q_0$, meaning it defines a holomorphic germ $\hat{\Phi} (t)$ at the point $(q_0,0) \in X \times \CC$, so 
\begin{eqntag}
	\hat{\Phi} (t) \in \cal{O}_{X,q_0} \cbrac{t} \cong \cal{O}_{\sfSigma,p_0} \cbrac{t}
\fullstop{,}
\end{eqntag}
where we have made use of the isomorphism $\pi^\ast : \cal{O}_{X,q_0} \iso \cal{O}_{\sfSigma, p_0}$.
Meanwhile, the fundamental groupoid $\Pi_1 (\sfSigma_\sfDelta)$ is foliated by source fibres, so in particular a neighbourhood of the point $1_{p_0} \in \Pi_1 (\sfSigma_\sfDelta)$ is foliated by open subsets of source fibres $\tilde{\sfSigma}_{\sfDelta,p}$ parameterised by points $p$ near $p_0$.
In turn, each of these open subsets are identified by the central charge $\Z_{p}$ with open subsets of the origin in $\CC$.
As a result, we may view the central charge $\Z$ as defining a biholomorphism from an open neighbourhood of the point $1_0 \in \Pi_1 (\sfSigma_\sfDelta)$ to an open neighbourhood of the point $(q_0, 0) \in X \times \CC$, yielding an isomorphism $\Z^\ast : \cal{O}_{\sfSigma, p_0} \set{t} \iso \cal{O}_{\Pi_1 (\sfSigma_\sfDelta), 1_{p_0}}$.
Altogether, we may write
\begin{eqntag}
\label{241019081222}
	\hat{\Phi} = \Z^\ast \pi^\ast \hat{\Phi} (t)
\fullstop
\end{eqntag}

\begin{prop}[endless continuation]
\label{241019113314}
Assume the hypothesis of \autoref{240801161941}.
Then $p_0$ has a neighbourhood $U \subset \sfSigma$ such that the Borel transform $\hat{\Phi} (q,t) \in \CC \cbrac{t}$ admits endless analytic continuation of log-algebraic type for every $q \in \pi(U) \subset X$.
Namely, let us regard $\hat{\Phi} (q,t) \in \CC \cbrac{t}$ as the holomorphic germ $\hat{\Phi} (p)$ at the point $1_{p}$ of the Borel surface $\tilde{\sfSigma}_{\sfDelta,p} \subset \Pi_1 (\sfSigma_\sfDelta)$ as explained in \autoref{241019111435}.
Then $\hat{\Phi} (p)$ can be analytically continuation along any path in the Borel surface $\tilde{\sfSigma}_{\sfDelta,p}$ that starts at $1_{p}$ and avoids the discrete subset $\Gamma_{p} \subset \tilde{\sfSigma}_{\sfDelta,p}$ consisting of all critical paths in $\sfSigma_\sfDelta$ starting at $p$.
In other words, the germ $\hat{\Phi} (p)$ extends to a global holomorphic function $\Phi (p)$ on the universal cover of the punctured Borel surface $\tilde{\sfSigma}_{\sfDelta,p} \smallsetminus \Gamma_{p}$ based at the origin $1_{p}$, which is canonically isomorphic to the universal cover $\tilde{\sfSigma}_{R, \sfDelta,p}$ of the punctured spectral curve $\sfSigma_{R,\sfDelta}$ based at $p$:
\begin{eqn}
	\tilde{\sfSigma}_{R, \sfDelta,p} 
		\cong \widetilde{\tilde{\sfSigma}_{\sfDelta,p} \!\!\smallsetminus\! \Gamma_{p}}
		\to \tilde{\sfSigma}_{\sfDelta,p} \smallsetminus \Gamma_{p}
\fullstop
\end{eqn}
Consequently, the holomorphic germ $\hat{\Phi}$ at $1_{p_0} \in \Pi_1 (\sfSigma_{R,\sfDelta})$ extends to a holomorphic function $\Phi$ on the open submanifold
\begin{eqn}
	\rm{s}^\inv (U) 
		= \mathop{\Coprod}_{p \in U} \tilde{\sfSigma}_{R, \sfDelta,p} 
		\subset \Pi_1 (\sfSigma_{R,\sfDelta})
\fullstop{,}
\end{eqn}
which is a $U$-family of endless Riemann surfaces of log-algebraic type.
\end{prop}

In other words, \autoref*{241019113314} says that the Borel singularities of the formal WKB solution $\hat{\psi} (q, \hbar)$ are properly understood as a discrete subset of the Borel surface $\tilde{\sfSigma}_{\sfDelta,p}$ instead of the traditional Borel $t$-plane $\CC$.
Meanwhile, the set of Borel singular \textit{values} $\Xi_p \subset \CC$ in the traditional Borel $t$-plane is formed by the central charges of all critical paths starting at $p$:
\begin{eqn}
	\Xi_p = \Z (\Gamma_p) \subset \CC
\end{eqn}
Note that whilst $\Gamma_p$ is discrete, the subset $\Xi_p$ may well have accumulation points.

\begin{prop}[resurgent Stokes diagram]
\label{241019115130}
Assume the hypothesis of \autoref{240801161941}.
Let $U \subset \sfSigma$ be a neighbourhood of $p_0$ such that $\hat{\Phi} (t) \in \mathcal{O}_\sfSigma (U) \cbrac{t}$.
Then, for every $p \in U$, the resurgent Stokes diagram for $\hat{\Phi} (p,t)$ coincides with the trajectory-based Stokes diagram at $p$.
Namely:
\begin{enumerate}
\item \textup{\bfseries(visible singularities).}
The set of visible singularities of $\hat{\Phi} (p,t)$ is the subset $\Gamma_{p}^{0} \subset \Gamma_{p}$ of all critical trajectories starting at $p$.
\item \textup{\bfseries(Stokes rays).}
The Stokes rays for $\hat{\Phi} (p,t)$ are the critical phases at $p$; the stability of a Stokes ray matches the stability of a critical phase.
\item \textup{\bfseries(regular rays).}
The regular rays for $\hat{\Phi} (p,t)$ are the regular phases at $p$; the stability of a regular ray matches the stability of a regular phase.
\item \textup{\bfseries(observable singularities).}
Finally, if $\alpha$ is a Stokes ray for $\hat{\Phi} (p,t)$ and $\gamma \in \Gamma_{p}^{0}$ is the corresponding critical trajectory with phase $\alpha$ (i.e., the visible singularity with phase $\alpha$), then the subset $\Gamma_{p,\alpha}^{+} \subset \Gamma_{p, \alpha}$ of observable singularities with phase $\alpha$ is in one-to-one correspondence with finite chains of saddle trajectories with phase $\alpha$ starting at the transition point where $\gamma$ terminates.
\end{enumerate}
\end{prop}

\paragraf
Next, we spell out in much greater detail the property of exponential type and hence Borel summability of WKB solutions.
Namely, we first focus on the Borel summability in regular directions before discussing Stokes directions.
For clarity and ease of application, we treat the cases of summability in a single direction and along an arc separately because of the subtle differences in assumptions and notation.

\begin{prop*}[Borel summability of WKB solutions in a direction]
\label{240923104232}
Assume the hypothesis of \autoref{240801161941}.
Suppose $\alpha$ is a (semi-)stable regular ray at $p_0$.
Then the formal WKB solution $\hat{\psi}$ is (semi-)stably Borel summable with phase $\alpha$, uniformly near $q_0$.
Namely, there is a domain $U \subset X$ around $q_0$ such that the Borel resummation
\begin{eqntag}
\label{241012184755}
	\psi_{\alpha}
		\coleq s_\alpha \big[\, \hat{\psi}_{\alpha} \,\big]
		= a_0 e^{-\S / \hbar} \A_{\alpha}
\qtext{where}
	\A_{\alpha} \coleq s_\alpha \big[\, \hat{\A} \,\big]
\end{eqntag}
is well-defined and (semi-)stable uniformly on $U$.
In particular, $\psi_{\alpha} \in \omega^\invsqroot_X$ is a holomorphic solution of the Schrödinger equation on the domain $U \times V$ where $V$ is a sector at the origin with aperture $\sfop{Arc} (\alpha)$.
Moreover, as $\hbar \to 0$ in this sector, $\psi_{\alpha}$ is asymptotic to $\hat{\psi}$ with uniform factorial growth, uniformly for all $q \in U$:
\begin{eqntag}
\label{241012184751}
	\psi_{\alpha} (q, \hbar) \simeq \hat{\psi} (q, \hbar)
\qquad
\text{as $\hbar \to 0$ unif. along $\sfop{Arc} (\alpha)$, unif. $\forall q \in U$}
\fullstop
\end{eqntag}
In fact, $\psi_{\alpha}$ is the unique holomorphic section of $\omega^\invsqroot_X$ over $U \times V$ which satisfies \eqref{241012184751}.
\end{prop*}

\begin{prop*}[Borel summability of WKB solutions along an arc]
\label{241019122801}
Assume the hypothesis of \autoref{240801161941}.
Suppose $A$ is an arc of regular rays at $p_0$.
Then the formal WKB solution $\hat{\psi}$ is Borel summable along $A$, uniformly near $q_0$.
Namely, there is a domain $U \subset X$ around $q_0$ such that the Borel resummation
\begin{eqntag}
\label{241012184831}
	\psi_{A} 
		\coleq s_A \big[\, \hat{\psi}_{A} \,\big]
		= a_0 e^{-\S / \hbar} \A_{A}
\qtext{where}
	\A_{A} \coleq s_A \big[\, \hat{\A} \,\big]
\end{eqntag}
is well-defined uniformly on $U$.
In particular, $\psi_{A} \in \omega^\invsqroot_X$ is a holomorphic solution of the Schrödinger equation on $U \times V$ where $V$ is a sector at the origin with aperture $\sfop{Arc} (A)$.
Moreover, as $\hbar \to 0$ in this sector, $\psi_{A}$ is asymptotic to $\hat{\psi}$ with factorial growth uniformly on $U$:
\begin{eqntag}
\label{241012184833}
	\psi_{A} (q, \hbar) \simeq \hat{\psi} (q, \hbar)
\qquad
\text{as $\hbar \to 0$ along $\sfop{Arc} (A)$, unif. $\forall q \in U$}
\fullstop
\end{eqntag}
In fact, $\psi_{A}$ is the unique holomorphic section of $\omega^\invsqroot_X$ over $U \times V$ which satisfies \eqref{241012184833}.
\end{prop*}

\paragraf
Lastly, we describe the lateral Borel summability of WKB solutions in Stokes directions.

\begin{prop*}[Lateral Borel summability of WKB solutions]
\label{241022040824}
Assume the hypothesis of \autoref{240801161941}.
Fix a phase $\alpha$ and a path $\wp$ on $\sfSigma_{R,\sfDelta}$ which starts at $p_0$ and terminates at some point $p_\ast$, and put $q_\ast \coleq \pi (p_\ast)$.
Suppose $\alpha$ is a (semi-)stable regular ray at every point $p \in \wp$ except $p_\ast$ where instead $\alpha$ is a (semi-)stable Stokes ray.
Then the formal WKB solution $\hat{\psi} (q_\ast, \hbar)$ is laterally Borel summable with phase $\alpha$ in either the left or right sense depending on the direction that $\wp$ crosses the critical trajectory with phase $\alpha$ passing through $p_\ast$:
\begin{eqntag}
\label{241012191457}
	\psi_{\alpha}^{\rm{L/R}} (q_\ast,\hbar)
		\coleq s_\alpha^{\rm{L/R}} \big[\, \hat{\psi}_{\alpha} \,\big] (q_\ast,\hbar)
		= a_0 (q_\ast) e^{-\S (q_\ast) / \hbar} \A_{\alpha}^{\rm{L/R}} (q_\ast, \hbar)
\fullstop{,}
\end{eqntag}
where $\A_{\alpha}^{\rm{L/R}} (q_\ast, \hbar) \coleq s_\alpha^{\rm{L/R}} \big[\, \hat{\A} \,\big] (q_\ast, \hbar)$.
In particular, $\A_{\alpha}^{\rm{L/R}} (q_\ast, \hbar) \in \cal{O} (V)$ is a holomorphic function in a sector $V$ at the origin with aperture $\sfop{Arc} (\alpha)$.
Moreover, as $\hbar \to 0$ in this sector, $\psi_{\alpha}^{\rm{L/R}} (q_\ast)$ is asymptotic to $\hat{\psi} (q_\ast)$ with factorial growth:
\begin{eqntag}
\label{241013135057}
	\psi_{\alpha}^{\rm{L/R}} (q_\ast,\hbar) \simeq \hat{\psi} (q_\ast, \hbar)
\qquad
\text{as $\hbar \to 0$ along $\sfop{Arc} (\alpha)$}
\fullstop
\end{eqntag}
\end{prop*}

We stress that $q_\ast$ in the equality \eqref{241012191457} and the asymptotic equivalence \eqref{241013135057} is fixed.
We also emphasise that the asymptotic equivalence \eqref{241013135057} does not hold \textit{uniformly} along the arc $\sfop{Arc} (\alpha)$, in contrast to the asymptotic equivalence \eqref{241012184751}.

\paragraf
Finally, we remark that if $\alpha$ is an unstable direction at $p_0$, then the only property of resurgence that fails is the exponential type.
Unstable directions correspond to divergent trajectories of quadratic differentials.
Consequently, in the absence of divergent trajectories, we can immediately conclude a stronger form of resurgence.

\begin{cor*}[Strong resurgence of WKB solutions]
\label{241007171114}
Consider a Schrödinger operator $\P$ on a marked curve $(X,D)$ satisfying \autoref{240924155219}.
Assume, in addition, that the associated quadratic differential $\phi_0$ has no divergent trajectories for any phase.
Then the formal WKB solution $\hat{\psi}$ based at any regular point $p_0 \in \sfSigma$ is strongly resurgent of log-algebraic type, uniformly for all $p$ near $p_0$.
\end{cor*}

\subsection{Resurgence of the WKB Differential}
\label{240815131241}

The first step in the proof of \autoref{240801161941} is to reformulate the problem in terms of the formal WKB differential $\hat{\Lambda}$.
In this subsection, we formulate resurgence properties $\hat{\Lambda}$ which are of great interest in their own right.
We have deliberately kept the structure and the phrasing in this subsection close to that of the previous subsection in order to make it clear that resurgence of the formal WKB differential implies the resurgence of the formal WKB solutions.

\begin{theorem}[Resurgence of the WKB differential]
\label{240918164353}
Consider a Schrödinger operator $\P$ on a marked curve $(X,D)$ satisfying \autoref{240924155219}.
Then the formal WKB differential $\hat{\Lambda}$ is resurgent of log-algebraic type, locally uniformly on the subset $\sfSigma_{R,\sfDelta}$ of regular points of the spectral curve, in all but the unstable directions.
\end{theorem}

As with \autoref*{240801161941}, we now break this Theorem down into several parts, each presented as a separate Proposition that elaborates on a specific aspect of resurgence of the formal WKB differential.
Namely, \autoref{241019125911} demonstrates the convergence of the Borel transform, \autoref{241019132652} describes the Borel transform's endless analytic continuation, \autoref{241019134004} clarifies the relationship between the resurgent and the trajectory-based Stokes diagrams, whilst Propositions \ref{241019135458}, \ref{241019144520}, and \ref{241019145013} spell out the Borel summability properties.
The proof of \autoref{240918164353} amounts to proving these Propositions.

Let us write $\hat{\Lambda} = \Lambda_0 + \hat{f} \sigma$ and describe the resurgent structure of $\hat{\Lambda}$ in terms of the resurgent structure of $\hat{f}$.

\begin{prop}[Borel transform]
\label{241019125911}
Assume the hypothesis of \autoref{240918164353}.
Then the formal power series $\hat{f} (p, \hbar) \in \CC \bbrac{\hbar}$ is of factorial type, locally uniformly for all $p \in \sfSigma_{R,\sfDelta}$.
In other words, the Borel transform of $\hat{f} (p, \hbar)$ in the variable $\hbar$ is a convergent power series in $t$, locally uniformly for all $p \in \sfSigma_{R,\sfDelta}$:
\begin{eqntag}
\label{240918203538}
	\hat{\phi} (t) \in \text{$\cal{O}_{\sfSigma_{R,\sfDelta}} \cbrac{t}$}
\qqtext{where}
	\hat{\phi} (p, t) \coleq \text{$\mathfrak{B}$} \big[ \, \hat{f} \, \big] (t) \in \CC \cbrac{t}
\fullstop
\end{eqntag}
More concretely, every regular point has a neighbourhood $U \subset \sfSigma_{R,\sfDelta}$ such that the Borel transform defines an element $\hat{\phi} (t) \in \cal{O}_{\sfSigma} (U) \cbrac{t}$.
Furthermore, for every $p \in \sfSigma_{R,\sfDelta}$, the convergent power series $\hat{\phi} (p, t) \in \CC \cbrac{t}$ can be canonically regarded, using the central charge $\Z_p$, as a holomorphic germ $\hat{\phi} (p)$ at the point $1_p$ of the source fibre $\tilde{\sfSigma}_{\sfDelta, p}$ of the fundamental groupoid $\Pi_1 (\sfSigma_\sfDelta)$, or equivalently at the point $1_p$ of the source fibre  $\tilde{\sfSigma}_{R,\sfDelta, p}$ of the fundamental groupoid $\Pi_1 (\sfSigma_{R,\sfDelta})$.
Consequently, the Borel transform of $\hat{f}$ is canonically a holomorphic germ $\hat{\phi}$ at the point $1_{p} \in \Pi_1 (\sfSigma_\sfDelta)$ for all $p \in \sfSigma_{R,\sfDelta}$; i.e., $\hat{\phi}$ is a holomorphic germ along the identity bisection $\sfSigma_{R,\sfDelta} \subset \Pi_1 (\sfSigma_{R,\sfDelta})$.
\end{prop}

\begin{proof*}
Recall that by \autoref{240915182325}, the formal series $\hat{f}$ is the unique $\hbar$-power series solution of the Riccati equation
\begin{eqntag}
\label{240919090251}
	\hbar \V (f) - f = \hbar (f^2 + w f + \W)
\fullstop
\end{eqntag}
One can derive a recursive formula for the coefficients of the formal power series $\hat{f}$, very similar to the WKB recursion \eqref{240915180751}.
This recursion can then be used to prove that the power series 
\begin{eqn}
	\hat{f} (p, \hbar) = \sum_{k=1}^\infty f_k (p) \hbar^k
\end{eqn}
is of factorial type uniformly near regular points on $\sfSigma$.
That is, any regular point has a neighbourhood $U \subset \sfSigma$ such that there is a real constant $\M > 0$ which provides the bound $\big| f_k (p) \big| \leq \M^k k!$ for all $p \in U$ and all $k \geq 0$.
Details can be found in, e.g., \cite[Corollary 3.12]{MY2008.06492}.
It follows now that the Borel transform $\hat{\phi}$ of $\hat{f}$ is uniformly convergent; i.e., $\hat{\phi} \in \cal{O}_\sfSigma (U) \cbrac{t}$.

Now, let us spell out in detail the relationship between the holomorphic germs
\begin{eqntag}
\label{241019155122}
	\hat{\phi} (p,t) \in \CC \cbrac{t}
\qqtext{and}
	\hat{\phi} (p) \in \cal{O}_{\tilde{\sfSigma}_{\sfDelta, p}, 1_{p}}
\fullstop
\end{eqntag}
By \autoref{240328174841}, the central charge $\Z_{p} : \tilde{\sfSigma}_{\sfDelta,p} \to \CC$ restricts to a biholomorphism from a neighbourhood of the point $1_{p}$ to a neighbourhood of the origin in $\CC$.
It induces an isomorphism
\begin{eqntag}
\label{241019155424}
	\Z_{p} : \CC \cbrac{t} \too \cal{O}_{\tilde{\sfSigma}_{\sfDelta,p}, 1_{p}}
\qtext{and so}
	\hat{\phi} (p) = \Z_{p}^\ast \hat{\phi} (p,t)
\fullstop{,}
\end{eqntag}
and this identity is true locally uniformly for all $p \in \sfSigma_{R,\sfDelta}$.
\end{proof*}

\begin{prop}[endless continuation]
\label{241019132652}
Assume the hypothesis of \autoref{240918164353}.
Then the Borel transform $\hat{\phi} (p,t) \in \CC \cbrac{t}$ admits endless analytic continuation of log-algebraic type for every $p \in \sfSigma_{R,\sfDelta}$.
Namely, let us regard $\hat{\phi} (p,t)$ as the holomorphic germ $\hat{\phi} (p)$ at $1_p \in \tilde{\sfSigma}_{\sfDelta, p} \subset \Pi_1 (\sfSigma_\sfDelta)$ as explained in \autoref{241019125911}.
Then $\hat{\phi} (p)$ can be analytically continued along any path in the source fibre $\tilde{\sfSigma}_{\sfDelta, p}$ that starts at $1_p$ and avoids the discrete subset $\Gamma_p \subset \tilde{\sfSigma}_{\sfDelta, p}$ consisting of critical paths in $\sfSigma_\sfDelta$ starting at $p$.
In other words, the germ $\hat{\phi} (p)$ extends to a global holomorphic function $\phi (p)$ on the universal cover $\tilde{\sfSigma}_{R, \sfDelta, p}$ of the punctured source fibre $\tilde{\sfSigma}_{\sfDelta, p} \smallsetminus \Gamma_p$ based at $1_p$:
\begin{eqn}
	\tilde{\sfSigma}_{R, \sfDelta,p} 
		= \widetilde{\tilde{\sfSigma}_{\sfDelta,p} \!\!\smallsetminus\! \Gamma_{p}}
		\to \tilde{\sfSigma}_{\sfDelta,p} \smallsetminus \Gamma_{p}
\fullstop
\end{eqn}
Consequently, the holomorphic germ $\hat{\phi}$ extends to a holomorphic function $\phi$ on the fundamental groupoid $\Pi_1 (\sfSigma_{R,\sfDelta})$ which is a two-dimensional holomorphic manifold.
\end{prop}

In other words, \autoref*{241019132652} says that the Borel `plane' for $\hat{f}$ and hence the formal WKB differential $\hat{\Lambda}$ is the source fibre $\tilde{\sfSigma}_{\sfDelta,p}$ of the groupoid $\Pi_1 (\sfSigma_\sfDelta)$, the set of Borel singularities is the discrete subset of critical paths $\Gamma_p \subset \tilde{\sfSigma}_{\sfDelta,p}$ starting at $p$, and the set of Borel singular values $\Xi_p \subset \CC$ is formed by the central charges of all critical paths starting at $p$; i.e., $\Xi_p = \Z (\Gamma_p)$.
Note that $\Gamma_p$ is discrete but $\Xi_p$ may have accumulation points or even be dense in some regions.

\begin{proof*}[Proof outline.]
The proof proceeds in several steps, organised into the subsections below.
Here we supply a summary outline.

We start in \autoref{240725210615} by considering the Borel-transformed global Riccati equation on the spectral curve which we attempt to solve using the method of successive approximations.
The technique is to rewrite the Borel-transformed Riccati equation as an infinite recursive system of initial value problems, each of which can be solved by flowing a canonical meromorphic vector field on the spectral curve.

In order to construct this flow, we describe in \autoref{240725160524} how this vector field may be naturally viewed as a left-invariant vector field on the fundamental groupoid $\Pi_1 (\sfSigma_{R,\sfDelta})$.
As such, we can use it propagate initial value data given along the identity bisection to the total space of the groupoid.
In \autoref{240801170203}, we put together this sequence of global initial value problem solutions on the groupoid and show that this sequence converges to a global holomorphic function on $\Pi_1 (\sfSigma_{R, \sfDelta})$ which is the analytic continuation of the Borel transform.
\end{proof*}

\begin{prop}[resurgent Stokes diagram]
\label{241019134004}	
Assume the hypothesis of \autoref{240918164353}.
For every $p \in \sfSigma_{R,\sfDelta}$, the resurgent Stokes diagram for convergent power series $\hat{\phi} (p,t) \in \CC \cbrac{t}$ coincides with the trajectory-based Stokes diagram at $p$.
Namely:
\begin{enumerate}
\item \textup{\bfseries(visible singularities).}
The set of visible singularities for $\hat{\phi} (p,t)$ is the subset $\Gamma_{p}^{0} \subset \Gamma_{p}$ of all critical trajectories starting at $p$.
\item \textup{\bfseries(Stokes rays).}
The Stokes rays for $\hat{\phi} (p,t)$ are the critical phases at $p$; the stability of a Stokes ray matches the stability of a critical phase.
\item \textup{\bfseries(regular rays).}
The regular rays for $\hat{\phi} (p,t)$ are the regular phases at $p$; the stability of a regular ray matches the stability of a regular phase.
\item \textup{\bfseries(observable singularities).}
Finally, if $\alpha$ is a Stokes ray for $\hat{\phi} (p,t)$ and $\gamma \in \Gamma_{p}^{0}$ is the corresponding critical trajectory with phase $\alpha$ (i.e., the visible singularity with phase $\alpha$), then $\hat{\phi} (p,t)$ has observable singularities with phase $\alpha$ other than $\gamma$ if and only if there is a saddle connection with phase $\alpha$ passing through the transition point where the critical trajectory $\gamma$ terminates.
\end{enumerate}
\end{prop}

\begin{proof*}
This Proposition follows immediately from \autoref{241019132652} and the geometry of the spectral curve described in detail in \autoref{240919123411}.
\end{proof*}

\begin{prop}[Borel summability of the WKB differential in a single direction]
\label{241019135458}
Assume the hypothesis of \autoref{240918164353}.
Pick a regular point $p_0 \in \sfSigma_{R, \sfDelta}$ and suppose $\alpha$ is a (semi-)stable regular ray at $p_0$.
Then the formal WKB differential $\hat{\Lambda}$ is (semi-)stably Borel summable with phase $\alpha$, uniformly for all $p$ near $p_0$.
Namely, there is a domain $U \subset \sfSigma_{R,\sfDelta}$ around $p_0$ such that the Borel resummation
\begin{eqntag}
\label{241019135750}
	\Lambda_\alpha 
		\coleq s_\alpha \big[ \, \hat{\Lambda} \, \big]
		= \Lambda_0 + f_\alpha \sigma
\qtext{where}
	f_\alpha \coleq s_\alpha \big[ \, \hat{f} \, \big]
\end{eqntag}
is well-defined and (semi-)stable uniformly on $U$.
In particular, $f_\alpha \in \cal{O} (U \times V)$ where $V \subset \CC$ is a sector at the origin with aperture $\sfop{Arc} (\alpha)$.
Moreover, as $\hbar \to 0$ in this sector, $f_\alpha$ is asymptotic to $\hat{f}$ with uniform factorial growth, uniformly on $U$:
\begin{eqntag}
\label{241019144334}
	f_\alpha (p, \hbar) \simeq \hat{f} (p, \hbar)
\qquad
	\text{as $\hbar \to 0$ unif. along $\sfop{Arc} (\alpha)$, unif. $\forall p \in U$.}	
\end{eqntag}
In fact, $f_\alpha$ is the unique holomorphic function on $U \times V$ which satisfies \eqref{241019144334}.
\end{prop}

\begin{prop}[Borel summability of the WKB differential along an arc]
\label{241019144520}
Assume the hypothesis of \autoref{240918164353}.
Pick a regular point $p_0 \in \sfSigma_{R, \sfDelta}$ and suppose $A$ is an arc of regular rays at $p_0$.
Then the formal WKB differential $\hat{\Lambda}$ is Borel summable along $A$, uniformly for all $p$ near $p_0$.
Namely, there is a domain $U \subset \sfSigma_{R,\sfDelta}$ around $p_0$ such that the Borel resummation
\begin{eqntag}
\label{241019144732}
	\Lambda_A 
		\coleq s_A \big[ \, \hat{\Lambda} \, \big]
		= \Lambda_0 + f_A \sigma
\qtext{where}
	f_A \coleq s_A \big[ \, \hat{f} \, \big]
\end{eqntag}
is well-defined uniformly on $U$.
In particular, $f_A \in \cal{O} (U \times V)$ where $V$ is a sector at the origin with aperture $\sfop{Arc} (A)$.
Moreover, as $\hbar \to 0$ in this sector, $f_A$ is asymptotic to $\hat{f}$ with factorial growth, uniformly for all $p \in U$:
\begin{eqntag}
\label{241019144926}
	f_A (p, \hbar) \simeq \hat{f} (p, \hbar)
\qquad
	\text{as $\hbar \to 0$ along $\sfop{Arc} (A)$, unif. $\forall p \in U$.}	
\end{eqntag}
In fact, $f_A$ is the unique holomorphic function on $U \times V$ which satisfies \eqref{241019144926}.
\end{prop}

\begin{prop}[lateral Borel summability of the WKB differential]
\label{241019145013}
Assume the hypothesis of \autoref{240918164353}.
Fix a point $p \in \sfSigma_{R,\sfDelta}$ and suppose $\alpha$ is a (semi-)stable Stokes ray at $p$.
Then the formal WKB differential $\hat{\Lambda}$ is laterally Borel summable at $p$ with phase $\alpha$; i.e., the left and right lateral Borel resummations
\begin{eqntag}
\label{241019145210}
	\Lambda^{\rm{L/R}}_\alpha (p) 
		\coleq s_\alpha^{\rm{L/R}} \big[ \, \Lambda (p) \, \big]
		= \Lambda_0 (p) + f^{\rm{L/R}}_\alpha (p) \sigma (p)
\fullstop{,}\quad
	f^{\rm{L/R}}_\alpha (p) \coleq s^{\rm{L/R}}_\alpha \big[ \, \hat{f}(p) \, \big]
\fullstop{,}
\end{eqntag}
where $\Lambda_0 (p), \sigma (p) \in T^\ast_p \sfSigma$, are well-defined.
In particular, $f^{\rm{L}}_\alpha (p), f^{\rm{R}}_\alpha (p) \in \cal{O} (V)$ are two holomorphic functions in a sector $V \subset \CC$ at the origin with opening $\sfop{Arc} (\alpha)$ both of which are asymptotic to $\hat{f}$ as $\hbar \to 0$ with factorial growth:
\begin{eqntag}
\label{241019145847}
	f^{\rm{L/R}}_\alpha (p, \hbar) \simeq \hat{f} (p, \hbar)
\qquad
	\text{as $\hbar \to 0$ along $\sfop{Arc} (\alpha)$.}
\end{eqntag}
\end{prop}

\begin{lem*}
\label{240918203006}
We have the following implications between the Propositions stated above:
\ref{241019125911} $\Rightarrow$ \ref{241019111435}, \ref{241019132652} $\Rightarrow$ \ref{241019113314}, \ref{241019134004} $\Rightarrow$ \ref{241019115130}, \ref{241019135458} $\Rightarrow$ \ref{240923104232}, \ref{241019144520} $\Rightarrow$ \ref{241019122801}, and \ref{241019145013} $\Rightarrow$ \ref{241022040824}.
\end{lem*}

\begin{proof*}
For all $p$ near $p_0$, the formal power series $\hat{\A}$ and $\hat{f}$ are related by the identity	
\begin{eqntag}
\label{241003114859}
	\hat{\A} (q, \hbar)
		= \exp \left( - \int_{p_0}^{p} \hat{\Lambda} (\hbar) \right)
		= \exp \left( - \int_{p_0}^{p} \Lambda_0 \right) 
			\exp \left( - \int_{p_0}^{p} \hat{f} (\hbar) \sigma \right)
\fullstop
\end{eqntag}
Consequently, the Borel transforms $\hat{\Phi}$ and $\hat{\phi}$ are related by the identity
\begin{eqntag}
\label{241003114848}
	\hat{\Phi} (q, t)
		= \exp \left( - \int_{p_0}^{p} \Lambda_0 \right) \exp^\ast \left( - \int_{p_0}^{p} \hat{\phi} (t) \sigma \right)
\fullstop{,}
\end{eqntag}
where $\exp^\ast$ is the convolution exponential:
\begin{eqn}
	\exp^\ast (z) \coleq \sum_{k=1}^\infty z^{\ast k}
\fullstop
\end{eqn}
Note that $z^{\ast 1} = z$ and $z^{\ast 0} = 0$ because we prefer to work with a non-unital convolution algebra; in particular, functions that are independent of $\hbar$ are in the kernel of the Borel transform.
Note also that the convolution product is with respect to the variable $t$.
Under the canonical identification of the neighbourhood of $t = 0$ with the neighbourhood of the constant path $1_p \in \tilde{\sfSigma}_{\sfDelta, p}$, the convolution product is taken along paths inside source fibres of the groupoid $\tilde{\sfSigma}_{R,\sfDelta}$, which are simply connected spaces.	

The Borel resummations of $\hat{\A}$ and $\hat{\Lambda}$ with phase $\alpha$ are related by the formula
\begin{eqntag}
\label{241003121439}
	\A_{\alpha} (q, \hbar)
		= s_\alpha \big[ \, \hat{\A} \, \big] (q, \hbar)
		= \exp \left( - \int_{p_0}^{p} \Lambda_\alpha (\hbar) \right)
\fullstop{,}
\end{eqntag}
provided $\hat{\Lambda}$ is uniformly Borel summable with phase $\alpha$ along the path connecting $p_0$ to $p$.
The same comment holds if $\alpha$ is replaced by an arc $A$.
\end{proof*}

\subsection{The Method of Successive Approximations}
\label{240725210615}

We now begin the construction of the endless analytic continuation of the Borel transform of the formal WKB differential in order to prove \autoref{241019132652}.

\paragraf
Recall that by \autoref{240915182325}, the formal series $\hat{f}$ is the unique $\hbar$-power series solution of the Riccati equation
\begin{eqntag}
\label{241019163611}
	\hbar \V (f) - f = \hbar (f^2 + w f + \W)
\fullstop
\end{eqntag}
and by \autoref{241019125911}, the Borel transform $\hat{\phi}$ is locally uniformly convergent.
Our goal now is to describe the analytic continuation $\phi$ of $\hat{\phi}$ by analysing the equation satisfied by $\phi$.
To this end, let $\omega$ be the Borel transform of $\W$; i.e., 
\begin{eqntag}
\label{240919102536}
	\omega (p, t) \coleq \Borel \big[~ \W ~\big] (p, t)
\qtext{so that}
	\W (p,\hbar) = \A (p, 0) + \Laplace \big[\, \omega \,\big] (p,\hbar)
\fullstop
\end{eqntag}
Note $\omega$ is a polynomial function in $t$ because $\W$ is a polynomial in $\hbar$, so the Laplace transform of $\omega$ needs no indication of direction.
Next, we consider the following integro-differential equation in $\phi$:
\begin{eqntag}
\label{240709154052}
	\left(\V - \del_t \right) \phi = \phi \ast \phi + w \phi + \omega
\fullstop
\end{eqntag}
The term $\phi \ast \phi$ is the convolution product in the variable $t$; explicitly,
\begin{eqn}
	\phi \ast \phi (p,t) \coleq \int_0^t \phi (p, t') \phi (p, t - t') \d{t'}
\fullstop
\end{eqn}
Of course, equation \eqref{240709154052} is simply the Borel transform of the Riccati equation \eqref{241019163611} (multiplied through by $\hbar^\inv$).
In particular, the Borel transform $\hat{\phi}$ of the formal power series $\hat{f}$ is the unique convergent $t$-power series solution of \eqref{240709154052}.
To be precise, let us spell out the relationship between the solutions of these equations, which is simply the content of the Borel-Laplace method phrased in our setting.

\begin{lem}
\label{240919104027}
If $\phi = \phi (p,t)$ is any holomorphic solution of \eqref{240709154052}, defined for $p$ in a neighbourhood $U \subset \sfSigma$ of $p_0$ and with a uniformly well-defined Laplace transform in the variable $t$ in a direction $\alpha$, then
\begin{eqntag}
\label{240919104203}
	f_\alpha (p, \hbar) 
		\coleq \Laplace_\alpha \big[ \, \phi \, \big] (p, \hbar)
\end{eqntag}
is a holomorphic solution of the Riccati equation \eqref{241019163611}, defined for $(p, \hbar) \in U \times V$ where $V$ is a sector in the $\hbar$-plane at the origin with aperture $\sfop{Arc} (\alpha)$, which admits a uniform asymptotic expansion
\begin{eqntag}
\label{240919105831}
	f_\alpha (p, \hbar) \sim \hat{f} (p, \hbar)
\qquad\text{as $\hbar \to 0$ along $\sfop{Arc} (\alpha)$ , unif. $\forall p$ near $p_0$}
\fullstop
\end{eqntag}
In other words, $f_\alpha \in \cal{O} (U \times V)$ is the uniform Borel resummation of $\hat{f} \in \cal{O}_\sfSigma (U) \bbrac{\hbar}$ in the direction $\alpha$:
\begin{eqntag}
\label{240919111252}
	f_\alpha (p, \hbar) = s_\alpha \big( \hat{f} \big) (p, \hbar)
\fullstop
\end{eqntag}
\end{lem}

\paragraf
Thus, the convergent Borel transform $\hat{\phi}$ gives a local solution of equation \eqref{240709154052}.
To obtain its global properties, we use a completely different method to construct a holomorphic solution $\phi$ to \eqref{240709154052} whose Taylor series expansion at $t = 0$ is the local solution $\hat{\phi}$.
Our approach is to use the method of successive approximations to rewrite \eqref{240709154052} as the following infinite recursive sequence of initial value problems:
\begin{eqntag}
\label{240712152412}
	\left(\V - \del_t \right) \phi_\pto{k} = \Phi_\pto{k}
\qqtext{and}
	\phi_\pto{k} \big|_{t = 0} = 0
\qqquad
	(k \geq 1)
\fullstop{,}
\end{eqntag}
where the recursion kernel is given by
\begin{eqntag}
\label{240923164138}
	\Phi_\pto{k} \coleq 
		\sum_{i + j = k - 2} \phi_\pto{i} \ast \phi_\pto{j}
		 + w \phi_\pto{k-1}
\qquad\text{$(k \geq 2)$}
\fullstop{,}
\end{eqntag}
and the starting data is
\begin{eqntag}
\label{240923164306}
	\phi_\pto{0} \coleq f_1 = \frac{\Lambda_1}{\sigma}
\qqtext{and}
	\Phi_\pto{1} \coleq w \phi_\pto{0} + \omega
\fullstop
\end{eqntag}
The advantage of this infinite recursive system over the single equation \eqref{240709154052} is the fact that each equation in \eqref{240712152412} can be solved by direct integration of the vector field $\V - \del_t$ appearing on the lefthand side.
Any solution of \eqref{240712152412} then yields a solution of \eqref{240709154052} as follows.

\begin{lem}
\label{240805131610}
If a solution $\set{ \phi_\pto{k} }_{k = 0}^\infty$ of \eqref{240712152412} is found with the property that the infinite series
\begin{eqntag}
\label{230517133455}
	\phi \coleq \sum_{k=0}^\infty \phi_\pto{k}
\end{eqntag}
is locally uniformly convergent, then $\phi$ is a solution of \eqref{240709154052}.
\end{lem}

The verification of this Lemma is a straightforward calculation which we omit; e.g.,see \cite[§C.3]{MY2008.06492}.
The recursive system is very easy to solve locally near $t = 0$, but it is not obvious how to construct the endless analytic continuation of this local solution.
Thus, the main difficulty in solving the recursive system \eqref{240712152412} is that we want to find an essentially global solution.
In fact, an actual global solution which is valid for all $t \in \CC$ does not exist because these equations necessarily have singularities in the $t$-plane.
The key to doing so is to understand the vector field $\V - \del_t$ from a global geometric point of view.

\subsection{Generalised Global Flow on the Spectral Cover}
\label{240725160524}

\paragraf
The vector field $\V$ on $\sfSigma$ has poles and zeros wherever $\sigma$ has zeros and poles, respectively.
Consequently, $\V$ is a global nonvanishing section of $\cal{T}_{\sfSigma} (R-\sfDelta)$.
In particular, it is holomorphic and nonvanishing on the subset of regular points $\sfSigma_{R,\sfDelta} \subset \sfSigma$ where it therefore generates a maximal complex flow
\begin{eqn}
	\varphi : \Omega \subset \sfSigma_{R,\sfDelta} \times \CC \too \sfSigma_{R,\sfDelta}
\fullstop
\end{eqn}
Recall that the flow domain $\Omega$ of a holomorphic vector field has a natural structure of a holomorphic Lie groupoid $\Omega \rightrightarrows \sfSigma_{R,\sfDelta}$ with source map $\rm{s} (p, t) = p$ and target map $\rm{t} (p, t) = \varphi (p, t)$.
If $\V$ were a complete vector field, then $\Omega$ would be nothing but the action groupoid $\sfSigma_{R,\sfDelta} \rtimes \CC$.
When $\V$ is not complete (as is the case for us), $\Omega$ is a strict subset of $\sfSigma_{R,\sfDelta} \times \CC$.
Apart from the fact that $\Omega$ contains an open neighbourhood of the identity bisection $\sfSigma_{R,\sfDelta} = \sfSigma_{R,\sfDelta} \times \set{0} \subset \Omega$, not much else can in general be concluded about $\Omega$.

\paragraf
On the other hand, consider the fundamental groupoid $\Pi_1 (\sfSigma_{R,\sfDelta}) \rightrightarrows \sfSigma_{R,\sfDelta}$ with source and target maps denoted by $\rm{s}, \rm{t}$.
The multiplicative function $\Z$ yields a holomorphic surjective map
\begin{eqn}
	\rho \coleq (\rm{s}, \Z) : \Pi_1 (\sfSigma_{R,\sfDelta}) \too \sfSigma_{R,\sfDelta} \times \CC
\end{eqn}
which we call the \dfn{flow anchor map}.
If the vector field $\V$ were complete, then $\rho$ would be a groupoid isomorphism onto the action groupoid $\sfSigma_{R,\sfDelta} \rtimes \CC$ determined by $\V$.
When $\V$ is not complete, the action groupoid $\sfSigma_{R,\sfDelta} \rtimes \CC$ does not exist, but the flow anchor map $\rho$ is still well-defined.
The main observation is that the fundamental groupoid $\Pi_1 (\sfSigma_{R,\sfDelta})$ and the flow groupoid $\Omega$ are isomorphic near the identity bisection $\sfSigma_{R,\sfDelta}$ via the flow anchor map which in particular intertwines the target map $\rm{t} : \Pi_1 (\sfSigma_{R,\sfDelta}) \to \sfSigma_{R,\sfDelta}$ and the complex flow $\varphi : \Omega \to \sfSigma_{R,\sfDelta}$.

\begin{lem}
\label{240418194933}
The flow anchor map $\rho$ defines an isomorphism of groupoid germs around the identity bisection:
\begin{eqn}
	\rho : \big[ \Pi_1 (\sfSigma_{R,\sfDelta}) \big]_{\sfSigma_{R,\sfDelta}} \iso [\Omega]_{\sfSigma_{R,\sfDelta}}
\fullstop
\end{eqn}
\end{lem}

\begin{proof*}
Choose any point $p \in \sfSigma_{R,\sfDelta}$.
We need to show that $p$, viewed as a point on the identity bisections of both $\Pi_1 (\sfSigma_{R,\sfDelta})$ and $\Omega$, has an open neighbourhood $W_1$ in $\Pi_1 (\sfSigma_{R,\sfDelta})$ and an open neighbourhood $W_2$ in $\Omega$ such that the flow anchor map restricts to a biholomorphism $W_1 \iso W_2$ that intertwines the groupoid structure maps.

Since the one-form $\sigma$ is nonvanishing in a neighbourhood of $p$, there is a simply connected neighbourhood $U \subset \sfSigma_{R,\sfDelta}$ of $p$ such that the map $\Z_p : x \mapsto t = \int_p^x \sigma \eqcol \Z_p (x)$ defines a biholomorphism from $U$ to a disc $V \subset \CC$ centred at $0$ of some radius $r>0$.
Fix a strictly smaller concentric disc $V_0 \subset V$ of radius $r_0 < r$, and let $U_0 \subset U$ be the preimage of $V_0$ under $\Z_p$.
Finally, let $\DD \subset \CC$ be a disc centred at the origin $0$ of radius $r - r_0$.
Note that for any point $t_0 \in V_0$ and any $t \in \DD$, the point $t_0 + t$ is contained in $V$.
Now we define
\begin{eqn}
	W_1 \coleq \rm{s}^\inv (U_0) \cap \Z (\DD) \big|^\text{c} \subset \Pi_1 (\sfSigma_{R,\sfDelta})
\qtext{and}
	W_2 \coleq U_0 \times \DD \subset \Omega
\fullstop{,}
\end{eqn}
where `` $\dumbdash |^\text{c}$ '' means connected component containing the constant path at $p$.
In words, $W_1$ is the set of all homotopy classes of paths $\gamma$ on $\sfSigma_{R,\sfDelta}$ that start in $U_0$ and whose $\Z$-length is less than $r - r_0$ (i.e., such that $\Z(\gamma)$ is contained in the small disc $\DD$).
Taking the connected component ensures that any such path $\gamma$ is in fact contained in $U$.

Since $U$ is simply connected, the restriction of the flow anchor map is clearly an isomorphism $\rho: W_1 \iso W_2$.
For suppose $\gamma, \gamma' \in \Pi_1 (\sfSigma_{R,\sfDelta})$ are two paths with the same source $p_0 \in U_0$ and such that $t = \Z (\gamma) = \Z (\gamma')$.
Then the fact that $t \in \DD$ implies that $\gamma$ and $\gamma'$ are homotopic and contained in $U$ because the point $t_0 + t$ where $t_0 = \Z_p (p_0)$ is still contained in $V$.
Conversely, if $(p_0, t) \in U_0 \times \DD$, then there is a unique point $q \in U$ such that $\Z_p (q) = \Z_p (p_0) + t$, so the path $\gamma$, which is contained in $U$ and goes from $p$ to $q$, is such that $\rho (\gamma) = (p_0, t)$.
\end{proof*}

\begin{rem}[generalised global flows]
\label{240725205848}
The argument in the proof of this Lemma can be expanded to include the divisors $R$ and $\sfDelta$, which then shows that the flow $\varphi$ of the vector field $\V$ actually extends to the fundamental groupoid $\Pi_1 (\sfSigma, \sfDelta) \rightrightarrows (\sfSigma, \sfDelta)$ (i.e., the \textit{twisted fundamental groupoid} in the sense of Gualtieri-Li-Pym \cite{MR3808258}), where it coincides with the groupoid target map $\rm{t} : \Pi_1 (\sfSigma, \sfDelta) \to \sfSigma$.
In this manner, by replacing the flow domain $\Omega$ with the fundamental groupoid $\Pi_1 (\sfSigma, \sfDelta)$, the flow of $\V$ can be thought of as a type of generalised morphism $\varphi : \sfSigma \times \CC \dashrightarrow \sfSigma$:
\begin{eqn}
\begin{tikzcd}
&	\Pi_1 (\sfSigma, \sfDelta)
		\ar[dl, "\rho"']
		\ar[dr, "\rm{t}"]
\\
	\sfSigma \times \CC
		\ar[rr, dashed, "\varphi"']
&&	\sfSigma
\fullstop
\end{tikzcd}
\end{eqn}
This is a powerful idea that allows one to define \textit{generalised} global flows of \textit{meromorphic} vector fields.
In this paper, we use this idea in a simplified setting where the above diagram is restricted to the regular locus $\sfSigma_{R,\sfDelta} \subset \sfSigma$.
\end{rem}

\paragraf
Let $\tilde{\V}$ be the unique left-invariant vector field on the fundamental groupoid $\Pi_1 (\sfSigma_{R,\sfDelta})$ which is the source-lift of the vector field $\V$ on $\sfSigma_{R,\sfDelta}$; i.e.,
\begin{eqn}
	\rm{s}_\ast \tilde{\V} = \V
\fullstop
\end{eqn}
Recall that this vector field is tangent to the target-fibres of $\Pi_1 (\sfSigma_{R,\sfDelta})$.
The key observation is that, under the local identification of \autoref{240418194933}, the vector field $\V - \del_t$ on $\sfSigma_{R,\sfDelta} \times \CC$ appearing on the lefthand side of the Borel-transformed Riccati equation \eqref{240709154052} is nothing but the left-invariant vector field $\tilde{\V}$ on $\Pi_1 (\sfSigma_{R,\sfDelta})$ in the following sense.

\begin{lem}
\label{240725162959}
Consider the vector field $\V - \del_t$ on $\sfSigma_{R,\sfDelta} \times \CC$.
Its germ along the identity bisection $\sfSigma_{R,\sfDelta} \cong \sfSigma_{R,\sfDelta} \times \set{0} \subset \sfSigma_{R,\sfDelta} \times \CC$ is the pushforward of the left-invariant vector field $\tilde{\V}$ on $\Pi_1 (\sfSigma_{R,\sfDelta})$ by the flow anchor map $\rho$:
\begin{eqn}
	\rho_\ast \tilde{\V} = \V - \del_t
\fullstop
\end{eqn}
\end{lem}

\begin{proof*}
The target fibre $\rm{t}^\inv (q)$ of any $q \in \sfSigma_{R,\sfDelta}$ is given in local coordinates $(x,t)$ near the identity bisection by the equation $t = \int_x^q \sigma$; i.e., $\Z_p (x) + t = 0$.
Thus, the left-invariant lift of $\V$ to $\Omega$ is precisely the vector field $\V - \del_t$.
\end{proof*}

\begin{prop}
\label{240708173915}
Let $\Phi$ and $g$ be any pair of holomorphic functions on $\Pi_1 (\sfSigma_{R,\sfDelta})$ and $\sfSigma_{R,\sfDelta}$ respectively.
Then the following initial value problem has a unique holomorphic solution $\phi$ on $\Pi_1 (\sfSigma_{R,\sfDelta})$:
\begin{eqntag}
\label{240709153501}
	\tilde{\V} \phi = \Phi
\qqtext{such that}
	\phi \big|_{\sfSigma_{\sfDelta}} = g
\end{eqntag}
In particular, for any $p \in \sfSigma_{R,\sfDelta}$, let $\hat{\phi}_p$ denote the germ of $\phi$ at $1_p$, and use the canonical identification between the germs of open neighbourhoods of $1_p$ in $\tilde{\sfSigma}_{R,\sfDelta,p}$ and in $\tilde{\sfSigma}_{\sfDelta,p}$.
Then the restriction
\begin{eqn}
	\phi_p \coleq \phi \big|_{\tilde{\sfSigma}_{R,\sfDelta,p}}
\end{eqn}
defines an endless analytic continuation of the germ $\hat{\phi}_p$ at $1_p \in \tilde{\sfSigma}_{\sfDelta,p}$ to the whole source fibre $\tilde{\sfSigma}_{\sfDelta,p}$ away from the subset of critical paths $\Gamma_p \subset \tilde{\sfSigma}_{\sfDelta,p}$ starting at $p$.
\end{prop}

\begin{proof*}
Take any path $\gamma : [0,1] \to \sfSigma_{R,\sfDelta}$ and denote its terminal point by $q \coleq \rm{t} (\gamma)$.
Let $\tau \coleq |\gamma| = \int_\gamma |\d{t}| \in \RR_+$ be its length, and let $(\gamma_r)$ be the path homotopy in the simply connected target fibre $\rm{t}^\inv (q)$ parameterised by $r \in [0,\tau]$ between the constant path $\gamma_0 = 1_q$ at $q$ and the path $\gamma_\tau = \gamma$ such that $|\gamma_r| = r$.
Then the solution $\phi$ is given by the formula
\begin{eqntag}
\label{240725173207}
	\phi (\gamma) = g(q) + \int_0^\tau \Phi (\gamma_r) \d{r}
\fullstop
\end{eqntag}
It just remains to remark that this expression only depends on the homotopy class of $\gamma$ in $\sfSigma_{R,\sfDelta}$.
\end{proof*}

\subsection{Endless Analytic Continuation of the Borel Transform}
\label{240801170203}

We can now reformulate the system of initial value problems \eqref{240712152412} as a global existence and uniqueness question on the groupoid $\Pi_1 (\sfSigma_{R,\sfDelta})$.

\paragraf
Let $\tilde{w} \coleq \rho^\ast w = \rm{s}^\ast w$ and $\tilde{\omega} \coleq \rho^\ast \omega$.
These are holomorphic functions on the fundamental groupoid $\Pi_1 (\sfSigma_{R,\sfDelta})$.
Now, consider the following infinite recursive system of initial value problems on $\Pi_1 (\sfSigma_{R,\sfDelta})$:
\begin{eqntag}
\label{240725194251}
	\tilde{\V} \phi_\pto{k} = \tilde{\Phi}_\pto{k}
\qqtext{and}
	\phi_\pto{k} \big|_{t = 0} = 0
\qqquad
	(k \geq 1)
\fullstop{,}
\end{eqntag}
where the recursion kernel is given by
\begin{eqntag}
\label{240805143427}
	\tilde{\Phi}_\pto{k} \coleq 
		\sum_{i + j = k - 2} \phi_\pto{i} \ast \phi_\pto{j}
		 + \tilde{w} \phi_\pto{k-1}
\qquad\text{$(k \geq 2)$}
\fullstop{,}
\end{eqntag}
and the starting data is
\begin{eqntag}
\label{240923164821}
	\phi_\pto{0} \coleq \rm{s}^\ast f_1 = \frac{\rm{s}^\ast \Lambda_1}{\rm{s}^\ast \sigma}
\qqtext{and}
	\tilde{\Phi}_\pto{1} \coleq \tilde{w} \phi_\pto{0} + \tilde{\omega}
\fullstop
\end{eqntag}
Note that the convolution product in \eqref{240805143427} is taken within each source fibre of $\Pi_1 (\sfSigma_{R,\sfDelta})$ which is a simply connected space with a chosen basepoint.
Namely, if $\gamma \in \sfSigma_{R,\sfDelta,p}$ and $\phi_1, \phi_2 \in \cal{O} (\sfSigma_{R,\sfDelta,p})$, then
\begin{eqntag}
\label{241011133323}
	\phi_1 \ast \phi_2 (\gamma)
		= \int_0^\tau \phi_1 (\gamma_r) \phi_2 (\gamma_{\tau-r}) \d{r}
\fullstop
\end{eqntag}

By \autoref{240708173915}, for each $k$, there is a unique holomorphic function $\phi_\pto{k}$ on $\Pi_1 (\sfSigma_{R,\sfDelta})$ that satisfies \eqref{240725194251}, and it is given by the formula
\begin{eqntag}
\label{240801201519}
	\phi_\pto{k} (\gamma) = \int_0^\tau \tilde{\Phi}_\pto{k} (\gamma_r) \d{r}
\fullstop
\end{eqntag}
We now use this formula to prove that the infinite series $\phi = \sum \phi_\pto{k}$ of \autoref{240805131610} converges to a holomorphic solution of the Borel-transformed Riccati equation \eqref{240709154052} which defines the analytic continuation of the Borel transform $\hat{\phi}$ of $\hat{f}$.

\begin{prop}
\label{240920100212}
The infinite series $\phi$ defines a holomorphic function on the fundamental groupoid $\Pi_1 (\sfSigma_{R,\sfDelta})$ whose germ along the identity bisection $\sfSigma_{R,\sfDelta}$ equals the Borel transform $\hat{\phi}$ of the divergent series $\hat{f}$.
In particular, let $W \subset \sfSigma_\sfDelta$ be any relatively compact open neighbourhood of the ramification locus $R$, and denote the complement by $\sfSigma'_{\sfDelta} \coleq \sfSigma_\sfDelta \smallsetminus W \subset \sfSigma_{\sfDelta}$.
Then there are real constants $\C, \K > 0$ such that for every $k \geq 0$ and every path $\gamma$ on $\sfSigma'_{\sfDelta}$ of length $\tau \coleq |\gamma|$, the function $\phi$ satisfies the exponential bound
\begin{eqntag}
\label{240805153231}
	\big| \phi (\gamma) \big| \leq \C e^{\K \tau}
\fullstop
\end{eqntag}
\end{prop}

All the assertions of \autoref{241019132652} now follow immediately from \autoref*{240920100212}.
The key to proving this Proposition is the following bound on each term $\phi_\pto{k}$.

\begin{lem}
\label{240801201822}
Let $W \subset \sfSigma_\sfDelta$ be any relatively compact open neighbourhood of the ramification locus $R$, and denote the complement by $\sfSigma'_{\sfDelta} \coleq \sfSigma_\sfDelta \smallsetminus W \subset \sfSigma_{\sfDelta}$.
Then there are real constants $\C, \M, \L > 0$ such that for every $k \geq 0$ and every path $\gamma$ on $\sfSigma'_{\sfDelta}$ of length $\tau \coleq |\gamma|$, the function $\phi_\pto{k}$ satisfies the bound
\begin{eqntag}
\label{240801202158}
	\big| \phi_\pto{k} (\gamma) \big| \leq \C \M^k \frac{\tau^k}{k!} e^{\L \tau}
\fullstop
\end{eqntag}
Furthermore, there is a universal constant $m > 0$ (i.e., independent of all choices including $W$) such that the constant $\M$ may be chosen to be $\M = m \C$.
\end{lem}

\begin{proof*}
First, note that there are constants $\C, \L > 0$ such that for all $(p,t) \in \sfSigma'_\sfDelta \times \CC$,
\begin{eqntag}
\label{240801204722}
	\big| \phi_\pto{0} (p) \big|,
	\big| w (p) \big| \leq \C
\qtext{and}
	\big| \omega (p,t) \big| \leq \C e^{\L |t|}
\fullstop
\end{eqntag}
It follows in particular that, for any $\gamma \in \Pi_1 (\sfSigma'_\sfDelta)$ of length $\tau = |\gamma|$,
\begin{eqntag}
\label{240805121639}
	\big| \tilde{w} (\gamma) \big| \leq \C
\qtext{and}
	\big| \tilde{\omega} (\gamma) \big| \leq \C e^{\L \tau}
\fullstop
\end{eqntag}
To demonstrate the bound \eqref{240801202158}, we first recursively construct a sequence of positive real numbers $\set{ \M_k }$ such that
\begin{eqntag}
\label{240805153106}
	\big| \phi_\pto{k} (\gamma) \big| \leq \M_k \frac{\tau^k}{k!} e^{\L \tau}
\fullstop
\end{eqntag}
Then we will show that there is a positive constant $\M$ such that $\M_k \leq \C \M^k$.
In fact, we will show that $\M = m \C$ for some universal constant $m > 0$ which in particular does not depend on $\C$.

First, we can take $\M_0 = \C$.
Then the constants $\M_k$ can then be constructed by induction on $k$ with the help of some elementary integral estimates (see, e.g., \cite[§C.3]{MY2008.06492}).
Specifically, the convolution term can be estimated as follows:
\begin{eqntag}
\label{240923155047}
	\big| \phi_\pto{i} \ast \phi_\pto{j} (\gamma) \big|
		\leq \M_i \M_j \frac{\tau^{i+j+1}}{(i+j+1)!} e^{\L \tau}
\fullstop
\end{eqntag}
Then estimating in \eqref{240805143427} yields
\begin{eqn}
	\big| \tilde{\Phi}_\pto{k} (\gamma) \big|
		\leq \M_{k} \frac{\tau^{k-1}}{(k-1)!} e^{\L \tau}
\qtext{where}
	\M_{k} \coleq \C \left( \sum_{i + j = k - 2} \M_i \M_j + \M_{k-1} \right)
\fullstop
\end{eqn}
Now \eqref{240805153106} follows from \eqref{240801201519}.

Let us now introduce a new sequence of real numbers $\set{m_k}$ defined by the relation
\begin{eqntag}
\label{240923184739}
	\M_k = \C \C^k m_k
\qquad\text{$(k \geq 0)$}.
\end{eqntag}
They satisfy the following recursive relation
\begin{eqntag}
\label{240923185237}
	m_0 = 1
\fullstop{,}
\qqquad m_1 = 1
\fullstop{,}
\qqquad m_k = \sum_{i+j=k-2} m_i m_j + m_{k-1}
\qquad\text{$(k \geq 2)$}
\fullstop
\end{eqntag}
The first few terms of this sequence are
\begin{eqn}
	1, 1, 2, 4, 9, 21, 51, 127, 323, 835, \ldots
\end{eqn}
Observe that this sequence is independent of all choices.

To show that there is a positive constant $m > 0$ such that $m_k \leq m^k$, we argue as follows.
Consider the following power series in an abstract variable $z$:
\begin{eqn}
	\hat{g} (z) \coleq \sum_{k=0}^\infty m_k z^k \in \CC \bbrac{z}
\fullstop
\end{eqn}
We will show that $\hat{g} (z)$ is in fact a convergent power series.
First, we observe that $\hat{g} (0) = m_0 = 1$ and that $\hat{g} (z)$ satisfies the following algebraic equation:
\begin{eqntag}
\label{240805150108}
	\hat{g} = 1 + z^2 \hat{g}^2 + z \hat{g}
\fullstop
\end{eqntag}
This identity can be verified by expanding and comparing the coefficients using the recursive formula \eqref{240923185237} for $m_k$.
Now, consider the holomorphic function $\G = \G (g, z)$ of two variables defined by
\begin{eqn}
	\G (g, z) \coleq - g + 1 + z^2 g^2 + zg
\fullstop	
\end{eqn}
It has the following properties:
\begin{eqn}
	\G (1, 0) = 0
\qqtext{and}
	\frac{\del \G}{\del g} \bigg|_{(g,z) = (1, 0)} = -1 \neq 0
\fullstop
\end{eqn}
By the Holomorphic Implicit Function Theorem, there exists a unique function $g(z)$ which is holomorphic at $z = 0$ and which satisfies $g(0) = 1$ and $\G \big( g(z), z \big) = 0$ for all $z$ sufficiently close to $0$.
Since $\hat{g} (0) = 1$ and $\G \big( \hat{g}(z), z \big) = 0$, the power series $\hat{g} (z)$ must be the Taylor expansion of $g(z)$ at $z = 0$.
As a result, $\hat{g} (z)$ is in fact a convergent power series, which means that its coefficients grow at most exponentially with $k$: i.e., there is a constant $m > 0$ such that $m_k \leq m^k$ for all $k$.
Combining this inequality with \eqref{240923184739}, we arrive at the desired result.
\end{proof*}

\begin{rem}[Motzkin numbers]
\label{240923185923}
In the final stages of preparation of this manuscript, out of curiosity we fed the sequence of numbers $\set{m_k}$ given by \eqref{240923185237} into the Online Encyclopedia of Integer Sequences to discover that they are called the \textit{Motzkin numbers}, which are closely related to Catalan numbers.
Amongst many other applications, they enumerate the different ways of drawing non-intersecting chords between points on a circle; e.g., see \cite{MR1691863}.
\end{rem}

\begin{proof}[Proof of \autoref{240920100212}.]
The bounds \eqref{240801202158} imply the exponential bound \eqref{240805153231} as follows:
\begin{eqntag}
\label{240801203211}
	\big| \phi (\gamma) \big| 
		\leq \sum_{k=0}^\infty \big| \phi_\pto{k} (\gamma) \big|
		\leq \sum_{k=0}^\infty \C \M^k \frac{\tau^k}{k!} e^{\L \tau}
		\leq \C e^{(\M + \L) \tau}
\fullstop
\end{eqntag}
Thus, we can set $\K \coleq \M + \L$.
Moreover, this bound holds uniformly for all $\gamma \in \Pi_1 (\sfSigma'_{R,\sfDelta})$ of length at most $\tau$, so the infinite series $\phi$ converges locally uniformly on $\Pi_1 (\sfSigma'_{R,\sfDelta})$ and therefore defines a holomorphic function on $\Pi_1 (\sfSigma'_{R,\sfDelta})$.
By taking smaller and smaller relatively compact neighbourhoods $W$ of $R$, we can exhaust $\Pi_1 (\sfSigma_{R,\sfDelta})$ by subsets of the form $\Pi_1 (\sfSigma'_{R,\sfDelta})$.
In conclusion, the infinite series $\phi$ converges locally uniformly on the groupoid $\Pi_1 (\sfSigma_{R,\sfDelta})$ and therefore defines a holomorphic function.

Furthermore, the germ of $\phi$ along the identity bisection $\sfSigma_{R,\sfDelta}$ is canonically a convergent power series in $t$ which by \autoref{240805131610} satisfies the Borel-transformed Riccati equation \eqref{240709154052}.
But the Borel transform $\hat{\phi}$ of the formal power series $\hat{f}$, as given by \eqref{240918203538}, is the unique convergent $t$-power series solution of \eqref{240709154052}.
Therefore, the germ of $\phi$ along the identity bisection must coincide with the Borel transform $\hat{\phi}$.
\end{proof}

\subsection{Exponential Type and Borel Resummation}
\label{241014142237}

The following Proposition is just a rephrasing of Propositions \ref{241019135458}, \ref{241019144520}, and \ref{241019145013}.

\begin{prop}
\label{240922111706}
Let $p \in \sfSigma$ be any regular point.
If $\alpha$ is a (semi-)stable regular ray at $p$, then the germ $\hat{\phi} (p, t)$ has a well-defined Laplace transform with phase $\alpha$.
If $\alpha$ is a (semi-)stable Stokes ray at $p$, then the germ $\hat{\phi} (p, t)$ has well-defined lateral Laplace transforms with phase $\alpha$.
\end{prop}

\begin{proof*}
Suppose first that $\alpha$ is regular.
By \autoref{240924153408}, the geodesic $\gamma$ starting at $p$ with phase $\alpha$ has an open neighbourhood $U$ which is mapped by the central charge $\Z_p$ to a halfstrip around the halfline $\RR_\alpha$.
In particular, every path in $U$ starting at $p$ is an element of the fundamental groupoid $\Pi_1 (\sfSigma_{R, \sfDelta})$.
At the same time, the complement $W \coleq \sfSigma_\sfDelta \smallsetminus \bar{U}$ is a relatively compact neighbourhood of the ramification locus, so by \autoref{240801201822} there are constants $\C, \K > 0$ such that, for every trajectory $\beta$ contained in $U$ of length $\tau = |\beta|$, the holomorphic function $\phi$ on the groupoid $\Pi_1 (\sfSigma_{R, \sfDelta})$ satisfies the exponential bound $|\phi (\beta)| \leq \C e^{\K \tau}$.
In particular, if we let $\beta_t$ be the canonical lift via $\Z_p$ of the line segment $[0,t e^{i\alpha}]$ to a geodesic starting at $p$, then $|\phi (\beta_t)| \leq \C e^{\K t}$.
The case of a Stokes ray $\alpha$ is argued similarly.
\end{proof*}

\begin{rem*}
\label{241019150107}
This completes the proof of all assertions in Propositions \ref*{241019125911}-\ref*{241019145013}, and hence Propositions \ref*{241019111435}-\ref*{241022040824}.
Therefore, the proofs of \autoref{240918164353} as well as this paper's main \autoref{240801161941} are now complete.
\end{rem*}

\begin{adjustwidth}{-2cm}{-1.5cm}
{\footnotesize
\bibliographystyle{nikolaev}
\bibliography{/Users/nikita/Documents/References}
}
\end{adjustwidth}
\end{document}